\documentclass[12pt]{amsart}

\usepackage{amssymb}
\usepackage{amsmath}
\usepackage{amsxtra}
\usepackage{amscd}
\usepackage{graphicx}
\usepackage{amsfonts}
\usepackage{pb-diagram}
\usepackage{float}

\numberwithin{equation}{section}

\newtheorem{thm}{Theorem}
\newtheorem{lem}{Lemma}
\newtheorem{cor}{Corollary}
\newtheorem{prop}{Proposition}

\newtheorem{defn}{Definition}

\newtheorem{rem}{Remark}
\newtheorem{note}{Note}

\newtheorem*{notations}{Notations}
\newtheorem*{notation}{Notation}
\begin{document}

\title{$p$--Adic framed braids II }

\author{J. Juyumaya}
\address{J. Juyumaya: Universidad de Valpara\'{\i}so \\
Gran Breta\~na 1091, Valpara\'{\i}so, Chile.}
\email{juyumaya@uvach.cl}

\author[S. Lambropoulou]{S. Lambropoulou\\with an appendix by \\ P.
 G\'{e}rardin}
\address{S. Lambropoulou: Department of Mathematics,
National Technical University of Athens,
Zografou campus, GR--157 80 Athens, Greece.}
\email{sofia@math.ntua.gr}
\urladdr{http://www.math.ntua.gr/$\tilde{~}$sofia}


\address{P. G\'{e}rardin: Case 7012 \\
2 place Jussieu, F--75251 Paris cedex 05.}
\email{gerardin@math.jussieu.fr}

\thanks{This research  has been co-financed by the European Union (European Social Fund - ESF) and Greek national funds through the Operational Program ``Education and Lifelong Learning" of the National Strategic Reference Framework (NSRF) - Research Funding Program: THALES: Reinforcement of the interdisciplinary and/or inter-institutional research and innovation. Moreover, both authors were partially supported  by Fondecyt 1085002, the National Technical U. Athens and Dipuv. }



\begin{abstract}
The Yokonuma--Hecke algebras are quotients of the modular framed braid group and they support Markov traces.  In this paper, which is sequel to \cite{jula1}, we explore further the structures of the $p$--adic framed braids and the $p$--adic Yokonuma--Hecke algebras constructed in \cite{jula1}, by means of dense sub--structures approximating $p$--adic elements.
We also construct a $p$--adic Markov trace on the $p$--adic Yokonuma--Hecke algebras  and we approximate the values of the $p$--adic trace on $p$--adic elements. Surprisingly, the Markov traces do not re--scale directly to yield isotopy invariants of framed links. This leads to imposing the `$E$--condition' on the trace parameters.  For solutions of the `$E$--system' we then define 2--variable isotopy invariants of modular framed links. These lift to $p$--adic isotopy invariants of classical framed links. The Yokonuma--Hecke algebras have topological interpretations in the context of framed knots, of classical knots  of singular knots and of transverse knots.
 \end{abstract}

\date{}
\maketitle

\begin{center}
\begin{minipage}{10cm}
{\scriptsize
\tableofcontents}
\end{minipage}
\end{center}

\vspace{3mm}

\begin{center}
{\sc Introduction}
\end{center}

The Yokonuma-Hecke  algebras ${\rm Y}_{d, n}(u)$ were introduced by Yokonuma \cite{yo} in the context of  Chevalley groups, as   generalizations  of the Iwahori-Hecke algebras. More precisely, the Iwahori-Hecke algebra associated to a finite Chevalley group $G$ is the centralizer algebra associated to the permutation representation of $G$ with respect to a Borel subgroup of $G$. The Yokonuma-Hecke  algebra is the centralizer algebra associated to the permutation representation of $G$ with respect to  a unipotent subgroup of $G$. Setting $d=1$ the algebra ${\rm Y}_{1, n}(u)$ coincides with the classical Iwahori--Hecke algebra of type $A$, ${\rm H}_n(u)$.

The algebra ${\rm Y}_{d, n}(u)$ may be viewed \cite{ju} as a quotient of the modular framed braid group ${\mathcal F}_{d,n}$ (where the framings are considered modulo $d$) over a quadratic relation (Eq.~\ref{quadr}) involving the framing generators $t_i$ in a subtle way, by means of weighted idempotents $e_{d,i} \in {\Bbb C}{\mathcal F}_{d,n}$ (Eq.~\ref{edi}). That is:
${\rm Y}_{d, n}(u) = \frac{{\Bbb C}{\mathcal F}_{d,n}}{ I_{d,n}}$,
where $I_{d,n}$ is the ideal generated by the expressions: $\sigma_i^2 - 1 - (u-1) \, e_{d,i} - (u-1) \, e_{d,i} \, \sigma_i$.
See Figures~\ref{ed1} and \ref{g1invrs} for some diagrammatic interpretations and also \cite{jula1} for more.

In \cite{ju} the first author constructed a linear Markov trace ${\rm tr}_d$ on the algebras ${\rm Y}_{d, n}(u)$. With the trace ${\rm tr}_d$ in hand we assumed we had an invariant for framed links, using the framed braid equivalence, derived in an analogous manner as the 2-variable Jones or HOMFLYPT polynomial via the Ocneanu trace on the Iwahori-Hecke algebras  of type $A$ \cite{jo}.
 It is well--known that framed links can be used for constructing  $3$--manifolds (closed, connected, oriented) using topological surgery. Then two $3$--manifolds are homeomorphic if and only if any two framed links in $S^3$ representing them are related through isotopy and the Kirby moves, or the equivalent Fenn--Rourke moves. Isotopy invariants of framed links  have been used  for constructing topological invariants of $3$--manifolds.
 One motivation for our work was paper \cite{KS} where framed braid equivalence is given for the Fenn--Rourke moves.

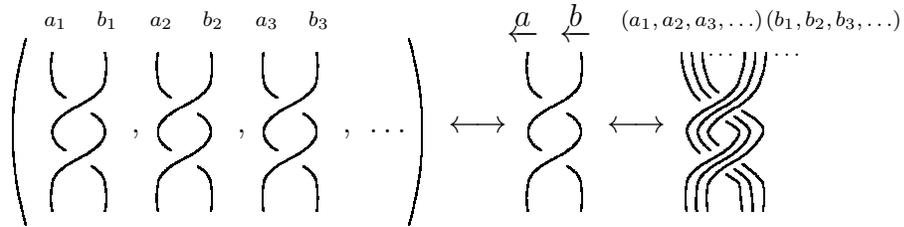
\begin{figure}[H]
\begin{picture}(320,95)

\qbezier(0,15)(-10,50)(0,85) 
\qbezier(10,70)(9,75)(10,80)
\qbezier(30,70)(31,75)(30,80)
\qbezier(10,20)(9,25)(10,30)
\qbezier(30,20)(31,25)(30,30)

\qbezier(10,30)(10,35)(20,40)
\qbezier(20,40)(30,45)(30,50)
\qbezier(30,30)(30,34)(25,37)
\qbezier(15,43)(10,46)(10,50)
\qbezier(10,50)(10,55)(20,60)
\qbezier(20,60)(30,65)(30,70)
\qbezier(30,50)(30,54)(25,57)
\qbezier(15,63)(10,66)(10,70)

\put(40,50){,}

\qbezier(50,70)(49,75)(50,80)
\qbezier(70,70)(71,75)(70,80)
\qbezier(50,20)(49,25)(50,30)
\qbezier(70,20)(71,25)(70,30)

\qbezier(50,30)(50,35)(60,40)
\qbezier(60,40)(70,45)(70,50)
\qbezier(70,30)(70,34)(65,37)
\qbezier(55,43)(50,46)(50,50)
\qbezier(50,50)(50,55)(60,60)
\qbezier(60,60)(70,65)(70,70)
\qbezier(70,50)(70,54)(65,57)
\qbezier(55,63)(50,66)(50,70)

\put(80,50){,}

\qbezier(90,70)(89,75)(90,80)
\qbezier(110,70)(111,75)(110,80)
\qbezier(90,20)(89,25)(90,30)
\qbezier(110,20)(111,25)(110,30)

\qbezier(90,30)(90,35)(100,40)
\qbezier(100,40)(110,45)(110,50)
\qbezier(110,30)(110,34)(105,37)
\qbezier(95,43)(90,46)(90,50)
\qbezier(90,50)(90,55)(100,60)
\qbezier(100,60)(110,65)(110,70)
\qbezier(110,50)(110,54)(105,57)
\qbezier(95,63)(90,66)(90,70)

\put(120,50){,}
\put(130, 50){$\ldots$}

\qbezier(145,15)(155,50)(145,85) 

\put(160,50){$\longleftrightarrow$}
\qbezier(190,70)(189,75)(190,80)
\qbezier(210,70)(211,75)(210,80)
\qbezier(190,20)(189,25)(190,30)
\qbezier(210,20)(211,25)(210,30)

\qbezier(190,30)(190,35)(200,40)
\qbezier(200,40)(210,45)(210,50)
\qbezier(210,30)(210,34)(205,37)
\qbezier(195,43)(190,46)(190,50)
\qbezier(190,50)(190,55)(200,60)
\qbezier(200,60)(210,65)(210,70)
\qbezier(210,50)(210,54)(205,57)
\qbezier(195,63)(190,66)(190,70)

\put(220,50){$\longleftrightarrow$}

\qbezier(250,50)(250,55)(260,60)
\qbezier(260,60)(269,65)(270,70)
\qbezier(254,50)(255,55)(265,60)
\qbezier(265,60)(273,65)(274,70)
\qbezier(258,50)(260,55)(270,60)
\qbezier(270,60)(277,65)(278,70)

\qbezier(255,62)(251.5,65)(250,70)
\qbezier(259,63)(255,66)(254,70)
\qbezier(262,65)(259,67)(258,70)

\qbezier(270,50)(270,51)(265,54)
\qbezier(274,50)(274,52)(268,56)
\qbezier(279,50)(280,53)(271,58)


\qbezier(250,70)(248,75)(248,80) 
\qbezier(254,70)(252,75)(252,80)
\qbezier(258,70)(256,75)(256,80)

\qbezier(270,70)(272,75)(272,80)
\qbezier(274,70)(276,75)(276,80)
\qbezier(278,70)(280,75)(280,80)


\qbezier(250,30)(250,35)(260,40)
\qbezier(260,40)(269,45)(270,50)
\qbezier(254,30)(255,35)(265,40)
\qbezier(265,40)(273,45)(274,50)
\qbezier(258,30)(260,35)(270,40)
\qbezier(270,40)(277,45)(279,50)

\qbezier(255,42)(251.5,45)(250,50)
\qbezier(259,43)(255,46)(254,50)
\qbezier(262,45)(259,47)(258,50)

\qbezier(270,30)(270,31)(265,34)
\qbezier(274,30)(274,32)(268,36)
\qbezier(279,30)(279,32)(271,38)


\qbezier(250,20)(249,25)(250,30)
\qbezier(254,20)(253,20)(254,30)
\qbezier(258,20)(257,25)(258,30)

\qbezier(271,20)(270.5,25)(270,30)
\qbezier(275,20)(275,25)(274,30)
\qbezier(279,20)(279.5,27)(279,30)

\put(7,90){{\tiny $a_1$}}
\put(27,90){{\tiny $b_1$}}

\put(47,90){{\tiny $a_2$}}
\put(67,90){{\tiny $b_2$}}

\put(87,90){{\tiny $a_3$}}
\put(107,90){{\tiny $b_3$}}


\put(182,90){$\underleftarrow{a}$}
\put(202,90){$\underleftarrow{b}$}


\put(225,90){{\tiny $(a_1, a_2, a_3, \ldots )$}}
\put(280,90){{\tiny $(b_1, b_2, b_3, \ldots )$}}

\put(258,79){{\tiny $\ldots $}}
\put(283,79){{\tiny $\ldots $}}


\end{picture}
\caption{The different facets of a $p$--adic framed braid}\label{facets}
\end{figure}

As the framings in ${\rm Y}_{d, n}(u)$ are subjected to the modular condition, we first sought  ways to improve on this point, given that in the main Kirby move the framings add up. On the other hand, the definition of the elements $e_{d,i}$, which appear in the quadratic relation of the algebra ${\rm Y}_{d, n}(u)$, rests on the fact that the framings are modular. Consequently, in \cite{jula1} we constructed the {\it $p$--adic framed braid group} ${\mathcal F}_{p^{\infty} ,n}$ as the inverse limit of the modular framed braid groups ${\mathcal F}_{p^r ,n}$:
\begin{center}
$
{\mathcal F}_{p^{\infty} ,n}  := \varprojlim_{r\in{\Bbb N}} {\mathcal F}_{p^r ,n}
$
\end{center}
In ${\mathcal F}_{p^{\infty},n}$ there are no modular relations for the framings.
A $p$--adic framed braid may be viewed as an infinite sequence of the same classical braid, with the modular framings of the corresponding strands forming a $p$--adic integer.  By certain group isomorphisms (see Eq. \ref{isoms}), a $p$--adic framed braid may also be viewed as a classical framed braid but with framings $p$--adic integers or, equivalently, as a classical braid with infinite cablings replacing each strand, whose corresponding modular framings form a $p$--adic integer. See Figure~\ref{facets} for the different facets of a $p$--adic framed braid. In the present paper we shall ignore the above identifications, which had to be observed closely in \cite{jula1}.

In \cite{jula1} we further constructed the {\it $p$--adic Yokonuma--Hecke algebra} ${\rm Y}_{p^{\infty}, n}(u)$ as the inverse limit of the algebras ${\rm Y}_{p^{r}, n}(u)$:
\begin{center}
$
{\rm Y}_{p^{\infty}, n}(u) := \varprojlim_{r\in{\Bbb N}} {\rm Y}_{p^r,n}(u)
$
\end{center}
In ${\rm Y}_{p^{\infty},n}(u)$ also there are no modular relations for the framings.

\smallbreak
In the course of our study  several interesting questions arose. One question was how  the algebras
\begin{center}
$
{\Bbb C}{\mathcal F}_{p^{\infty} ,n}  \qquad \mbox{ and } \qquad  \varprojlim_{r} {\Bbb C}{\mathcal F}_{p^r,n}
$
\end{center}
compare. In this paper, which is sequel to \cite{jula1},  we show (Proposition~\ref{algebras} and Theorem~\ref{densesets}) that the algebra ${\Bbb C}{\mathcal F}_{p^{\infty} ,n}$ can be regarded as a proper subalgebra of  $\varprojlim_{r} {\Bbb C}{\mathcal F}_{p^r,n}$ and it is, in fact, dense.  Another interesting question is how the following algebras compare:
 $$
 \varprojlim_{r} \frac{{\Bbb C}{\mathcal F}_{p^r,n}}{ I_{p^r,n}} = {\rm Y}_{p^{\infty}, n}(u) \qquad \mbox{ and } \qquad \frac{\varprojlim_{r} {\Bbb C}{\mathcal F}_{p^r,n}}{ \varprojlim_{r}  I_{p^r,n}}
 $$
In \cite{jula1} we showed that $\frac{\varprojlim_{r} {\Bbb C}{\mathcal F}_{p^r,n}}{ \varprojlim_{r}  I_{p^r,n}}$ is dense in ${\rm Y}_{p^{\infty}, n}(u)$. In this paper we show that the two algebras are in fact isomorphic (see Proposition~\ref{rhosurj}, Subsection~1.7). 
In Section~1 we recall briefly our constructions in \cite{jula1} of the $p$--adic framed braids and the $p$--adic Yokonuma--Hecke algebras, giving emphasis to the relations and the main properties in each structure. Results stated as lemmas or propositions are not contained in \cite{jula1}. For details and proofs of previous results we refer the reader to \cite{jula1}.

\smallbreak
The most important questions are about the existence of dense substructures in order to approximate $p$--adic elements by known objects. In Theorem~1\cite{jula1} we showed that inside ${\mathcal F}_{p^{\infty},n}$ lies a dense copy of the classical framed braid group ${\mathcal F}_n$. In Theorem 3\cite{jula1} we also gave a set of topological generators of ${\rm Y}_{p^{\infty},n}(u)$ and a list of relations that they satisfy. We further showed that a quadratic relation similar to the one of the algebra ${\rm Y}_{d, n}(u)$ holds in  ${\rm Y}_{p^{\infty},n}(u)$ (see Eq.~\ref{pquadr}), involving the elements $e_{i}$ in $\varprojlim_r {\Bbb C}{\mathcal F}_{p^r,n}$ and  ${\rm Y}_{p^{\infty},n}(u)$, which are lifts of the idempotents $e_{p^r,i}$  (Eq.~\ref{ei}).
 The elements $e_{i}$  are still idempotents  but they are no more weighted sums and, as we show here, in Section~2, the elements $e_i$ are purely $p$--adic. For this reason we did not manage in \cite{jula1} to figure out a dense subalgebra of ${\rm Y}_{p^{\infty},n}(u)$ given by generators and relations.

In the present paper we construct a dense subalgebra ${\widetilde{{\rm Y}_n}(u)}$ of ${\rm Y}_{p^{\infty},n}(u)$ by means of a presentation (see Definition~\ref{densextn} and Theorem~\ref{denseyhtilde}):
\[
\widetilde{{\rm Y}_n}(u) =\frac{\widetilde{{\Bbb C}{\mathcal F}_n}}
{\langle g_i^2 - 1 - (u-1)e_{i} - (u-1)e_{i}\, g_i \rangle}
\]
where $\widetilde{{\Bbb C}{\mathcal F}_n}$ is the extension of the algebra ${\Bbb C}{\mathcal F}_n$ by the elements $e_i$.

With the dense subalgebras available, we then give approximations for purely $p$--adic elements by sequences of constant elements in the dense subalgebras, after having defined the notions of {\it $s$--truncation} and {\it $r$--expansion} ($r\geq s$) for any $p$--adic element (see Subsections~2.5 and~2.6). In  particular, we show that the elements $e_{i}$ can be approximated by elements ${\bf e}_{p^r,i}$ in ${\Bbb C}{\mathcal F}_n$. We also construct  approximations for arbitrary $p$--adic elements (see Theorem~\ref{pprox}). Approximating $p$--adic elements is particularly useful for understanding deeper the $p$--adic structures. In \cite{jula1} we also discussed approximations but here we really expand on the subject.

\smallbreak
Further, in Section~3 we construct a $p$--adic Markov trace $\tau_{p^{\infty}}$ on the algebras ${\rm Y}_{p^{\infty}, n}(u)$ (see Theorem~\ref{ptrace}). This trace arises as the inverse limit of the traces ${\rm tr}_{p^r}$ (see Theorem~\ref{trace}) and it takes values in the inverse limit of certain polynomial rings.  Since ${\mathcal F}_{n}$ is dense in ${\mathcal F}_{p^{\infty},n}$,  $\tau_{p^{\infty}}$ can be used, in principle, for constructing an invariant of framed links with no restriction on the framings. Using the approximation results of Section~2 we show how to approximate the trace of a $p$--adic element by constant sequences of polynomials, and we give some computations (see Subsection~3.3).

\smallbreak
The next consideration then is to try to normalize and re--scale the traces ${\rm tr}_d$ and  $\tau_{p^{\infty}}$  according to the  framed braid equivalence, which  comprises conjugation in the classical framed braid groups ${\mathcal F}_{n}$ and positive and negative stabilization moves, in order to obtain isotopy invariants of oriented framed links.
 As we show in Section~4, trying to re--scale the trace ${\rm tr}_d$ so that $\alpha g_n$ and $\alpha g_n^{-1}$ will be assigned the same trace value for any $\alpha \in {\rm Y}_{d, n}(u)$, leads to imposing the {\it $E$--condition} on the trace parameters $x_1, \ldots ,x_{d-1}$ (Definition~\ref{defcon}). That is, $x_1, \ldots ,x_{d-1}$ have to satisfy a non linear system of equations, the {\it $E$--system} (see Eq.~\ref{Esystem}). The traces ${\rm tr}_d$ are the only known Markov traces which do not re--scale directly to yield isotopy invariants of knots.

 Surprisingly, there are always non--trivial solutions of the $E$--system in the set of complex numbers. In Subsection~4.3 we explore various solutions with emphasis to ones that do not imply the somewhat trivial condition
\begin{center}
$
{\rm tr}_d (e_{d,i}) =1.
$
\end{center}
 In the Appendix Paul G\'{e}rardin gives a method for finding the complete set of solutions of the $E$--system and shows that the solutions are parametrized by the non--empty subsets of ${\Bbb Z}/d{\Bbb Z}$. Finally, in Subsection~4.4. we show that a solution of the $E$--system lifts to the $p$--adic level.

\smallbreak
Given now a solution of the $E$--system, parametrized by a subset $S$ of ${\Bbb Z}/d{\Bbb Z}$, we re--scale and normalize the traces ${\rm tr}_d$ and we define 2--variable isotopy invariants $\Gamma_{d,S}$  (on variables: $u$ and the trace parameter $z$) of oriented framed links (see Definition~\ref{gamma} and Theorem~\ref{gammainv} in Section~5). In Proposition~\ref{skein} we also give a skein relation satisfied by the invariants $\Gamma_{d,S}$.

\begin{figure}[H]
\setlength{\unitlength}{0.7pt}
\begin{picture}(250,95)

\qbezier(0,70)(-1,75)(0,80)
\qbezier(20,70)(21,75)(20,80)
\qbezier(0,20)(-1,25)(0,30)
\qbezier(20,20)(21,25)(20,30)

\qbezier(0,30)(0,35)(10,40)
\qbezier(10,40)(20,45)(20,50)
\qbezier(20,30)(20,34)(15,37)
\qbezier(5,43)(0,46)(0,50)
\qbezier(0,50)(0,55)(10,60)
\qbezier(10,60)(20,65)(20,70)
\qbezier(20,50)(20,54)(15,57)
\qbezier(5,63)(0,66)(0,70)

\put(60,50){$\underrightarrow{\rm closure}$}

\qbezier(165.29,85.838)(174.443,91.514)(184.986,89.311)
\qbezier(152.123,70.338)(155.486,80.569)(164.979,85.658)
\qbezier(152,50)(148,59.999)(152.001,70)
\qbezier(164.979,34.342)(155.487,39.43)(152.124,49.663)
\qbezier(184.986,30.69)(174.444,29.514)(165.29,34.163)
\qbezier(202.66,40.751)(196,32.288)(185.341,30.751) 
\qbezier(209.731,59.82)(210.07,49.056)(202.892,41.027)
\qbezier(202.891,78.973)(210.07,70.945)(209.732,60.18)

\qbezier(224.71,85.838)(215.557,91.514)(205.014,89.311)
\qbezier(237.877,70.338)(234.514,80.569)(225.021,85.658)
\qbezier(238,50)(242,59.999)(237.999,70)
\qbezier(225.021,34.342)(234.513,39.43)(237.876,49.663)
\qbezier(205.014,30.69)(215.556,29.514)(224.71,34.163)

\qbezier(180.269,59.82)(179.93,49.056)(187.108,41.027)
\qbezier(187.109,78.973)(179.93,70.945)(180.268,60.18)
\qbezier(204.66,89.248)(194,87.713)(187.34,79.249) 

\put(-8,90){$\underleftarrow{a}$}
\put(13,90){$\underleftarrow{b}$}

\put(131,70){$\underleftarrow{a}$}
\put(245,70){$\underleftarrow{b}$}

\put(176.5,55){$\spcheck$}
\put(205.5,55){$\spcheck$}

\end{picture}
\caption{ A $p$--adic framed braid and a $p$--adic framed link}\label{plink}
\end{figure}
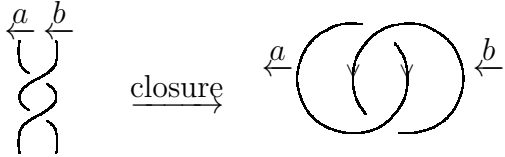

Theoretically, closing a $p$--adic framed braid gives rise to an oriented {\it $p$--adic framed link}, see Figure~\ref{plink}. So, normalizing and re--scaling the $p$--adic trace $\tau_{p^{\infty}}$ according to the $p$--adic framed braid equivalence, one could obtain isotopy invariants of $p$--adic  framed links. As we show (Theorem~\ref{padicinv}), the invariants $\Gamma_{p^r,S}$ lift to a $p$--adic invariant $\Gamma_{p^{\infty},S}$ of $p$--adic framed links. The important thing is that $\Gamma_{p^{\infty},S}$ is also an invariant of classical framed links {\it with no modular restriction on the framings}.  All the above are to be found in Section~5, where we also give computations on concrete examples.

\smallbreak

In a further development, in \cite{jula3} we represented the classical braid group into the Yokonuma--Hecke algebra ${\rm Y}_{d, n}(u)$ by treating the modular generators $t_j$ as formal generators and ignoring their topological interpretation as framing generators. This led, up to the $E$--condition and using the classical Markov braid equivalence, to $S$--parametrized 2--variable polynomial invariants of classical oriented links, which satisfy a `closed' cubic relation (closed in the sense of involving only the braiding generators). In \cite{ChLa} we show that these invariants do not coincide with the HOMFLYPT polynomial that is constructed from the Iwahori-Hecke algebras ${\rm H}_n(u)$, except in a few trivial cases, that is, $u=1$ or $q=1$ or ${\rm tr}(e_i)=1$. Yet, our computational data  \cite{CJKL} seem to indicate that these invariants do not distinguish more or less knot pairs than the HOMFLYPT polynomial.

Further, in \cite{jula4} we constructed a monoid representation of the singular braid monoid algebra into ${\rm Y}_{d, n}(u)$.  Then, given the $E$--condition and the Markov braid equivalence for singular braids, we defined invariants of oriented singular links. Finally, it has been observed (by S.~Chmutov and also by V.~Turaev) that the  traces ${\rm tr}_d$ on the Yokonuma--Hecke algebras adapt more naturally to the theory of transverse knots. Indeed, the transverse braid equivalence comprises framed braid conjugation and only positive stabilization moves. This means that the traces ${\rm tr}_d$ do not need to re--scale, so we obtain naturally topological invariants of transverse knots via algebraic means. Yet, here again, computations seem to indicate that these invariants are not stronger than known geometrical ones (see \cite{CJKL}).

To summarize, the Yokonuma--Hecke algebras comprise the only known examples of algebras  with topological applications to different knot categories, which support Markov traces that do not re--scale directly according to topological braid equivalence. For these reasons, in \cite{jula5} the framization of the BMW algebra and of other knot algebras has been established, while in \cite{GJL} the Yokonuma--Temperley--Lieb algebra has been constructed as a quotient of the Yokonuma--Hecke algebra over linear relations identical to the linear relations used for obtaining the classical Temperley--Lieb algebra as a quotient of the classical Iwahori--Hecke algebra. In \cite{GJL} there is another thing to report: as we show, the traces ${\rm tr}_d$ pass through to the quotient algebras only for the specific values of the parameter $z$ that V.F.R. Jones has computed in the classical case.

\smallbreak
This paper is sequel to \cite{jula1}, but we meant to keep it self--contained, so that the reader can appreciate the subtle differences between the different $p$--adic objects discussed in the paper.

\smallbreak
We would like to thank  Drossos Gintides and Johannes Grassberger for helping us find interesting solutions to the $E$--system using mathematical software.  Also, we are thankful to Paul G\'{e}rardin for computing the general solution of the $E$--system. Finally, we thank the Referee for very valuable remarks that led to improving the quality of the paper.

\section{Framed braids, quotient algebras and $p$--adic objects}

\subsection{\it Framed braid groups}

The classical braid group on $n$ strands, $B_n$, is generated
by the elementary braids $\sigma_1,\ldots ,\sigma_{n-1}$, where
$\sigma_i$ is the positive crossing between the $i$th and the
$(i+1)$st strand, satisfying the well--known braid relations:
 $\sigma_{i}\sigma_{i+1}\sigma_{i}= \sigma_{i+1}\sigma_{i}\sigma_{i+1}$ and $\sigma_i\sigma_j = \sigma_j\sigma_i$ for $\vert i-j\vert >1$. On the other hand the group
${\Bbb Z}^n$ is generated by the  \lq elementary framings\rq\,
$(0, \ldots,0,1,0, \ldots,0)$ with  $1$ in the $i$th
position. In the multiplicative notation an element $a=(a_1, \ldots,a_n) \in {\Bbb Z}^n$ can be expressed as $a= t_1^{a_1}\ldots t_{n}^{a_n}$, where $t_1,\ldots, t_{n}$ are the standard multiplicative generators of ${\Bbb Z}^n$. The {\it framed braid group} on $n$ strands is then defined as:
\begin{equation}\label{action}
{\mathcal F}_{n} = {\Bbb Z}^n \rtimes  B_n
\end{equation}
where the action of $B_n$ on ${\Bbb Z}^n$ is given by the permutation induced by a braid on the indices: $\sigma_it_j=t_{\sigma_i(j)}\sigma_i$. A word $w$ in ${\mathcal F}_{n}$ has, thus, the `splitting property', i.e. it splits into the \lq framing\rq \, part and the \lq
braiding\rq \, part:
$w = t_1^{a_1}\ldots t_n^{a_n} \, \sigma$, \ where $\sigma \in B_n$. So $w$ is a
classical braid with an integer, its framing, attached to each
strand. Especially, an element of ${\Bbb Z}^n$ is identified with
a framed identity braid on $n$ strands, while a classical braid in $B_n$ is
viewed as a framed braid with all framings 0. The multiplication in
${\mathcal F}_n$ is defined by concatenating the braids and
collecting the total framing of each strand to the top. For a good treatment of the group ${\mathcal F}_n$ see, for example, \cite{KS}.

 Further, for a positive integer $d$, the {\it $d$--modular framed braid
group} on $n$ strands, ${\mathcal F}_{d,n}$, is defined as the
quotient of ${\mathcal F}_n$ over the {\it modular relations}:
\begin{equation}\label{modular}
t_i^d= 1 \quad (i =1, \ldots, n)
\end{equation}
Thus, ${\mathcal F}_{d,n}=({\Bbb Z}/d{\Bbb Z})^n \rtimes  B_n$.
Framed braids in ${\mathcal F}_{d,n}$  have framings modulo $d$.

\subsection{\it Yokonuma--Hecke algebras}

In the group algebra ${\Bbb C} {\mathcal F}_{d,n}$ we have the following elements, which are idempotents (see for example Lemma~4\cite{jula1}).
\begin{equation}\label{edi}
e_{d,i} := \frac{1}{d} \sum_{m=0}^{d-1}t_i^m t_{i+1}^{d-m}
\qquad(i=1,\ldots , n-1)
\end{equation}
In fact $e_{d,i} \in {\Bbb C}({\Bbb Z}/d{\Bbb Z})^n$. Figure~\ref{ed1} illustrates a diagrammatic interpretation of $e_{d,1}  \in {\Bbb C} {\mathcal F}_{d,3}$.

\smallbreak
\begin{figure}[H]
\begin{picture}(320,60)

\put(0,38){$e_{d,1} =$}
\put(37, 41){$1$}
\qbezier(36,40)(41,40)(46,40)
\put(36,30){$d$}
\qbezier(55,20)(45,40)(55,60) 
\qbezier(65,20)(65,40)(65,60)
\qbezier(80,20)(80,40)(80,60)
\qbezier(95,20)(95,40)(95,60)
\put(105,40){$+$}
\qbezier(125,20)(125,40)(125,60)
\qbezier(140,20)(140,40)(140,60)
\qbezier(155,20)(155,40)(155,60)
\put(165,40){$+$}
\qbezier(185,20)(185,40)(185,60)
\qbezier(200,20)(200,40)(200,60)
\qbezier(215,20)(215,40)(215,60)
\put(225,40){$+$}
\put(240,40){$\cdots$}
\put(260,40){$+$}

\qbezier(280,20)(280,40)(280,60)
\qbezier(295,20)(295,40)(295,60)
\qbezier(310,20)(310,40)(310,60)
\qbezier(320,20)(330,48)(320,60) 
\put(62,65){\tiny{$0$}}
\put(77,65){\tiny{$0$}}
\put(92,65){\tiny{$0$}}
\put(122,65){\tiny{$1$}}
\put(130,65){\tiny{$d-1$}}
\put(153,65){\tiny{$0$}}
\put(182,65){\tiny{$2$}}
\put(190,65){\tiny{$d-2$}}
\put(213,65){\tiny{$0$}}
\put(270,65){\tiny{$d-1$}}
\put(292,65){\tiny{$1$}}
\put(307,65){\tiny{$0$}}

\end{picture}
\caption{The element $e_{d,1} \in {\Bbb C} {\mathcal F}_{d,3}$} \label{ed1}
\end{figure}

In the following we fix  $u \in {\Bbb C}\backslash \{0,1\}$. The {\it Yokonuma--Hecke algebra} ${\rm Y}_{d,n}(u)$ is defined as the quotient of the group algebra
${\Bbb C} {\mathcal F}_{d,n}$ over the ideal $I_{d,n}$ generated by the expressions $\sigma_i^2 - 1 - (u-1) \, e_{d,i} - (u-1) \, e_{d,i} \, \sigma_i$, which give rise to the following quadratic relations (corresponding $\sigma_i$ to $g_i$):
\begin{equation}\label{quadr}
g_i^2 = 1 + (u-1) \, e_{d,i} + (u-1) \, e_{d,i} \, g_i
\end{equation}
 (see \cite{jula1} for diagrammatic interpretations). Since the quadratic relations do not change the framing, we have ${\Bbb C}({\Bbb Z}/d{\Bbb Z})^n\subset {\rm Y}_{d,n}(u)$ and we keep the same notation for the elements of ${\Bbb C}({\Bbb Z}/d{\Bbb Z})^n$ in ${\rm Y}_{d,n}(u)$. In particular we use the same notation for the elements $e_{d,i}$ in ${\rm Y}_{d,n}(u)$.
 The elements $g_i$ are invertible (see Figure~\ref{g1invrs} for a diagrammatic interpretation):
\begin{equation}\label{invrs}
g_i^{-1} = g_i + (u^{-1} - 1)\, e_{d,i} + (u^{-1} - 1)\, e_{d,i}\, g_i
\end{equation}

\smallbreak
\begin{figure}[H]
\setlength{\unitlength}{.8pt}
\begin{picture}(400,150)

\qbezier(-50,100)(-51,105)(-50,110)
\qbezier(-30,100)(-29,105)(-30,110)

\qbezier(-50,110)(-50,114)(-45,117)
\qbezier(-35,123)(-30,126)(-30,130)
\qbezier(-30,110)(-30,115)(-40,120)
\qbezier(-40,120)(-50,125)(-50,130)

\qbezier(-50,130)(-49,135)(-50,140)
\qbezier(-30,130)(-29,135)(-30,140)
\qbezier(-10,100)(-10,120)(-10,140)
\put(0,117){$=$}

\qbezier(20,100)(18,105)(20,110)
\qbezier(40,100)(41,105)(40,110)
\qbezier(20,110)(20,115)(30,120)
\qbezier(30,120)(40,125)(40,130)
\qbezier(40,110)(40,114)(35,117)
\qbezier(25,123)(20,126)(20,130)
\qbezier(20,130)(19.5,135)(20,140)
\qbezier(40,130)(41,135)(40,140)

\qbezier(57,100)(57,120)(57,140)

\put(64,117){\small{$+\frac{u^{-1}-1}{d}$}}
\qbezier(117,100)(107,120)(117,140) 
\qbezier(130,100)(130,120)(130,140)
\qbezier(148,100)(148,120)(148,140)
\qbezier(166,100)(166,120)(166,140)
\put(180,117){$+$}
\qbezier(203,100)(203,120)(203,140)
\qbezier(221,100)(221,120)(221,140)
\qbezier(239,100)(239,120)(239,140)
\put(250,117){$+$}
\qbezier(273,100)(273,120)(273,140)
\qbezier(291,100)(291,120)(291,140)
\qbezier(309,100)(309,120)(309,140)
\put(323,117){$+\cdots +$}
\qbezier(380,100)(380,120)(380,140)
\qbezier(398,100)(398,120)(398,140)
\qbezier(416,100)(416,120)(416,140)
\qbezier(425,100)(435,120)(425,140) 

\put(-3,37){\small{$+\frac{u^{-1}-1}{d}$}}

\qbezier(50,20)(40,40)(50,60) 

\qbezier(60,50)(59,55)(60,60)
\qbezier(80,50)(81,55)(80,60)

\qbezier(60,20)(59,25)(60,30)
\qbezier(80,20)(81,25)(80,30)

\qbezier(60,30)(60,35)(70,40)
\qbezier(70,40)(80,45)(80,50)
\qbezier(80,30)(80,34)(75,37)
\qbezier(65,43)(60,46)(60,50)

\qbezier(95,20)(95,40)(95,60)

\put(115,37){$+$}

\qbezier(150,50)(149,55)(150,60)
\qbezier(170,50)(171,55)(170,60)

\qbezier(150,20)(149,25)(150,30)
\qbezier(170,20)(171,25)(170,30)

\qbezier(150,30)(150,35)(160,40)
\qbezier(160,40)(170,45)(170,50)
\qbezier(170,30)(170,34)(165,37)
\qbezier(155,43)(150,46)(150,50)
\qbezier(185,20)(185,40)(185,60)

\put(205,37){$+$}

\qbezier(230,50)(229,55)(230,60)
\qbezier(250,50)(251,55)(250,60)

\qbezier(230,20)(229,25)(230,30)
\qbezier(250,20)(251,25)(250,30)

\qbezier(230,30)(230,35)(240,40)
\qbezier(240,40)(250,45)(250,50)
\qbezier(250,30)(250,34)(245,37)
\qbezier(235,43)(230,46)(230,50)
\qbezier(265,20)(265,40)(265,60)
\put(285,37){$+\cdots +$}

\qbezier(340,50)(339,55)(340,60)
\qbezier(360,50)(361,55)(360,60)

\qbezier(340,20)(339,25)(340,30)
\qbezier(360,20)(361,25)(360,30)

\qbezier(340,30)(340,35)(350,40)
\qbezier(350,40)(360,45)(360,50)
\qbezier(360,30)(360,34)(355,37)
\qbezier(345,43)(340,46)(340,50)

\qbezier(375,20)(375,40)(375,60)

\qbezier(385,20)(395,40)(385,60) 
\put(-50,142){\tiny{$0$}}
\put(-32,142){\tiny{$0$}}
\put(-11,142){\tiny{$0$}}

\put(18,142){{\tiny $0$}}
\put(38,142){\tiny{$0$}}
\put(55,142){\tiny{$0$}}

\put(127,142){\tiny{$0$}}
\put(146,142){\tiny{$0$}}
\put(164,142){\tiny{$0$}}

\put(200,142){\tiny{$1$}}
\put(213,142){\tiny{$d\!-\!1$}}
\put(237,142){\tiny{$0$}}

\put(270,142){\tiny{$2$}}
\put(282,142){\tiny{$d\!-\!2$}}
\put(307,142){\tiny{$0$}}

\put(366,142){\tiny{$d\!-\!1$}}
\put(395,142){\tiny{$1$}}
\put(414,142){\tiny{$0$}}

\put(58,63){\tiny{$0$}}
\put(77,63){\tiny{$0$}}
\put(93,63){\tiny{$0$}}

\put(148,63){\tiny{$1$}}
\put(160,63){\tiny{$d\!-\!1$}}
\put(183,63){\tiny{$0$}}

\put(228,63){\tiny{$2$}}
\put(240,63){\tiny{$d\!-\!2$}}
\put(264,63){\tiny{$0$}}

\put(331,63){\tiny{$d\!-\!1$}}
\put(358,63){\tiny{$1$}}
\put(372,63){\tiny{$0$}}
\end{picture}
\caption{The element $g_1^{-1} \in {\rm Y}_{d,3}(u)$}\label{g1invrs}
\end{figure}

\begin{note} \rm We note that in \cite{yo}, \cite{ju}, \cite{jula1,jula2,jula3,jula4} instead of Eq.~\ref{quadr} the following quadratic relation is used:
\begin{equation}\label{oldquadr}
g_i^2 = 1+(u-1)e_i - (u-1)e_ig_i.
\end{equation}
\end{note}

From the above, a presentation of  ${\rm Y}_{d,n}(u)$ is given by the generators $t_1,\ldots,t_n$, $g_1,\ldots,g_{n-1}$, satisfying: the braid relations and the quadratic relations (\ref{quadr}) for the $g_i$'s, the modular relations (\ref{modular}) and commuting relations for the $t_j$'s, together with the mixed relations below, deriving from (\ref{action}):
$$
\begin{array}{rcl}
g_i t_i & = & t_{i+1} g_i \\
g_i t_{i+1} & = & t_i g_i \\
g_i t_j & = & t_j g_i \quad \text{ for } j\not= i, i+1
\end{array}
$$
Note that, omitting the quadratic relations (\ref{quadr}), we have a presentation for ${\mathcal F}_{d,n}$.

\begin{rem}\label{YtoH}\rm
There is an epimorphism of the Yokonuma--Hecke algebra ${\rm Y}_{d,n}(u)$ onto the Iwahori--Hecke algebra ${\rm H}_n(q)$ via the map
$$
\begin{array}{lll}
g_i & \mapsto & G_i \\
t_j & \mapsto & 1
\end{array}
$$
where $G_i$ are the standard generators of ${\rm H}_n(q)$. Further, for $d=1$ we have all $t_j= 1$ and  ${\rm Y}_{1,n}(u)$ coincides with the algebra ${\rm H}_n(u)$. Also, the mapping $g_i  \mapsto (i,i+1)$ and $t_j \mapsto 1$  defines an epimorphism of ${\rm Y}_{d,n}(1)$ onto the group algebra of the symmetric group.
\end{rem}

In ${\rm Y}_{d,n}(u)$ the following relations hold (see Lemma 4, Proposition 5\cite{jula1}):
\begin{equation}\label{edikrels}
\begin{array}{rcll}
t_j e_{d,i} & = & e_{d,i} t_j & \\
e_{d,j} e_{d,i} & = & e_{d,i} e_{d,j} & \\
g_j e_{d,i} & = & e_{d,i} g_j & \text{ for } j\not= i-1, i+1 \\
g_{i-1} e_{d,i} & = & e_{d,i-1,i+1} g_{i-1}
& \text{ and }\quad e_{d,i} g_{i-1} =
g_{i-1} e_{d,i-1,i+1}\\
g_{i+1} e_{d,i} & = & e_{d,i,i+2} g_{i+1}
& \text{ and }\quad  e_{d,i} g_{i+1}
= g_{i+1} e_{d,i,i+2}
\end{array}
\end{equation}
where:
\begin{equation}\label{edik}
e_{d,i,k} = \frac{1}{d} \sum_{1\leq s \leq d}t_i^st_{k}^{d-s}
\end{equation}
for any $i, k$ with $k\not=
i$, abbreviating $e_{d, i, i+1}$ to $e_{d, i}$.  Note that, using (\ref{invrs}), relations (\ref{edikrels}) are also valid if all the $g_r$'s are replaced by their inverses $g_r^{-1}$. Clearly
$e_{d,i,k} = e_{d,k,i}$ and it can be easily deduced that $e_{d,i,k}^2 = e_{d,i,k}$.

\begin{lem}\label{powers}
Let $m \in {\Bbb Z}$. Then the following relations hold in ${\rm Y}_{d,n}(u)$.
\begin{enumerate}
\item For $m$ positive we have:
\begin{eqnarray*}
g_i^{m}
&  = &
1 + \frac{u^m -1}{u+1} \, e_{d,i}\, g_i + \frac{u^m -1}{u+1} \, e_{d,i},
\quad  \text{ \ if \ } m= \text{even }\\
g_i^{m}
&  = &
g_i + \frac{u^{m} -u}{u+1} \, e_{d,i} g_i + \frac{u^{m} -u}{u+1} \, e_{d,i},
\quad  \text{ \ if \ } m= \text{odd }
\end{eqnarray*}
\item For $m$ negative we have:
\begin{eqnarray*}
g_i^{m}
&  = &
1 + \frac{u^{m-1} - u^{-1}}{u^{-1} +1} \, e_{d,i}\, g_i + \frac{u^{m-1} - u^{-1}}{u^{-1} +1}\, e_{d,i},
\quad  \text{ \ if \ } m= \text{even }\\
g_i^{m}
&  = &
g_i + \frac{u^{m-1} -1}{u^{-1} +1} \, e_{d,i} g_i + \frac{u^{m-1} -1}{u^{-1} +1}\, e_{d,i},
\quad  \text{ \ if \ } m= \text{odd }
\end{eqnarray*}
\end{enumerate}
\end{lem}

\begin{proof} By induction. We shall only check the case when $m$ positive. The case $m$ negative follows similarly. For  $m=1$ the statement is clearly true. For $m=2$ we have Eq.~\ref{quadr} and we note that $u-1 = \frac{u^2 -1}{u+1}$, so the statement is also true. Suppose the statement holds for all $m$ up to $2k-1$, $k \in {\Bbb N}$. Then, for $m=2k$ we have:
\begin{eqnarray*}
g_i^{2k}
&  = &
g_i^{2k-1} g_i \\
& = &
\left[g_i + \frac{u^{2k-1} -u}{u+1}  \, e_{d,i}\, g_i +   \frac{u^{2k-1} -u}{u+1}  \, e_{d,i} \right] g_i \\
& \stackrel{(\ref{quadr})}{=} & \left[1 + \frac{u^{2k-1} -u}{u+1}  \, e_{d,i}\right] \left[1 + (u-1) \, e_{d,i} + (u-1) \, e_{d,i} \, g_i \right]  +  \frac{u^{2k-1} -u}{u+1}  \, e_{d,i} \, g_i
 \\
& = &
1 +  \frac{u^{2k} -1}{u+1} \, e_{d,i}\, g_i + \frac{u^{2k} -1}{u+1} \, e_{d,i}.
\end{eqnarray*}

\noindent Also, for $m=2k+1$ we have:

\begin{eqnarray*}
g_i^{2k+1}
&  = &
g_i^{2k} g_i \\
& = &
\left[1 +  \frac{u^{2k} -1}{u+1} \, e_{d,i}\, g_i + \frac{u^{2k} -1}{u+1} \, e_{d,i} \right] g_i \\
& \stackrel{(\ref{quadr})}{=} &
g_i + \frac{u^{2k} -1}{u+1} \, e_{d,i}\, g_i + \frac{u^{2k} -1}{u+1} \, e_{d,i}\, \left[1 + (u-1) \, e_{d,i} + (u-1) \, e_{d,i} \, g_i \right] \\
& = &
g_i + \frac{u^{2k+1} -u}{u+1} \, e_{d,i}\, g_i + \frac{u^{2k+1} - u}{u+1} \, e_{d,i}.
\end{eqnarray*}
\end{proof}

\subsection{\it Inverse limits and the $p$--adic integers}

Our references for inverse limits are mainly \cite{riza} and \cite{wi}.  An {\it inverse system} $(X_i, \phi_j^i)$ of topological spaces  indexed by a directed set $I$, consists of a family $(X_i \ ; \ i\in I)$ of topological spaces (groups, rings, algebras, et cetera) and a family
$(\phi_j^i: X_i \longrightarrow X_j \ ; \  i,j \in I, \ i\geq j)$ of continuous homomorphisms, such that
$$
\phi_i^i = {\rm id}_{X_i} \quad{\rm and} \quad  \phi_k^j \circ \phi_j^i  = \phi_k^i
\quad {\rm whenever}  \quad
i \geq j \geq k
$$
The maps $\phi_j^i$ are also called {\it connecting homomorphisms}.
If no other topology is specified on the sets $X_i$ they are regarded as topological spaces  with the discrete topology. In particular, finite sets are compact Hausdorff spaces. The {\it inverse limit} $\varprojlim X_i$ of the inverse  system $(X_i, \phi_j^i)$  is defined as:
$$
\varprojlim X_i := \{ z \in \prod_{i\in I} X_i\, ; \,  (\phi_j^i\circ \varpi_i )(z) = \varpi_j (z)
\quad \text{ whenever}  \quad j\geq i\}
$$
where the map $\varpi_i$ denotes the natural projection of
$\prod X_i$ onto $X_i$. Recall that, if $X_i =X$ for all $i$ and $\phi_j^i$ is the identity for all $i,j$ then
$\varprojlim X$ can be identified naturally with $X$ (identifying a constant sequence $(x,x,\ldots)$ with $x\in X$).

\begin{notation} \rm
In the following we fix a prime number $p$ and we   denote by ${\Bbb N}$ the set of positive integers regarded as a directed set with the usual order. Further, for $r\geq s$ we denote  $\vartheta^{r}_s$ the natural epimorphism:
\begin{equation}\label{varthetars}
\begin{array}{cccc}
 \vartheta^{r}_s : & {\Bbb Z}/{p^r}{\Bbb Z} & \longrightarrow & {\Bbb Z}/p^s{\Bbb Z}\\
 & m  & \mapsto &   m \,( {\rm mod}\, p^s)
\end{array}
\end{equation}
\end{notation}

We now denote
$$
C_r  =  \langle t_r \  ; \, t_r^{p^r} = 1 \rangle \cong {\Bbb Z} / p^r{\Bbb Z}
$$
 the cyclic group of order $p^r$ in the
multiplicative notation. For $r\geq s$ we denote  $\theta_s^r$ the following natural connecting epimorphism,
\begin{equation}\label{thetars}
\begin{array}{cccc}
\theta_s^r : &  C_r &  \longrightarrow & C_s \\
            \,&  t_r^m & \mapsto &   t_s^{\vartheta^r_s(m) }
\end{array}
\end{equation}
Thus we obtain an inverse system of groups $(C_r , \theta_s^r)$, whose inverse limit is the group of {\it $p$--adic integers} ${\Bbb Z}_p$,
\begin{center}
 $
 {\Bbb Z}_p:= \varprojlim_{r\in{\Bbb N}} C_r.
 $
 \end{center}
 The group ${\Bbb Z}_p$  can be regarded as:
$$
{\Bbb Z}_p =  \{ {\bold t}^{\underleftarrow{a}} := (t_1^{a_1},
t_2^{a_2}, \ldots) \, \in \prod_{i\in {\Bbb N}} C_i \ ; \ a_r \in  {\Bbb Z}, \ a_r
\equiv a_s \ (\text{mod} \, p^s) \ \ \text{whenever} \ r\geq s\}.
$$
Notice that the multiplication in ${\Bbb Z}_p$ is then defined as:
$$
{\bold t}^{\underleftarrow{a}} {\bold t}^{\underleftarrow{b}} = {\bold t}^{\underleftarrow{a} + \underleftarrow{b}} = (t_1^{a_1 + b_1},
t_2^{a_2 + b_2}, \ldots)
$$
In ${\Bbb Z}_p$ the `classical integers' are the coherent sequences which are eventually constant. The element ${\bold t} : =  (t_1, t_2, \ldots) \in {\Bbb Z}_p$ corresponds to  $(1, 1, \ldots)$ in the additive notation, so it generates in ${\Bbb Z}_p$ a copy of ${\Bbb Z}$ and we can write ${\Bbb Z} = \langle {\bold t}\rangle$.

\begin{defn}\label{unique}
\rm An element $\underleftarrow{a} = (a_1, a_2, \ldots) \in {\Bbb Z}_p$ is said to be in {\it reduced form} if each entry $a_r \in {\Bbb Z} / p^r{\Bbb Z}$ is expressed in its (unique) reduced $p$--adic expansion:
$$
a_r = k_0 + k_1 p + k_2 p^2 + \cdots + k_{r-1} p^{r-1} + p^r{\Bbb Z}
$$
where $ k_0,\ldots,k_{r-1} \in \{0,1,\ldots,p-1\}$.
In the multiplicative notation this means that the exponents $a_r$ of $t_r^{a_r}$ are in the above reduced form.
\end{defn}

\subsection{\it $p$--adic framed braids}

Consider now the  group $C_r^n$. We define in $C_r^n$ the elements:
$$
t_{r,i} := (1,\ldots,1, t_r,1,\ldots,1)
$$
with $t_r$ in the $i$th position. Then we have
$$
C_r^n=  \langle t_{r,1}, t_{r,2}, \ldots , t_{r,n} \  ; \, t_{r,i}^{p^r} = 1,\,
 t_{r,i}t_{r,j} =  t_{r,j}t_{r,i}\,\,\text{for all}\,\, i, j \rangle
$$
(notice that  $C_r^n \cong ({\Bbb Z} / p^r{\Bbb Z})^n$). By componentwise multiplication, the epimorphisms (\ref{varthetars}) define the connecting epimorphisms:
$$
\begin{array}{cccc}
\pi_s^r : &  C_r^n  &  \longrightarrow & C_s^n \\
         &  t_{r,i}^m  & \mapsto &  t_{s,i}^{\vartheta^r_s(m)}
\end{array}
$$
for all   $r\geq s$.
 Extending to the $B_n$--part by the identity map gives rise to the connecting epimorphisms:
\begin{equation}\label{pirsid}
\begin{array}{cccc}
 \pi^r_s \cdot {\rm id} : & {\mathcal F}_{p^r,n} & \longrightarrow & {\mathcal F}_{p^s,n} \\
    &  t_{r,i}^m  & \mapsto &  t_{s,i}^{\vartheta^r_s(m)} \\
      &  \sigma_i  & \mapsto &  \sigma_i
\end{array}
\end{equation}
In \cite{jula1} we defined the {\it $p$--adic framed braid group} on $n$ strands ${\mathcal
F}_{p^{\infty} ,n}$ as:
$$
{\mathcal F}_{p^{\infty} ,n} := \varprojlim_{r\in{\Bbb N}} {\mathcal F}_{p^r,n}
$$
Notice that ${\mathcal F}_{p^r,n} = C_r^n \rtimes B_n$. So, a $p$--adic framed braid is an infinite sequence of modular framed braids with the same braiding part and such that the framings of the $i$th strand in each position of the sequence
give rise to a $p$--adic integer. See Figure~\ref{facets}. Elements of ${\mathcal F}_{p^{\infty} ,n}$ are  denoted $\underleftarrow{\beta}$.
 We recall now from Proposition~4\cite{jula1} that there are group isomorphisms:
\begin{equation}\label{isoms}
{\mathcal F}_{p^{\infty} ,n} \cong {\Bbb Z}_p^n \rtimes B_n \cong (\varprojlim_{r\in{\Bbb N}} C_r^n) \rtimes B_n
\end{equation}
Hence, a $p$--adic framed braid may also be viewed as a classical framed braid but with framings $p$--adic integers or, equivalently, as a classical braid with infinite cablings replacing each strand, whose corresponding modular framings form a $p$--adic integer.
The $n$--tuples of constant sequences form the subgroup ${\Bbb Z}^n = \langle {\bold t}_1, \ldots, {\bold t}_n \rangle$, where ${\bold t}_i := ({\bold 1},\ldots,
{\bold 1}, {\bold t}, {\bold 1}, \ldots, {\bold 1})$ with ${\bold t}$ in the $i$th position. Note that the element ${\bold 1} := ({\bold 1},\ldots, {\bold 1})$ corresponds to the identity framed braid with all framings zero.
 Note further that  $(\sigma_i, \sigma_i, \ldots ) \in {\mathcal F}_{p^{\infty} ,n}$ gets identified with $\sigma_i \in B_n$ in the first isomorphism and that ${\bold t}_i \in { \Bbb Z} \subset { \Bbb Z}_p^n$ gets identified with $(t_{r,i})_r \in \varprojlim_{r}C_r^n$ in the second isomorphism.

In view of the first isomorphism, a $p$--adic framed braid splits into the  \lq $p$--adic
framing\rq \, part and the \lq braiding\rq \, part: ${\bold
t}_1^{\underleftarrow{a_1}}
\ldots {\bold t}_n^{\underleftarrow{a_n}} \, \sigma$, that is, to
each strand of the braid $\sigma \in B_n$ we attach a $p$--adic
integer (see Figure~\ref{facets}). $p$--adic framed braids are
multiplied by concatenating their braiding parts and collecting
the total $p$--adic framing of each strand to the top:
\[
({\bold t}_1^{\underleftarrow{a_1}}\ldots {\bold t}_{n}^{\underleftarrow{a_n}}\,
\sigma)({\bold t}_1^{\underleftarrow{b_1}}\ldots {\bold t}_{n}^{\underleftarrow{b_n}}\, \tau) := {\bold t}_1^{\underleftarrow{a_1} +
\underleftarrow{b_{\sigma(1)}}} \ldots {\bold t}_{n}^{\underleftarrow{a_n} + \underleftarrow{b_{\sigma(n)}}} \, \sigma\tau
\]
As already mentioned, isomorphisms (\ref{isoms}) imply that a $p$--adic framed braid can be interpreted as a classical braid with framings $p$--adic integers or as a classical braid, but with infinite cablings replacing each strand, such that the framings of each infinite cable form a $p$--adic integer (Figure~\ref{facets}). In the sequel we shall not distinguish either  between the isomorphic forms of ${\mathcal F}_{p^{\infty},n}$ or between the different interpretations of corresponding elements in them.

Note, finally, that the natural inclusions ${\mathcal F}_n  \subset {\mathcal F}_{n+1}$ of the classical framed braid groups induce natural inclusions ${\mathcal F}_{d,n} \subset {\mathcal F}_{d,n+1}$ of the modular framed braid groups and these induce the natural inclusions  ${\mathcal F}_{p^{\infty},n}  \subset {\mathcal F}_{p^{\infty},n+1}$ on the level of the $p$--adic framed braid groups.

\subsection{\it The $p$--adic Yokonuma--Hecke algebra}

For all $r\geq s$ the linear extension of the map  (\ref{pirsid}) yields a connecting algebra epimorphism:
\begin{equation}\label{varphi}
\varphi ^{r}_s : {\Bbb C} {\mathcal F}_{p^r,n}\longrightarrow {\Bbb C} {\mathcal F}_{p^s,n}
\end{equation}
Now, passing to the quotient algebras, we obtain the following connecting algebra epimorphism
$$
\phi_s^r : {\rm Y}_{p^r,n}(u) \longrightarrow {\rm Y}_{p^s,n}(u)
$$
(cf. \cite{jula1} for details of the construction).
So we obtain the inverse system $({\rm Y}_{p^r, n}(u), \phi_s^r)$.
In \cite{jula1} the {\it $p$--adic Yokonuma--Hecke algebra} ${\rm Y}_{p^{\infty}, n}(u)$ was defined  as:
$$
{\rm Y}_{p^{\infty} ,n}(u): = \varprojlim_{r\in{\Bbb N}} {\rm Y}_{p^r,n}(u).
$$

\begin{note} \rm By construction, $\varprojlim_{r} {\Bbb C}{C_r^n} \subset \varprojlim_{r} {\Bbb C}{\mathcal F}_{p^r,n}$ and  $\varprojlim_{r} {\Bbb C}{C_r^n} \subset {\rm Y}_{p^{\infty}, n}(u)$.
 Note also that there are no modular relations in ${\mathcal F}_{p^{\infty} ,n}$ as well as in ${\Bbb C} {\mathcal F}_{p^{\infty} ,n}$, in $\varprojlim_{r} {\Bbb C}{\mathcal F}_{p^r,n}$ and in  ${\rm Y}_{p^{\infty},n}(u)$.
\end{note}

\begin{defn}\label{reduced}\rm
We shall say that an element in ${\mathcal F}_{p^r,n}$ is in its (unique) {\it reduced form } if its modular framings are in their (unique) reduced  $p$--adic expansion.
Then, by the linear extension, an element in ${\Bbb C}{C_r^n}$ or in ${\Bbb C} {\mathcal F}_{p^r,n}$ has a (unique) {\it reduced form}. Further, an element $y+I_{p^r,n}$ in ${\rm Y}_{p^r,n}(u)$ is in {\it reduced form} if the element $y \in {\Bbb C} {\mathcal F}_{p^r,n}$ is written in its (unique) reduced form.
\end{defn}

\begin{defn}\label{preduced}\rm
An element $\underleftarrow{\beta}=(\beta_1, \beta_2, \ldots) \in {\mathcal F}_{p^{\infty},n}$ is said to be in its (unique) {\it reduced form} if every entry $\beta_r \in {\mathcal F}_{p^r,n}$ is reduced according to Definition~\ref{reduced}. In view of the first isomorphism in (\ref{isoms}), we may also say that $\underleftarrow{\beta} = {\bold t}_1^{\underleftarrow{a_1}}\ldots {\bold t}_{n}^{\underleftarrow{a_n}}\, \sigma$ is in reduced form if its $p$--adic framings $\underleftarrow{a_1}, \ldots, \underleftarrow{a_n}$ are reduced according to Definition~\ref{unique}. Further, by the linear expansion on ${\mathcal F}_{p^{\infty} ,n}$, an element in ${\Bbb C} {\mathcal F}_{p^{\infty} ,n}$ has a (unique) {\it reduced form}.
 An element in  $\varprojlim_{r}{\Bbb C}{C_r^n}$, in  $\varprojlim_{r} {\Bbb C}{\mathcal F}_{p^r,n}$ or  in ${\rm Y}_{p^{\infty} ,n}(u)$ is said to be in {\it reduced form} if every entry is reduced in ${\Bbb C}{C_r^n}$, in ${\Bbb C} {\mathcal F}_{p^r,n}$ or in ${\rm Y}_{p^r,n}(u)$ respectively, according to Definition~\ref{reduced}.
\end{defn}

\subsection{\it The elements $e_{i}$}\label{eltsei}

We define now the elements
\begin{equation}\label{ei}
e_{i}:= (e_{p,i}, e_{p^2,i}, \ldots )
\end{equation}
  where
\begin{equation}\label{epri}
e_{p^r,i} = \frac{1}{p^r} \sum_{m=0}^{p^r-1}t_{r,i}^m t_{r,i+1}^{-m} \in {\Bbb C}{C_r^n} \subset {\Bbb C}{\mathcal F}_{p^r,n}
\end{equation}

\begin{lem}\label{eprtoeps}
 $ e_{i}\in \varprojlim_{r} {\Bbb C}{C_r^n} \subset \varprojlim_{r} {\Bbb C}{\mathcal F}_{p^r,n}$ and $e_{i} \in {\rm Y}_{p^{\infty},n}(u)$ for $i= 1, \ldots, n-1$.
\end{lem}

\begin{proof}
We shall show the coherency of the terms in $e_{i}$, that is: $\varphi ^{r}_s  (e_{p^r, i}) = e_{p^s, i}$ $(r\geq s)$. Note first that by the maps (\ref{varphi}) and (\ref{pirsid}):
$$
\varphi ^{r}_s ( t_{r,i}^m t_{r,i+1}^{-m} ) = \pi_s^r ( t_{r,i}^m t_{r,i+1}^{-m} ) = t_{s,i}^m t_{s,i+1}^{-m}
$$
But, by (\ref{epri}), $e_{p^r,i}$ is a sum of $p^r$ terms with linear coefficients
$\frac{1}{p^r}$ and $e_{p^s,i}$ is a sum of $p^s$ terms with linear
coefficients $\frac{1}{p^s}$. Yet,
\[
e_{p^r,i} =\frac{1}{p^{r-s}} \left( \sum_{m=0}^{p^s-1}
\frac{1}{p^s} t_{r,i}^m t_{r,i+1}^{-m} + \sum_{m=p^s}^{2p^s-1}
\frac{1}{p^s} t_{r,i}^m t_{r,i+1}^{-m} +\cdots+
\sum_{m=p^r-p^s}^{p^r-1} \frac{1}{p^s} t_{r,i}^m t_{r,i+1}^{-m}
\right).
\]
The element in ${\Bbb C}{C_r^n}$ in each one of the
$p^{r-s}$ sums maps to $e_{p^s,i}$ in ${\Bbb C}{C_s^n}$, so $\varphi ^{r}_s  (e_{p^r, i}) = e_{p^s, i}$. Moreover, by the definition of the Yokonuma--Hecke algebra it follows also that $\phi_s^r (e_{p^r, i}) = e_{p^s, i}$. Thus $e_{i} \in {\rm Y}_{p^{\infty},n}(u)$.
\end{proof}

 The elements $e_i$ are no more averaged sums but they are still idempotents. Further, setting by construction $g_i:= (g_i, g_i,\ldots)$ and $1:= (1, 1, \ldots )$ we have in ${\rm Y}_{p^{\infty},n}(u)$ the braid relations for the $g_i$'s and the relations:
\begin{equation}\label{pquadr}
g_i^2 = 1 + (u-1)e_{i} + (u-1)e_{i}\, g_i
\end{equation}
and
$$
g_i^{-1} = g_i + (u^{-1} - 1)\, e_{i} + (u^{-1} - 1)\, e_{i}\, g_i
$$
Moreover, for powers of $g_i$ relations analogous to the ones in Lemma~\ref{powers} are valid in ${\rm Y}_{p^{\infty}, n}(u)$, after replacing $e_{d,i}$ by $e_{i}$.
Finally, using the elements $e_{d,i,k}$ in Eq.~\ref{edik}, we can define for $i= 1, \ldots , n-1$ and $k\not=i$ the elements:
\begin{equation}
e_{i,k}:= (e_{p,i,k}, e_{p^2,i,k}, \ldots ) \in \varprojlim_{r\in {\Bbb N}} {\Bbb C}{C_r^n} \subset \varprojlim_{r\in{\Bbb N}} {\Bbb
C}{\mathcal F}_{p^r,n}
\end{equation}
abbreviating $e_{i,i+1}$ to $e_{i}$. Clearly, $e_{i,k} = e_{k,i}$ and $e_{i,k}^2 = e_{i,k}$.
 These elements satisfy relations analogous to relations (\ref{edikrels}), after replacing $e_{d,i}$ by $e_{i}$, $e_{d,i,k}$ by $e_{i,k}$ and $t_j$ by ${\bold t}_j$ (cf. Lemma~7\cite{jula1} and Proposition~10\cite{jula1}).
These relations are also valid if all $g_k$'s are replaced by their inverses $g_k^{-1}$.

\subsection{\it Comparing braid algebras}

It is worth stressing at this point that, despite the definition ${\mathcal F}_{p^{\infty} ,n} = \varprojlim_{r} {\mathcal F}_{p^r,n}$, the algebras ${\Bbb C}{\mathcal F}_{p^{\infty} ,n}$ and $\varprojlim_{r} {\Bbb C}{\mathcal F}_{p^r,n}$ are not isomorphic. Before stating our result let us take a closer look at the two algebras. ${\Bbb C}{\mathcal F}_{p^{\infty} ,n}$ consists in all finite linear expressions of the form
$$
\lambda_1\underleftarrow{b}_1 + \ldots + \lambda_k\underleftarrow{b}_k
$$
where $k \in {\Bbb N}$, $\lambda_1, \ldots, \lambda_k \in {\Bbb C}$ and $\underleftarrow{b}_1, \ldots, \underleftarrow{b}_k \in {\mathcal F}_{p^{\infty} ,n}$. On the other hand, elements in $\varprojlim_{r} {\Bbb C}{\mathcal F}_{p^r,n}$ are coherent sequences of elements in the group algebras ${\Bbb C}{\mathcal F}_{p^r,n}$ in the sense of map (\ref{varphi}). That is, given elements $\beta_1, \ldots, \beta_k \in {\mathcal F}_{p^r,n}$ and $c_1, \ldots, c_k \in {\Bbb C}$ then
$$
\varphi ^{r}_s (c_1\beta_1 + \ldots + c_k\beta_k) = c_1 \, (\pi^r_s \cdot {\rm id}) (\beta_1) + \ldots + c_k\, (\pi^r_s \cdot {\rm id}) (\beta_k)\in {\Bbb C}{\mathcal F}_{p^s,n}
$$

By the definition of map $\varphi ^{r}_s$ it appears as though all positions of a $p$--adic element in $\varprojlim_{r} {\Bbb C}{\mathcal F}_{p^r,n}$  have the same number of coherent monomials with the same coefficients. This form of a $p$--adic element is always possible by construction, but it may be hidden, since the linear coefficients may allow for various manipulations. To appreciate this subtlety in the form of elements in $\varprojlim_{r} {\Bbb C}{\mathcal F}_{p^r,n}$ we shall consider as typical example the elements $e_{i}$ and the manipulation of their coherent monomials in the proof of Lemma~\ref{eprtoeps}.

\smallbreak
On the level of $p$--adic braids, ${\mathcal F}_{p^{\infty} ,n}$ sits naturally in
$\varprojlim_{r} {\Bbb C}{\mathcal F}_{p^r,n}$. So, we obtain a natural linear map

\begin{equation}\label{mapf}
f : {\Bbb C}{\mathcal F}_{p^{\infty} ,n} \longrightarrow \varprojlim_{r\in {\Bbb N}} {\Bbb C}{\mathcal F}_{p^r,n}
\end{equation}
which is constant, by construction, on ${\mathcal F}_{p^{\infty} ,n}$, that is,  $f(\underleftarrow{b})=\underleftarrow{b}$. By linearity, the image of $f$ is generated by all
$f(\underleftarrow{b})$ where $\underleftarrow{b} \in {\mathcal F}_{p^{\infty} ,n}$. Let us see how exactly the map $f$ works. Consider $\lambda_1, \ldots, \lambda_k \in {\Bbb C}$ and $\underleftarrow{b}_1, \ldots, \underleftarrow{b}_k$ different elements in ${\mathcal F}_{p^{\infty} ,n}$, where $\underleftarrow{b}_i = (b_{ri})_r$ with $b_{ri} \in {\mathcal F}_{p^r,n}$. Then we have:

$$
\begin{array}{lcl}
 {\Bbb C}{\mathcal F}_{p^{\infty}, n} \ni & \lambda_1\underleftarrow{b}_1 + \ldots + \lambda_k\underleftarrow{b}_k = & \\
 & \lambda_1 (b_{11}, b_{21}, \ldots ) + \ldots + \lambda_k (b_{1k}, b_{2k}, \ldots ) \stackrel{f}{\rightarrow} & \\
 & (\lambda_1b_{11}, \lambda_1b_{21}, \ldots ) + \ldots + (\lambda_kb_{1k}, \lambda_kb_{2k}, \ldots ) = & \\
& (\lambda_1b_{11} + \ldots + \lambda_kb_{1k}, \ \lambda_1b_{21}+ \ldots + \lambda_kb_{2k}, \ \ldots ) & \in \varprojlim_{r} {\Bbb C}{\mathcal F}_{p^r,n}.
\end{array}
$$

\begin{lem}\label{einotinim}
$e_{i}\notin f({\Bbb C}{\mathcal F}_{p^{\infty} ,n})$.
\end{lem}
\begin{proof}
 Suppose that $e_{i}\in f({\Bbb C}{\mathcal F}_{p^{\infty} ,n})$. Then, from the above, $e_{i} = a_1\underleftarrow{b}_1 + \ldots + a_k\underleftarrow{b}_k$ for some $a_1, \ldots, a_k \in {\Bbb C}$ and $\underleftarrow{b}_1, \ldots, \underleftarrow{b}_k \in {\mathcal F}_{p^{\infty} ,n}$. Then, from the structure of $\varprojlim_{r} {\Bbb C}{\mathcal F}_{p^r,n}$ as linear space, we have in $\varprojlim_{r} {\Bbb C}{\mathcal F}_{p^r,n}$ the equality:
$$
(e_{p,i}, e_{p^2,i}, \ldots ) =
(a_1b_{11} + \ldots + a_kb_{1k}, \ a_1b_{21}+ \ldots + a_kb_{2k}, \ \ldots )
$$
Equivalently, in each ${\Bbb C}{\mathcal F}_{p^r,n}, \ r=1,2,\ldots$, we have the equality:
\begin{equation}\label{proofnotin}
\sum_{m=1}^{p^r}t_{r,i}^m t_{r,i+1}^{-m} = \sum_{j=1}^{k} p^ra_j \, b_{rj}
\end{equation}
Since $\underleftarrow{b}_1, \ldots, \underleftarrow{b}_k$ are different in ${\mathcal F}_{p^{\infty} ,n}$ there must exist some $s\in {\Bbb N}$ such that $b_{s1}, \ldots, b_{sk}$ are different elements in ${\mathcal F}_{p^s,n}$ and this is then true for any $r\geq s$. So, there exists some $r\geq s$ such that $k<p^r$. But then (\ref{proofnotin}) states equality of two linear expressions of linearly independent elements in ${\Bbb C}{\mathcal F}_{p^r,n}$. Since $k<p^r$, it follows that all coefficients $p^ra_j$ must be equal to $1$. Subtracting we obtain a summation of terms $t_{r,i}^m t_{r,i+1}^{-m}$ equal to zero, which is a contradiction since they are linearly independent. Therefore $e_{i}\notin f({\Bbb C}{\mathcal F}_{p^{\infty} ,n})$.
\end{proof}

\begin{cor}\label{einotin}
$e_{i}\notin {\Bbb C}{\mathcal F}_{p^{\infty} ,n}$.
\end{cor}
\begin{proof}
Suppose $e_{i}\in {\Bbb C}{\mathcal F}_{p^{\infty} ,n}$. Then $e_{i} = a_1\underleftarrow{b}_1 + \ldots + a_k\underleftarrow{b}_k$ for some $a_1, \ldots, a_k \in {\Bbb C}$ and $\underleftarrow{b}_1, \ldots, \underleftarrow{b}_k$ as above. But then $f(e_{i}) = e_{i}$, not possible by Lemma~\ref{einotinim}. Hence $e_{i}\notin {\Bbb C}{\mathcal F}_{p^{\infty} ,n}$.
\end{proof}

\begin{prop}\label{algebras}
The linear map $f$ of Eq.~\ref{mapf} is injective but not surjective. Hence, ${\Bbb C}{\mathcal F}_{p^{\infty} ,n} \cong f({\Bbb C}{\mathcal F}_{p^{\infty} ,n})$ and so the algebra ${\Bbb C}{\mathcal F}_{p^{\infty} ,n}$ can be regarded as a proper subalgebra of \ $\varprojlim_{r} {\Bbb C}{\mathcal F}_{p^r,n}$.
\end{prop}
\begin{proof}
We will show that $f$ is injective. With the above notations let
$$
f(\lambda_1\underleftarrow{b}_1 + \ldots + \lambda_k\underleftarrow{b}_k) = (0,0,\ldots)
$$
 Equivalently, $\lambda_1b_{r1} + \ldots + \lambda_kb_{rk}=0$ in ${\Bbb C}{\mathcal F}_{p^r,n}$ for all $r=1,2,\ldots$. As in the proof of Lemma~\ref{einotinim}, since $\underleftarrow{b}_1, \ldots, \underleftarrow{b}_k$ are different in ${\mathcal F}_{p^{\infty} ,n}$ there must exist some $r\in {\Bbb N}$ such that $b_{r1}, \ldots, b_{rk}$ are different elements in ${\mathcal F}_{p^r,n}$. Hence they are linearly independent in ${\Bbb C}{\mathcal F}_{p^r,n}$, hence $\lambda_1 = \ldots = \lambda_k = 0$. Therefore ${\rm \text Ker} f=\{ \underleftarrow{0}\}$ and so $f$ is injective.

The fact that $f$ is not surjective follows immediately from Lemmas~\ref{eprtoeps} and \ref{einotinim}, since $e_{i}\in \varprojlim_{r} {\Bbb C}{\mathcal F}_{p^r,n}$ but $e_{i}\notin f({\Bbb C}{\mathcal F}_{p^{\infty} ,n})$.
\end{proof}

\begin{defn} \rm The elements in $\varprojlim_{r} {\Bbb C}{\mathcal F}_{p^r,n}$  which are not in ${\Bbb C}{\mathcal F}_{p^{\infty} ,n}$ shall be called {\it purely $p$--adic elements}.
\end{defn}

All purely $p$--adic elements are considered to be in their unique reduced form as defined in Definition~\ref{preduced}.

\section{Dense subsets and approximations of $p$--adic elements}

\subsection{\it A general lemma}

Our method for finding dense subsets in our  $p$--adic structures is by means of the following known result.
\begin{lem}{\rm(cf. \cite{riza}, Lemma 1.1.7.)}\label{dense}
Let $\rho_i$ denote the restriction  of the canonical projection of $\varprojlim X_i$ onto
$X_i$ on  a subset $A \subset \varprojlim X_i$.  Recall that $\varprojlim A$ can be identified with $A$. If $\rho_i(A) = X_i$ for all $i\in I$,  then $\rho (\varprojlim
 A)$  is dense in  $\varprojlim X_i$, where $\rho = \varprojlim\rho_i :   \varprojlim A \longrightarrow  \varprojlim X_i$, the induced mapping.
\end{lem}
\begin{defn}\label{generators}
\rm (cf. \cite{riza} \S \, 2.4) Let $G_i$ be a group (ring, algebra, et cetera) for all $i\in I$. A subset $S \subset \varprojlim G_i$  is a set of {\it topological generators} of $\varprojlim G_i$ if the span $\langle S\rangle$ is dense in $\varprojlim G_i$. If moreover
$S$ is finite then $\varprojlim G_i$ is  said to be {\it finitely generated}.
\end{defn}
Our method for finding approximating sequences of $p$--adic elements is by strict inclusions of open neighborhoods. As a topological space, $\prod X_i$ is endowed with the product topology, so $\varprojlim X_i$ inherits the induced topology. It can be then verified  that $\varprojlim X_i$ is closed in  $\prod X_i$. A  basis of open sets in $\varprojlim X_i$  contains elements of the form
$$
  \varpi_{i}^{-1}(U_i)  \cap \varprojlim X_i
$$
where  $U_i$ open in $X_i$. Then, any open set in $\varprojlim X_i$ is a union of sets of the form
$$
\varpi_{i_1}^{-1}(U_1) \cap \ldots \cap \varpi_{i_n}^{-1}(U_n)  \cap \varprojlim X_i
$$
where $i_1, \ldots,i_n
\in I$ and $U_r$ open in $X_{i_r}$ for each $r$. (Cf. \cite{riza}, p.7.)

\subsection{\it  Approximations in ${\Bbb Z}_p$}

Since ${\Bbb Z}$ projects onto each factor ${\Bbb Z} / p^r{\Bbb Z}$, by Lemma~\ref{dense}, the image of ${\Bbb Z}$ under the induced map on the inverse limits is dense in ${\Bbb Z}_p$. Now $\varprojlim_{r} {\Bbb Z} = {\Bbb Z}$ and the induced map acts on an element $(x, x, \ldots) \in {\Bbb Z}$ by sending $x$ to $x\,(\text{mod}\, p^r) \in {\Bbb Z} / p^r{\Bbb Z}$ for every $r$. But, after some point $x$ will be unchanged by the modulus, so $\left(x\,(\text{mod}\, p), x\,(\text{mod}\, p^2), \ldots\right) = (x, x, \ldots)$. Therefore, the image of ${\Bbb Z}$ under the induced map on the inverse limits is ${\Bbb Z}$, and so ${\Bbb Z}$ is dense  in ${\Bbb Z}_p$.

 Further, ${\Bbb Z} = \langle {\bold t}\rangle$ so ${\bold t}$ is a  topological generator of ${\Bbb Z}_p$. Thus, an element ${\bold t}^{\underleftarrow{a}} =
(t_1^{a_1}, t_2^{a_2}, \ldots) $ in ${\Bbb Z}_p$ is approximated
by constant sequences, which are identified with integers.
We shall explain how to find such an approximating sequence for a $p$--adic
integer, in order to draw the strategy for the larger $p$--adic
structures we are dealing with.

The inherited topology of ${\Bbb Z}_p$ builds up from the discrete topology of each factor ${\Bbb Z} / p^r{\Bbb Z}$.
 Thus, a basic open set $U$ in ${\Bbb Z}_p$ is of the form:
 $U=\varpi_i^{-1} (U_i) \ ; \ U_i \subseteq {\Bbb Z} / p^i{\Bbb Z}$. For $U_i$  not a singleton,
$\varpi_i^{-1} (U_i) = \cup_{u \in U_i}\varpi_i^{-1} (\{u\})$.

Recall now Definition~\ref{unique}. It is then easy to verify the lemma below.
\begin{lem}\label{nbds}
Let $\underleftarrow{a} = (a_1, a_2, a_3, \ldots) \in {\Bbb Z}_p$ in reduced form
and let $U_i \subseteq {\Bbb Z} / p^i{\Bbb Z}$ for some $i$. Then
$\underleftarrow{a} \in U=\varpi_i^{-1} (U_i)$ if and only if $a_i
\in U_i$. Hence, a basic open neighborhood of $\underleftarrow{a}$ in ${\Bbb Z}_p$
 is of the form
$U=\varpi_i^{-1}(\{a_i\})$ for some $i$. Moreover, we have a nested sequence of neighborhoods with strict inclusions:
\[
\varpi_1^{-1}(\{ a_1 \}) \supsetneq \varpi_2^{-1}(\{a_2 \})
\supsetneq \cdots. \]
\end{lem}

By the strict inclusions of neighborhoods, the sequence of constant sequences $((a_k))_{k\in {\Bbb N}}$  in ${\Bbb Z}$ approximates $\underleftarrow{a} \in {\Bbb Z}_p$ and we write
$\underleftarrow{a} = \lim_{k} (a_k)$, or, in the multiplicative notation:
$$
{\bold t}^{\underleftarrow{a}} = \lim_{k} {\bold t}^{a_k}
$$
Indeed, subtracting each constant sequence
successively from $\underleftarrow{a}$ the
differences tend to the zero sequence:
$$
\begin{array}{ccl}
(a_1, a_2, a_3, \ldots) - (a_1, a_1, a_1, \ldots) & = &  (0, a_2-a_1, a_3-a_1, \ldots)  \\
 (a_1, a_2, a_3, \ldots) - (a_2, a_2, a_2, \ldots) &  = &  (0,  0 , a_3-a_2, \ldots) \\
   \, & \vdots & \, \\
\end{array}
$$

 In order to reach a general scheme for finding approximating sequences for purely $p$--adic elements we shall introduce the operations truncation and expansion for entries of $p$--adic integers.

\begin{defn}\rm
Let $\underleftarrow{a} = (a_1, a_2, \ldots) \in {\Bbb Z}_p$ in reduced form. For any indices $r,s$ with $r\geq s$ we define the {\it $s$--truncation of $a_r$} as the element
$$
a_{rs} = k_0 + k_1 p + \cdots + k_{s-1} p^{s-1} + p^r{\Bbb Z} \ \in {\Bbb Z} / p^r{\Bbb Z}
$$
Note that $\theta_s^r(a_{rs}) = a_s$, that is, the elements $a_{rs}$ and $a_s$ are coherent via the map (\ref{thetars}). Note also that $a_{rr} = a_r$. Similarly, we define the {\it $r$--expansion of $a_s$} as the element
$$
a_{sr} = k_0 + k_1 p + \cdots + k_{r-1} p^{r-1} + p^s{\Bbb Z} \ \in  {\Bbb Z} / p^s{\Bbb Z}
$$
Note that $a_{sr} \equiv a_s ({ \rm \text mod \, } p^s)$, so $\theta_s^r(a_r) = a_{sr}$.
 In the multiplicative notation $t_r^{a_r}$ is substituted by $t_r^{a_s}$ in the first case and $t_s^{a_s}$ is substituted by $t_s^{a_r}$ in the second case. The fact that $\underleftarrow{a}$ is in reduced form ensures that truncations and expansions of its entries are well--defined.
\end{defn}

In the above terminology, the constant sequences $(a_k)\in {\Bbb Z}$ approximating
$\underleftarrow{a}$ are found from $\underleftarrow{a}$ by truncating each term  after $a_k$
with respect to $a_k$ and expanding each term before $a_k$ with respect to $a_k$.

\subsection{\it  Approximations in ${\mathcal F}_{p^{\infty},n}$ and ${\Bbb C}{\mathcal F}_{p^{\infty} ,n}$}

Applying the canonical epimorphism (\ref{varthetars}) componentwise yields a canonical epimorphism of ${\Bbb Z}^n$ on each factor $ C_r^n$. So, by Lemma~\ref{dense}, ${\Bbb Z}^n$ is dense in ${\Bbb Z}_p^n$. Then, for example, for $\underleftarrow{a} = (a_1, a_2, \ldots)$ and
$\underleftarrow{b} = (b_1,b_2, \ldots) \in {\Bbb Z}_p$ in reduced form, the
element $(\underleftarrow{a}, \underleftarrow{b}) \in {\Bbb
Z}_p^2$ is approximated by the sequence $((a_k,b_k))_{k\in {\Bbb N}}$ of constant sequences, that
is, with terms in the dense subgroup ${\Bbb Z}^2$. In multiplicative notation:
$({\bold t}^{\underleftarrow{a}},{\bold t}^{\underleftarrow{b}}) =
{\bold t}_1^{\underleftarrow{a}}{\bold t}_2^{\underleftarrow{b}}
\in {\Bbb Z}_p^2 \ $ is approximated by the sequence
 $(({\bold t}^{a_k},{\bold t}^{b_k}))_{k\in {\Bbb N}} = (({\bold t}_1^{a_k}{\bold
t}_2^{b_k}))_{k\in {\Bbb N}} \in {\Bbb Z}^2$.

\begin{figure}[H]
\begin{picture}(220,95)

\qbezier(0,15)(-10,50)(0,85) 
\qbezier(10,70)(9,75)(10,80)
\qbezier(30,70)(31,75)(30,80)
\qbezier(10,20)(9,25)(10,30)
\qbezier(30,20)(31,25)(30,30)

\qbezier(10,30)(10,35)(20,40)
\qbezier(20,40)(30,45)(30,50)
\qbezier(30,30)(30,34)(25,37)
\qbezier(15,43)(10,46)(10,50)
\qbezier(10,50)(10,55)(20,60)
\qbezier(20,60)(30,65)(30,70)
\qbezier(30,50)(30,54)(25,57)
\qbezier(15,63)(10,66)(10,70)

\put(40,50){,}

\qbezier(50,70)(49,75)(50,80)
\qbezier(70,70)(71,75)(70,80)
\qbezier(50,20)(49,25)(50,30)
\qbezier(70,20)(71,25)(70,30)

\qbezier(50,30)(50,35)(60,40)
\qbezier(60,40)(70,45)(70,50)
\qbezier(70,30)(70,34)(65,37)
\qbezier(55,43)(50,46)(50,50)
\qbezier(50,50)(50,55)(60,60)
\qbezier(60,60)(70,65)(70,70)
\qbezier(70,50)(70,54)(65,57)
\qbezier(55,63)(50,66)(50,70)

\put(80,50){,}

\qbezier(90,70)(89,75)(90,80)
\qbezier(110,70)(111,75)(110,80)
\qbezier(90,20)(89,25)(90,30)
\qbezier(110,20)(111,25)(110,30)

\qbezier(90,30)(90,35)(100,40)
\qbezier(100,40)(110,45)(110,50)
\qbezier(110,30)(110,34)(105,37)
\qbezier(95,43)(90,46)(90,50)
\qbezier(90,50)(90,55)(100,60)
\qbezier(100,60)(110,65)(110,70)
\qbezier(110,50)(110,54)(105,57)
\qbezier(95,63)(90,66)(90,70)

\put(120,50){,}
\put(130, 50){$\ldots$}

\qbezier(145,15)(155,50)(145,85) 

\put(160,50){$=$}
\put(175,50){$\lim$}
\put(182,44){{\tiny $r$}}

\qbezier(200,70)(199,75)(200,80)
\qbezier(220,70)(221,75)(220,80)
\qbezier(200,20)(199,25)(200,30)
\qbezier(220,20)(221,25)(220,30)

\qbezier(200,30)(200,35)(210,40)
\qbezier(210,40)(220,45)(220,50)
\qbezier(220,30)(220,34)(215,37)
\qbezier(205,43)(200,46)(200,50)
\qbezier(200,50)(200,55)(210,60)
\qbezier(210,60)(220,65)(220,70)
\qbezier(220,50)(220,54)(215,57)
\qbezier(205,63)(200,66)(200,70)

\put(7,85){{\tiny $a_1$}}
\put(27,85){{\tiny $b_1$}}

\put(47,85){{\tiny $a_2$}}
\put(67,85){{\tiny $b_2$}}

\put(87,85){{\tiny $a_3$}}
\put(107,85){{\tiny $b_3$}}

\put(195,85){$a_r$}
\put(216,85){$b_r$}
\end{picture}
\caption{Approximating a $p$--adic braid by classical braids}
\label{brappx}
\end{figure}
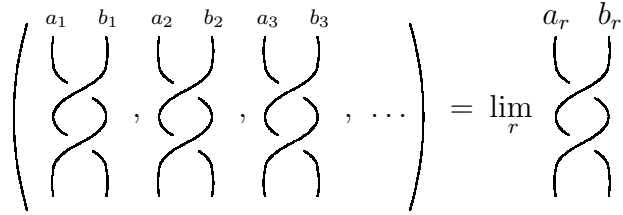

We then extend the projection of ${\Bbb Z}^n$ on each factor $ C_r^n$ by the
identity map on $B_n$. So, we obtain an epimorphism of the classical framed braid group
${\mathcal F}_n = {\Bbb Z}^n \rtimes B_n$ on each factor ${\mathcal F}_{p^r ,n}$.
Hence, by Lemma~\ref{dense}, ${\mathcal F}_n$ is dense in ${\mathcal F}_{p^{\infty} ,n}$. The set
$\{{\bold t}_1,\ldots, {\bold t}_n, {\sigma}_1, \ldots, {\sigma}_{n-1}\} $
is a set of topological  generators for ${\mathcal F}_{p^{\infty} ,n}$ satisfying relations analogous to the relations of ${\mathcal F}_n$.  Moreover, for $p$--adic integers  $\underleftarrow{a_i} = (a_{ri})_r$ in reduced form,
 an element $\underleftarrow{\beta} = {\bold t}_1^{\underleftarrow{a_1}} \ldots {\bold
t}_n^{\underleftarrow{a_n}} \cdot \sigma \in {\mathcal F}_{p^{\infty} ,n}$ has the approximation:
$$
\underleftarrow{\beta} =  \lim_k ({\bold t}_1^{a_{k1}} \ldots {\bold t}_n^{a_{kn}} \cdot \sigma)
$$
with ${\bold t}_1^{a_{k1}} \ldots {\bold t}_n^{a_{kn}} \cdot \sigma  \in {\mathcal F}_n$. An example is illustrated in Figure~\ref{brappx}.

\smallbreak

Passing to algebras, for an element in ${\Bbb C}{\mathcal F}_{p^{\infty} ,n}$ in reduced form (recall Definition~\ref{preduced}) it is easy to find an approximating sequence with elements in the dense algebra ${\Bbb C}{\mathcal F}_n$, since ${\mathcal F}_n$ is dense in ${\mathcal F}_{p^{\infty} ,n}$. Indeed, we simply extend linearly the approximations of its monomials in ${\mathcal F}_{p^{\infty},n}$, as described above. Note that, by construction, our $p$--adic element may always be written in the form where the linear combinations in each place have the same number of coherent monomials with the same coefficients. So, we have the following:

\begin{prop}\label{simpledense}  The algebra ${\Bbb C}{\mathcal F}_n$ is dense in ${\Bbb C}{\mathcal F}_{p^{\infty} ,n}$.
\end{prop}

 We would like now to find approximating sequences for purely $p$--adic elements in $\varprojlim_{r} {\Bbb C}{C_r^n}$, in $\varprojlim_{r} {\Bbb C}{\mathcal F}_{p^r,n}$ and in ${\rm Y}_{p^{\infty}, n}(u)$. The tactics used for ${\Bbb C}{\mathcal F}_{p^{\infty} ,n}$, that is, approximating each monomial cannot  be applied here for purely $p$--adic elements (such as $e_{i}$) since they cannot be written in the form where the linear combinations in each place have the same number of coherent terms with the same coefficients. So, finding an approximating sequence for purely $p$--adic elements is more tricky. In any case, we need first to find dense subalgebras of constant elements, in which the approximating terms should live.

\subsection{\it  Dense subsets in the $p$--adic algebras}

 Extending linearly the epimorphism of ${\mathcal F}_n$ on each factor ${\mathcal F}_{p^r ,n}$ defines an epimorphism $\eta_r$ of the algebra ${\Bbb C}{\mathcal F}_n$ on the algebra ${\Bbb C}{\mathcal F}_{p^r ,n}$. Moreover, the map $\eta_r$ composed with the canonical epimorphism $\rho_r$ defines an epimorphism $\mu_r := \rho_r \circ \eta_r$ of ${\Bbb C}{\mathcal F}_n$ on the algebra ${\rm Y}_{p^r, n}(u)$:
\begin{equation}\label{algepis}
\begin{array}{rcccc}
{\Bbb C}{\mathcal F}_n  &  \stackrel{\eta_r}{\longrightarrow} &  {\Bbb C}{\mathcal F}_{p^r,n}  & \stackrel{\rho_r}{\longrightarrow} &  {\rm Y}_{p^r, n}(u) \\
 \sigma_i   &  \mapsto & \sigma_i &  \mapsto & g_i \\
 t_{j}^m   & \mapsto & t_{r,j}^{ m \,( {\rm mod}\, p^r)} & \mapsto & t_{r,j}^{ m \,( {\rm mod}\, p^r)}
\end{array}
\end{equation}

Define further $\eta:=\varprojlim_{r} {\eta_r}, \ \rho:=\varprojlim_{r} {\rho_r}$ and $\mu:=\varprojlim_{r} {\mu_r}$, the corresponding induced maps on the inverse limits of the maps (\ref{algepis}). That is:
\begin{equation}\label{eta}
\begin{array}{ccc}
\eta : \ {\Bbb C}{\mathcal F}_n & \longrightarrow & \varprojlim_{r} {\Bbb C}{\mathcal F}_{p^r,n}
\end{array}
\end{equation}
\begin{equation}\label{rho}
\rho : \ \varprojlim_{r} {\Bbb C}{\mathcal F}_{p^r,n} \longrightarrow {\rm Y}_{p^{\infty}, n}(u)
\end{equation}
\begin{equation}\label{mu}
\begin{array}{ccc}
\mu : \ {\Bbb C}{\mathcal F}_n  & \longrightarrow &  {\rm Y}_{p^{\infty}, n}(u)
\end{array}
\end{equation}

Recall that $\varprojlim_{r} {\Bbb C}{\mathcal F}_n = {\Bbb C}{\mathcal F}_n$. One can easily check that $\eta(\sigma_i)=\sigma_i$ and $\eta({\bold t}_j)={\bold t}_j$, so $\eta$ is an injection. Also that $\rho({\sigma}_i) = g_i$  and $\rho({\bold t}_j) = {\bold t}_j$. Below, in Proposition~\ref{rhosurj}, we also show that $\rho$ is a surjection. Moreover, by the construction of the maps (\ref{algepis}) it follows that $\mu = \rho\circ \eta $ and that:
$$
\mu({\sigma}_i) = g_i \mbox{ and } \mu({\bold t}_j) = {\bold t}_j.
$$

\begin{notation}\label{denseyn} \rm
 We denote the subalgebra $\mu({\Bbb C}{\mathcal F}_n)$ of ${\rm Y}_{p^{\infty},n}(u)$ as:
 $$
 {\rm Y}_n(u) := \mu({\Bbb C}{\mathcal F}_n)
 $$
\end{notation}

In \cite{jula1} we had found topological generators for $\varprojlim_{r}C_r^n$ (Lemma~2\cite{jula1}), for ${\mathcal F}_{p^{\infty} ,n}$ (Theorem~1\cite{jula1}), for $\varprojlim_{r} {\Bbb C}{\mathcal F}_{p^r,n}$ (Proposition~6\cite{jula1}) and for ${\rm Y}_{p^{\infty}, n}(u)$ (Theorem~3\cite{jula1}).
 We prove here the following theorem.

\begin{thm}\label{densesets} The following hold:

\vspace{.15cm}
\noindent (i) The algebra ${\Bbb C}{\mathcal F}_n$ is dense in $\varprojlim_{r} {\Bbb C}{\mathcal F}_{p^r,n}$.

\vspace{.15cm}
\noindent (ii) The algebra ${\Bbb C}{\mathcal F}_{p^{\infty},n}$ is dense in $\varprojlim_{r} {\Bbb C}{\mathcal F}_{p^r,n}$.

\vspace{.15cm}
\noindent (iii) The set $X=\{{\bold 1}, {\bold t}_1,\ldots, {\bold t}_n, {\sigma}_1, \ldots, {\sigma}_{n-1}\}$ is a set of topological generators for the algebras ${\Bbb C}{\mathcal F}_{p^{\infty},n}$ and $\varprojlim_{r} {\Bbb C}{\mathcal F}_{p^r,n}$. Together with the braid relations for the ${\sigma}_i$'s, the commuting relations for the ${\bold t}_j$'s and the relations:
\[{\sigma}_i {\bold t}_i = {\bold t}_{i+1} {\sigma}_i , \
{\sigma}_i {\bold t}_{i+1} = {\bold t}_i {\sigma}_i \ {\text and } \
 {\sigma}_i {\bold t}_j = {\bold t}_j {\sigma}_i \ {\text for } \ j\not= i, i+1
\]
they furnish a topological presentation for ${\Bbb C}{\mathcal F}_{p^{\infty},n}$ and  $\varprojlim_{r} {\Bbb C}{\mathcal F}_{p^r,n}$.

\vspace{.15cm}
\noindent (iv) The algebra ${\rm Y}_n(u)=\mu({\Bbb C}{\mathcal F}_n)$ is dense in ${\rm Y}_{p^{\infty}, n}(u)$. Moreover, the set:
$$D=\{{\bold 1},{\bold t}_1,\ldots,{\bold t}_n, g_1,\ldots,g_{n-1}\}$$
is a set of topological generators for the algebra ${\rm Y}_{p^{\infty}, n}(u)$, satisfying the analogous relations of (iii) in ${\rm Y}_n(u)$.

\vspace{.15cm}
\noindent (v) The relations ${\bold t}_j e_{i} = e_{i} {\bold t}_j$ (see end of Subsection~\ref{eltsei}) are valid in $\varprojlim_{r} {\Bbb C}{\mathcal F}_{p^r,n}$ and in $ {\rm Y}_{p^{\infty} ,n}(u)$, but not in the dense subalgebras. The same is true about the quadratic relations (\ref{pquadr}) that are valid in ${\rm Y}_{p^{\infty}, n}(u)$.
\end{thm}

\begin{proof}
 Claim (i) is an application of the surjections (\ref{algepis}), Lemma~\ref{dense} and of the following observation:  after some point the image of the exponent $m$ in (\ref{algepis}) will not change, so $\eta ({\Bbb C}{\mathcal F}_n) = {\Bbb C}{\mathcal F}_n$. Then, since $\eta_r$ is a surjection, it follows by Lemma~\ref{dense} that the algebra ${\Bbb C}{\mathcal F}_n$ is dense in $\varprojlim_{r} {\Bbb C}{\mathcal F}_{p^r,n}$. Claim (ii) follows immediately from (i), Proposition~\ref{algebras} and Proposition~\ref{simpledense}.

It follows now from (i) that the set $X$ is a set of topological generators, satisfying the listed relations. Moreover, by the standard presentation of the classical framed braid group ${\mathcal F}_n$ (recall (\ref{action})), the relations given in claim (iii) are the only ones satisfied in the set $X$.

The fact that ${\rm Y}_n(u)=\mu({\Bbb C}{\mathcal F}_n)$ is dense in ${\rm Y}_{p^{\infty}, n}(u)$ is clear by a direct application of the surjections (\ref{algepis}) and Lemma~\ref{dense}. Further, $\mu({\bold t}_j) = {\bold t}_j$ and $\mu({\sigma}_i) = g_i$. So, claim (iv) follows. Finally, claim (v) follows from claims (iii) and (iv), from Corollary~\ref{einotin} and from the fact that $\mu(e_{i}) = e_{i}$.
\end{proof}

Focusing now a little more on ${\rm Y}_{p^{\infty}, n}(u)$, an apparently dense subset in ${\rm Y}_{p^{\infty}, n}(u)$, discussed in \cite{jula1}, comes from the following construction. For any $r$ we have the following exact sequence:
\begin{equation}\label{exact}
\begin{CD}
0 @>>>I_{p^r, n}  @>\iota_r>> {\Bbb C}{\mathcal F}_{p^r,n} @>\rho_r>>{\rm Y}_{p^r ,n}(u)@>>> 0
\end{CD}
\end{equation}
where $I_{p^r, n}$ is the ideal generated by the linear expressions in Eq.~\ref{quadr}. This induces the exact sequence:
$$
\begin{CD}
0 @>>>\varprojlim_{r} I_{p^r, n}  @>\iota>>  \varprojlim_{r}{\Bbb C}{\mathcal F}_{p^r,n}
 @>\rho>>  {\rm Y}_{p^{\infty} ,n}(u)
\end{CD}
$$
where $\iota:= \varprojlim_{r} \iota_r$ and $\rho:= \varprojlim_{r}\rho_r$. Hence,
 and since $\varprojlim_{r} I_{p^r, n}$ is an ideal in $\varprojlim_{r}{\Bbb C}{\mathcal F}_{p^r,n}$,
we have:
$$
\frac{\varprojlim_{r}{\Bbb C}{\mathcal F}_{p^r,n}}{\varprojlim_{r}I_{p^r, n}} \cong
\rho ( \varprojlim_{r}{\Bbb C}{\mathcal F}_{p^r,n})
$$
At the time of writing of \cite{jula1} it was not clear whether the map $\rho$ is a surjection or not. Yet, by application of Lemma~\ref{dense} we were able to derive the result that $\rho ( \varprojlim_{r}{\Bbb C}{\mathcal F}_{p^r,n})$ is dense in $ {\rm Y}_{p^{\infty} ,n}(u)$ (Proposition~8\cite{jula1}). We are now in a position to prove surjection for $\rho$. Before that we need to recall the following definition. An inverse system $X = (X_i, \phi_j^i)$ indexed by ${\Bbb N}$ is said to satisfy the {\it ML--condition} (Mittag--Leffler condition) if for any index $m$ there exists $n\geq m$ such that for all $n^\prime\geq n$ we have ${\rm \text Im}(\phi_m^n) = {\rm \text Im}(\phi_m^{n^\prime})$. Note that if all $\phi_j^i$ are surjective then $X$ satisfies the ML--condition.

\begin{prop}\label{rhosurj}
The map $\rho$ of Eq.~\ref{rho} is a surjection. Hence:
$$
{\rm Y}_{p^{\infty} ,n}(u) \cong \frac{\varprojlim_{r}{\Bbb C}{\mathcal F}_{p^r,n}}{\varprojlim_{r} I_{p^r, n}}.
$$
\end{prop}
\begin{proof}
The exact sequence (\ref{exact}) induces the following exact sequence of inverse systems.
$$
\begin{CD}
0 @>>>(I_{p^r, n},\varphi ^{r}_s)  @>(\iota_r)>> ({\Bbb C}{\mathcal F}_{p^r,n},\varphi ^{r}_s) @>(\rho_r)>>({\rm Y}_{p^r ,n},\phi_s^r)@>>> 0
\end{CD}
$$
Now, by Lemma~6\cite{jula1}, $\varphi ^{r}_s(I_{p^r,n})=  I_{p^s,n}$ for any $r,s$ with $r\geq s$. Hence the inverse system $(I_{p^r, n},\varphi ^{r}_s)$ satisfies the ML--condition. Then, by a well--known result of Grothendieck the following exact sequence is induced:
$$
\begin{CD}
0 @>>>\varprojlim_{r} I_{p^r, n}  @>\iota>>  \varprojlim_{r}{\Bbb C}{\mathcal F}_{p^r,n}
 @>\rho>>  {\rm Y}_{p^{\infty} ,n}(u)@>>> 0
\end{CD}
$$
Hence $\rho$ is surjection.
\end{proof}

We shall now recapitulate what we know about our $p$--adic objects, by means of  a concise diagram. For that we need to introduce three more intermediate structures.

\begin{defn}\label{densextn}\rm
\noindent (i) We define the dense subalgebra $\widetilde{{\Bbb C}{\mathcal F}_n}$ of $\varprojlim_{r}{\Bbb C}{\mathcal F}_{p^r,n}$ as the extension of the subalgebra ${\Bbb C}{\mathcal F}_n$ by the elements $e_{i}$:
$$
\widetilde{{\Bbb C}{\mathcal F}_n} :=\langle {\Bbb C}{\mathcal F}_{n},e_{1},\ldots, e_{n-1}\rangle
$$
(ii) We define the dense subalgebra
$\widetilde{{\Bbb C}{\mathcal F}_{p^{\infty} ,n}}$ of $\varprojlim_{r}{\Bbb C}{\mathcal F}_{p^r,n}$ as the extension of the subalgebra ${\Bbb C}{\mathcal F}_{p^{\infty} ,n}$ by the elements $e_{i}$:
$$
\widetilde{{\Bbb C}{\mathcal F}_{p^{\infty} ,n}} :=\langle {\Bbb C}{\mathcal F}_{p^{\infty} ,n}, e_{1},\ldots, e_{n-1}\rangle
$$
(iii) We define the dense subalgebra $\widetilde{{\rm Y}_n}(u)$ of ${\rm Y}_{p^{\infty},n}(u)$ as the extension of ${\rm Y}_n(u) = \mu({\Bbb C}{\mathcal F}_n)$ by the elements
$e_{i}$:
$$
{\widetilde{{\rm Y}_n}(u)} :=\langle {\rm Y}_n(u),e_{1},\ldots, e_{n-1} \rangle
$$
\end{defn}

Clearly $\widetilde{{\Bbb C}{\mathcal F}_n}$ is a proper subset of $\widetilde{{\Bbb C}{\mathcal F}_{p^{\infty} ,n}}$ and of $\varprojlim_{r}{\Bbb C}{\mathcal F}_{p^r,n}$. For example, $\widetilde{{\Bbb C}{\mathcal F}_n}$ does not contain the $p$--adic integers. By the same reason $\widetilde{{\rm Y}_n}(u)$ is also a proper subset of ${\rm Y}_{p^{\infty},n}(u)$. Also, by  Proposition~\ref{algebras}, $\widetilde{{\Bbb C}{\mathcal F}_{p^{\infty} ,n}}$ injects in  $\varprojlim_{r}{\Bbb C}{\mathcal F}_{p^r,n}$.
 Moreover, if we denote $\widetilde{\mu}$ the extension of $\mu$ on $\widetilde{{\Bbb C}{\mathcal F}_n}$ by defining $\widetilde{\mu}(e_{i}) = e_{i}$, we have that $\widetilde{\mu} (\widetilde{{\Bbb C}{\mathcal F}_n}) = \widetilde{{\rm Y}_n}(u)$. So, we have the following commuting diagram:
$$
\begin{diagram}
\node{e_{i}\notin  {\Bbb C}{\mathcal F}_{n}}
\arrow{e,t}{\eta =\iota}
\arrow{s,r}{\mu}
\node{\widetilde{{\Bbb C}{\mathcal F}_{n}} }
\arrow{e,tb}{\iota}{\not\subseteq}
\arrow{s,r}{\widetilde{\mu}}
\node{\widetilde{{\Bbb C}{\mathcal F}_{p^{\infty} ,n}}}
\arrow{e,tb}{\iota}{\not\subseteq}
\node{\varprojlim_{r} {\Bbb C}{\mathcal F}_{p^r,n}}
\arrow{sw,b}{\rho}\\
\node{e_{i}\notin Y_n(u)}
\arrow{e,t}{\iota}
\node{ \widetilde{ {\rm Y}_{n}(u)} }
\arrow{e,tb}{\iota}{\not\subseteq}
\node{{\rm Y}_{p^{\infty},n}(u)}
\node{}
\end{diagram}
$$

Further, the quadratic relations (\ref{pquadr}) are valid in $\widetilde{{\rm Y}_n}(u)$ and, by construction,  $\widetilde{{\rm Y}_n}(u)$ is the smallest subalgebra of  ${\rm Y}_{p^{\infty},n}(u)$ which is closed under the quadratic relations (\ref{pquadr}). The relations in Theorem~\ref{densesets}(iv) together with the quadratic relations (\ref{pquadr}) form then a complete set of relations for $\widetilde{{\rm Y}_n}(u)$. Thus, we have  the following:

\begin{thm}\label{denseyhtilde}
The dense subalgebra $\widetilde{{\rm Y}_n}(u)$ in Definition~\ref{densextn} can be viewed as the quotient:
\[
\widetilde{{\rm Y}_n}(u) =\frac{\widetilde{{\Bbb C}{\mathcal F}_n}}
{\langle g_i^2 - 1 - (u-1)e_{i} - (u-1)e_{i}\, g_i \rangle}
\]
and, so, ${\rm Y}_{p^{\infty},n}(u)$ can be regarded as a topological deformation of the above quotient algebra.
\end{thm}

With our dense subalgebras in hand, we shall next discuss an approximation for the elements $e_{i}$ before going to the general case of approximating purely $p$--adic elements.

\subsection{\it  Approximating $e_{i}$ }

Recall (Eq.~\ref{epri}) that an entry $e_{p^r,i}$ of the element  $e_{i}= (e_{p,i}, e_{p^2,i}, \ldots )  \in \varprojlim_{r} {\Bbb C}{C_r^n} \subset \varprojlim_{r} {\Bbb C}{\mathcal F} _{p^r,n}$ has $p^r$ terms with linear coefficients $\frac{1}{p^r}$ and this is in reduced form, according to Definition~\ref{preduced}. We are looking for an approximating sequence for $e_{i}$, consisting of constant terms. Recall from the proof of Lemma~\ref{eprtoeps} that $e_{p^r,i}$ can be arranged in the form:
\begin{equation}\label{trunc}
e_{p^r,i} =\frac{1}{p^{r-s}} \left( \sum_{m=0}^{p^s-1}
\frac{1}{p^s} t_{r,i}^m t_{r,i+1}^{-m} + \sum_{m=p^s}^{2p^s-1}
\frac{1}{p^s} t_{r,i}^m t_{r,i+1}^{-m} +\cdots+
\sum_{m=p^r-p^s}^{p^r-1} \frac{1}{p^s} t_{r,i}^m t_{r,i+1}^{-m}
\right)
\end{equation}
 of $p^{r-s}$ packets, each of which projects on $e_{p^s,i}$ by the coherency map
 $\varphi^{r}_s$  $(r\geq s)$. That is, $e_{p^r,i}$ `wraps' $p^{r-s}$ times on $e_{p^s,i}$ by the  map $\varphi^{r}_s$. We shall define the first packet as the $s$--truncation of $e_{p^r,i}$.

On the other hand, in order to make the linear expression for $e_{p^s,i}$ agree formally with that for $e_{p^r,i}$,  we rewrite $e_{p^s,i}$ as a sum of a total of $p^r$ terms, arranged in $p^{r-s}$ packets, each of which is equal to $e_{p^s,i}$:
\begin{equation}\label{expan}
e_{p^s,i} =\frac{1}{p^{r-s}} \left( \sum_{m=0}^{p^s-1}
\frac{1}{p^s} t_{s,i}^m t_{s,i+1}^{-m} + \sum_{m=p^s}^{2p^s-1}
\frac{1}{p^s} t_{s,i}^m t_{s,i+1}^{-m} +\cdots+
\sum_{m=p^r-p^s}^{p^r-1} \frac{1}{p^s} t_{s,i}^m t_{s,i+1}^{-m}
\right)
\end{equation}
 We shall define this as the $r$--expansion of $e_{p^s,i}$.

\begin{defn}\label{truncexp} \rm
For any indices $r,s$ with $r\geq s$ we define the element $\zeta_{r,s,i}$ in ${\Bbb C}{C_r^n} \subset  {\Bbb C}{\mathcal F} _{p^r,n}$ as the formal expression of $e_{p^s,i}$, but with the generators $t_{s,i}$ and $t_{s,i+1}$ replaced by the generators $t_{r,i}$ and $t_{r,i+1}$ respectively. That is:
$$
\zeta_{r,s,i} :=\frac{1}{p^s}\sum_{m=0}^{p^s-1}t_{r,i}^m t_{r, i+1}^{-m} \  \in {\Bbb C}{\mathcal F} _{p^r,n}
$$
The element $\zeta_{r,s,i}$ is called  {\it the $s$--truncation of $e_{p^r,i}$}. Note that $\zeta_{r,r,i} = e_{{p^r},i}$. Clearly, the elements $\zeta_{r,s,i}$ and $e_{p^s,i}$ are coherent: $\varphi ^{r}_s(\zeta_{r,s,i}) = e_{p^s,i}$. In fact, $\zeta_{r,s,i}$ `wraps' only once on $e_{p^s,i}$ via the map $\varphi^{r}_s$.

Further, we define the element $\zeta_{s,r,i}$ in ${\Bbb C}{C_s^n} \subset {\Bbb C}{\mathcal F} _{p^s,n}$ as the formal expression of $e_{p^r,i}$, but with the generators $t_{r,i}$ and
 $t_{r,i+1}$ replaced by the generators $t_{s,i}$ and $t_{s,i+1}$ respectively. That is:
$$
\zeta_{s,r,i} :=\frac{1}{p^r}\sum_{m=0}^{p^r-1}t_{s,i}^m t_{s, i+1}^{-m} \  \in {\Bbb C}{\mathcal F} _{p^s,n}
$$
The element $\zeta_{s,r,i}$ is called  {\it the $r$--expansion of $e_{p^s,i}$}. Clearly, $\zeta_{s,s,i} = e_{p^s,i}$. In fact, $e_{p^r,i}$ `wraps' only once on $\zeta_{s,r,i}$ via the map $\varphi ^{r}_s$.
\end{defn}

We define now the element ${\bf e}_{p^r,i}$ by $r$--truncating each term in $e_{i}$ after the $r$th position and by $r$--expanding each term before the $r$th position. That is,
$$
{\bf e}_{p^r,i} := (\zeta_{1,r,i}, \zeta_{2,r,i}, \ldots,\zeta_{r-1,r,i}, e_{p^r,i}, \zeta_{r+1,r,i}, \zeta_{r+2,r,i}, \ldots)
$$

\begin{prop}
For any index $r$ the element ${\bf e}_{p^r,i}$ is a constant sequence in $\varprojlim_{r} {\Bbb C}{C_r^n} \subset\varprojlim_{r} {\Bbb C}{\mathcal F} _{p^r,n}$. More precisely:
$$
{\bf e}_{p^r,i} = (\frac{1}{p^r}\sum_{m=0}^{p^r-1}t_{k,i}^m t_{k,i+1}^{-m})_k  = \frac{1}{p^r}\sum_{m=0}^{p^r-1}{\bf
t}_i^m{\bf t}_{ i+1}^{-m} \ \in {\Bbb C}{\mathcal F}_{n}
$$
Moreover, we have the approximation:
$$
e_{i} = \lim_{r} {\bf e}_{p^r,i}
$$
\end{prop}
\begin{proof}
By Definition~\ref{truncexp} the sequence ${\bf e}_{p^r,i}$ is coherent, so ${\bf e}_{p^r,i} \in \varprojlim_{r} {\Bbb C}{C_r^n} \subset \varprojlim_{r} {\Bbb C}{\mathcal F} _{p^r,n}$. Moreover, all terms in ${\bf e}_{p^r,i}$ have the same formal expression, that of $e_{p^r,i}$, so ${\bf e}_{p^r,i}$ is the constant sequence given in the statement. Finally, recall from (\ref{isoms}) that $(t_{r,i})_r  = {\bold t}_i \in { \Bbb Z}_p^n$. So, separating terms in ${\bf e}_{p^r,i}$ we obtain:
$$
{\bf e}_{p^r,i} = \frac{1}{p^r} \sum_{m=0}^{p^r-1} \left(t_{1,i}^m, t_{2,i}^m, \ldots\right)\left(t_{1,i+1}^{-m},  t_{2,i+1}^{-m}, \ldots\right) = \frac{1}{p^r} \sum_{m=0}^{p^r-1} {\bf t}_i^m{\bf t}_{ i+1}^{-m}
$$
Subtracting, now, from $e_{i}$ each element of the sequence $({\bf e}_{p^r,i})_r$ successively, we obtain:
\[
\begin{array}{l}
(e_{p,i}, e_{p^2,i}, e_{p^3,i}, \ldots ) - \ (e_{p,i},\zeta_{2,1,i},\zeta_{3,1,i},\ldots) = (0, e_{p^2,i}-\zeta_{2,1,i}, e_{p^3,i}-\zeta_{3,1,i}, \ldots)  \\
\ \ \ \ \\
 (e_{p,i}, e_{p^2,i}, e_{p^3,i}, \ldots ) - (\zeta_{1,2,i},e_{p^2,i},\zeta_{3,2,i},\ldots) = (0, \ \ \ \ \ \ 0 \ \ \ \ \ \ \ , e_{p^3,i}-\zeta_{3,2,i}, \ldots) \\
  \ \ \ \ \ \ \ \  \ \ \ \ \ \ \ \ \ \ \ \  \ \ \ \ \ \ \ \ \ \ \ \ \  \ \ \ \ \ \ \ \ \ \ \ \ \vdots   \\
\end{array}
\]
showing the approximation of $({\bf e}_{p^r,i})_r$  to $e_{i}$.
 \end{proof}

\subsection{\it  Approximating purely $p$--adic elements }

The approximation of $e_{i}$ indicates the method for approximating purely $p$--adic elements in $\varprojlim_{r} {\Bbb C}{C_r^n}$, in  $\varprojlim_{r} {\Bbb C}{\mathcal F} _{p^r,n}$ and in ${\rm Y}_{p^{\infty},n}(u)$. Indeed, we give first the following definition, imitation Definition~\ref{truncexp}:

\begin{defn}\label{truncexpabstr} \rm
Let $\underleftarrow{y}=(y_1, y_2, \ldots)$ an element in $\varprojlim_{r} {\Bbb C}{\mathcal F} _{p^r,n}$ resp. in ${\rm Y}_{p^{\infty},n}(u)$, in reduced form according to Definition~\ref{preduced}. For any indices $r,s$ with $r\geq s$ we define the element $y_{r,s}$ in ${\Bbb C}{\mathcal F} _{p^r,n}$  resp. ${\rm Y}_{p^r,n}(u)$ as the formal expression of $y_s$, but with the generators $t_{s,i}$ replaced by the generators $t_{r,i}$, for all $i$.
 The element $y_{r,s}$ is called  {\it the $s$--truncation of $y_r$}.

Further, we define the element $y_{s,r}$ in ${\Bbb C}{\mathcal F} _{p^s,n}$  resp. ${\rm Y}_{p^n}(u)$ as the formal expression of $y_r$, but with the generators $t_{r,i}$ replaced by the generators $t_{s,i}$, for all $i$.
The element $y_{s,r}$ is called  {\it the $r$--expansion of $y_s$}.
\end{defn}

Note that, either way, $y_{r,r} = y_r$. We define now the element ${\bf y}_r$ by $r$--truncating each term in $\underleftarrow{y}$ after the $r$th position and by $r$--expanding each term before the $r$th position. That is,
\[
{\bf y}_r = (y_{1,r}, y_{2,r}, y_r, y_{r+1,r}, y_{r+2,r}, \ldots)
\]

\begin{thm}\label{pprox}
For any index $r$ the element ${\bf y}_r$ is a constant sequence in $\varprojlim_{r} {\Bbb C}{\mathcal F}_n = {\Bbb C}{\mathcal F}_n$ resp. in ${\rm Y}_{n}(u)$. Moreover we have the approximation:
\begin{equation}\label{eiappx}
\underleftarrow{y} = \lim_{r} {\bf y}_r
\end{equation}
\end{thm}


\begin{proof}
Note that $\varphi ^r_s(y_r) = y_s$. Also, by construction, $\varphi ^{r}_s(y_{r,s}) = y_s$, that is, the elements $y_{r,s}$ and $y_s$ are coherent. In fact, $y_{r,s}$ `wraps' only once on $y_s$ via the map $\varphi ^r_s$. Moreover, since $\varphi ^r_s(y_r) = y_s$, it is also true by construction that $\varphi ^r_s(y_r) = y_{s,r}$. In fact, $y_r$ `wraps' only once on $y_{s,r}$ via the map $\varphi ^r_s$.
 Completely analogous facts are valid for the map $\phi ^r_s$ in place of $\varphi ^r_s$. Hence, the entries of ${\bf y}_r$ are coherent.

 The element  ${\bf y}_r$ is by construction a constant sequence. To see now that these constant sequences approximate $\underleftarrow{y}$ we subtract them successively from $\underleftarrow{y}$ and we confirm that the zero--sequence is gradually forming.
 \end{proof}

\section{A topological Markov trace}

In \cite{ju} the first author constructed linear Markov traces on the Yokonuma--Hecke algebras. The aim of this section is to extend these traces to a $p$--adic Markov trace on the algebra ${\rm Y}_{p^{\infty},n}(u)$.

\subsection{\it A  Markov trace on the Yokonuma--Hecke algebra}

The natural inclusions ${\mathcal F}_{d,n} \subset {\mathcal F}_{d,n+1}$ of the modular framed braid groups induce the algebra inclusions ${\Bbb C}{\mathcal F}_{d,n} \subset {\Bbb C}{\mathcal F}_{d,n+1}$ (setting ${\Bbb C}{\mathcal F}_{d,0}:={\Bbb C}$), which induce  the tower of algebras:
\begin{equation}\label{tower}
{\rm Y}_{d,0}(u)  \, \subset {\rm Y}_{d,1}(u) \subset {\rm Y}_{d,2}(u) \subset \ldots
\end{equation}
(setting ${\rm Y}_{d,0}(u):={\Bbb C}$). Thus, given $d$, we have the inductive system $\left( {\rm Y}_{d,n}(u) \right)_{n\in{\Bbb N}}$. Let $ {\rm Y}_{d,\infty}(u)$ be the corresponding inductive limit. Then we have the following.

\begin{thm}[cf. Theorem 12 in \cite{ju}]\label{trace}
Let $d$ a positive integer. For indeterminates $z$, $x_1$, $\ldots, x_{d-1}$ there exists a unique linear Markov trace ${\rm tr}_d = ({\rm tr}_{d, n})_{n\in {\Bbb N}}$
$$
 {\rm tr}_d:  {\rm Y}_{d,\infty}(u) \longrightarrow   {\Bbb C}[z, x_1, \ldots, x_{d-1}]
$$
 defined inductively on $n$ by the following rules:
$$
\begin{array}{rcll}
{\rm tr}_{d, n}(ab) & = & {\rm tr}_{d, n}(ba)  \qquad &  \\
{\rm tr}_{d, n}(1) & = & 1 & \\
{\rm tr}_{d,n+1}(ag_n) & = & z\, {\rm tr}_{d, n}(a) \qquad & (\text{Markov property} )\\
{\rm tr}_{d, n+1}(at_{n+1}^m) & = & x_m {\rm tr}_{d,n}(a)\qquad  & (  m = 1, \ldots , d-1)
\end{array}
$$
where $a,b \in {\rm Y}_{d,n}(u)$.
\end{thm}
Lifting to framed braids, in the second rule is meant the trace of the identity braid on $n$  strands with all framings zero. The third rule is the so--called {\it Markov property} of the trace. See Figure~6 for topological interpretations of the last two rules.

\smallbreak
\begin{figure}[H]
 \begin{picture}(330,70)

\put(-23,36){${\rm tr}_d$}
\qbezier(0,0)(-9,37)(0,74)  
\put(13,36){$a$}
\qbezier(8,0)(8,12)(8,24)
\qbezier(8,54)(8,62)(8,70)
\qbezier(22,54)(22,62)(22,70)
\qbezier(50,24)(52,47)(50,70)
\qbezier(30,0)(29,2)(30,4)
\qbezier(50,0)(51,2)(50,4)
\qbezier(58,0)(67,37)(58,74)  
\put(68,36){$=$}
\put(78,36){$z\,{\rm tr}_d$}
\qbezier(107,0)(96,37)(107,74)  
\qbezier(0,24)(0,39)(0,54)
\qbezier(30,24)(30,39)(30,54)

\qbezier(0,54)(15,54)(30,54)
\qbezier(0,24)(15,24)(30,24)
\qbezier(30,4)(30,9)(40,14)
\qbezier(40,14)(50,19)(50,24)
\qbezier(50,4)(50,8)(45,11)
\qbezier(35,17)(30,20)(30,24)
\qbezier(110,24)(110,39)(110,54)
\qbezier(140,24)(140,39)(140,54)

\qbezier(110,54)(125,54)(140,54)
\qbezier(110,24)(125,24)(140,24)
\qbezier(118,0)(118,12)(118,24)
\qbezier(118,54)(118,62)(118,70)
\qbezier(132,54)(132,62)(132,70)
\qbezier(132,0)(132,12)(132,24)
\put(122,36){$a$}
\qbezier(146,0)(156,37)(146,74) 
\put(160,36){$,$}
\put(175,36){${\rm tr}_d$}
\qbezier(197,0)(188,37)(197,74)  
\qbezier(202,24)(202,39)(202,54)
\qbezier(232,24)(231,39)(232,54)

\qbezier(202,54)(217,54)(232,54)
\qbezier(202,24)(217,24)(232,24)
\qbezier(210,0)(210,12)(210,24)
\qbezier(224,54)(224,62)(224,70)
\qbezier(210,54)(210,62)(210,70)
\qbezier(224,0)(224,12)(224,24)

\qbezier(240,0)(240,34)(240,68)
\put(238, 70){\tiny{ $m$}}
\put(215,36){$a$}
\qbezier(253,0)(261,37)(253,74) 
\put(261,36){$=$}
\put(272,36){$x_m {\rm tr}_d$}
\qbezier(306,0)(298,37)(306,74)  
\qbezier(310,24)(310,39)(310,54)
\qbezier(340,24)(340,39)(340,54)

\qbezier(310,54)(325,54)(340,54)
\qbezier(310,24)(325,24)(340,24)
\qbezier(318,0)(318,12)(318,24)
\qbezier(318,54)(318,62)(318,70)
\qbezier(332,54)(332,62)(332,70)
\qbezier(332,0)(332,12)(332,24)

\put(321,36){$a$}
\qbezier(345,0)(353,37)(345,74) 

\end{picture}
\caption{Topological interpretations of the trace rules}
\label{fig6}
\end{figure}

The key in the construction of ${\rm tr}_d$ is that ${\rm Y}_{d,n+1}(u)$ has a `nice' inductive linear basis. Indeed, every element of ${\rm Y}_{d,n+1}(u)$ is a unique linear combination of words, each of one of the following types:
\begin{equation}\label{basis}
 w_ng_ng_{n-1} \ldots g_i t_i^k \text{ \ \ \ \ \ or \ \ \ \ \ } w_nt_{n+1}^k, \quad  k\in {\Bbb Z}/d{\Bbb Z}
\end{equation}
where $w_n \in {\rm Y}_{d,n}(u)$.  Thus, the above
words furnish an inductive basis for ${\rm Y}_{d,n+1}(u)$, and each one involves  $g_n$ or a power of $t_{n+1}$ at most once. Cf.~\cite{ju} for details.

\begin{rem}\label{ocneanu}\rm
In the case $d=1$, when the algebra ${\rm Y}_{1,n}(u)$ coincides with the Iwahori--Hecke algebra ${\rm H}_n(u)$, the trace ${\rm tr}_1$ coincides with the Ocneanu trace (cf., for example, \cite{jo}).
\end{rem}

We shall use the notation ${X}_d$ for the set of indeterminates $\{x_1, x_2, \ldots ,x_{d-1}\}$ of the trace  ${\rm tr}_d$:
\begin{equation}\label{xd}
{X}_d := \{x_1, x_2, \ldots ,x_{d-1}\}.
\end{equation}

\subsection{\it The $p$--adic Markov trace}

Let $R$ denote the polynomial ring ${\Bbb C}[z]$ and let $r$ be a positive integer.
We denote $R\left[\mathfrak{X}_r\right]$ the polynomial ring on  the
indeterminates  of the set $\mathfrak{X}_r$,
$$
\mathfrak{X}_r: = \{x_a\,;\, a\in {\Bbb Z}/p^r{\Bbb Z}\}
$$

For all positive integers $r$, $s$ such that $r\geq s$ we have the  ring homomorphism
\begin{equation}\label{delta}
\delta^{r}_{s} :R\left[\mathfrak{X}_{r}\right] \longrightarrow
R\left[\mathfrak{X}_s\right]
\end{equation}
which is defined via the mapping:  $x_a\mapsto x_b$, where $b= \vartheta_s^r(a)$ (recall $\vartheta_s^r$ from Eq.~\ref{varthetars}).
It is then a  routine to  prove the following lemma.

\begin{lem}\label{invrings}
The family $\left( R[\mathfrak{X}_r], \delta_s^r \right)$ is an
inverse system of polynomial rings indexed by ${\Bbb N}$.
\end{lem}

\begin{notations}\rm
We shall use the notations $ \tau_r$ in  place of \, ${\rm tr}_{p^r}$ and $ \tau_{r, n}$ in  place of ${\rm tr}_{p^r, n}$. With these notations we have: $\tau_r = (\tau_{r, n})_{n\in {\Bbb N}}$.
\end{notations}

\begin{lem}\label{commute}
 The diagram below is commutative.
$$
\begin{diagram}
 \node{{\rm Y}_{p^r,n}(u)}\arrow{e,t}{\phi_s^{r}} \arrow{s,l}{\tau_{r,n}} \node{{\rm Y}_{p^s,n} (u)}\arrow{s,r}{\tau_{s, n }}\\
\node{R[\mathfrak{X}_{r}]} \arrow{e,b}{\delta_s^{r}}\node{R[\mathfrak{X}_s]}
\end{diagram}
$$
\end{lem}

\begin{proof}
The proof is by induction on $n$. The lemma is immediate for $n=1$. Assume the lemma is true for some $n$. In order to prove it for $n+1$ we must check that
$(\tau_{s,n+1} \circ \phi_s^{r})(x) = (\delta_s^{r}\circ\tau_{r, n+1})(x)$ for all $x\in {\rm Y}_{p^r,n+1}(u)$. Since, by definition, the maps $\phi_s^{r}, \tau_{s,n+1}$ and $\tau_{r,n+1}$ are linear, it suffices to prove that $(\tau_{s,n+1} \circ \phi_s^{r})(\alpha) = (\delta_s^{r}\circ\tau_{r,n+1})(\alpha)$, for any $\alpha$ in the inductive basis of ${\rm Y}_{p^r,n+1}(u)$.

Assume first that $\alpha= w_ng_ng_{n-1} \ldots g_i t_{r,i}^k $, where $w_n\in {\rm Y}_{p^r,n}(u)$. Then:
\begin{eqnarray*}
 \tau_{s,n+1}  (\phi_s^{r}(\alpha))
& = &
 \tau_{s,n+1}\left(\phi_s^{r}(w_n) g_n g_{n-1} \ldots g_i t_{r,i}^k\right)
 \quad (k \text{ regarded modulo }  p^s)\\
& = &
z \, \tau_{s,n}\left(\phi_s^{r}(w_n) g_{n-1} \ldots g_i t_{r,i}^k\right)\\
& = &
z \, \tau_{s,n}\left(\phi_s^{r}(w_n g_{n-1} \ldots g_i t_{r,i}^k)\right)\\
& = &
z\, \delta_s^{r}\left(  \tau_{r,n} (w_ng_{n-1} \ldots g_i t_{r,i}^k)\right)\quad (\text{induction hypothesis} )\\
& = &
\delta_s^{r}\left( z \, \tau_{r,n} (w_ng_{n-1} \ldots g_i t_{r,i}^k)\right)\\
& = &
\delta_s^{r}\left(  \tau_{r,n+1} (w_ng_ng_{n-1} \ldots g_i t_{r,i}^k)\right)\quad (\text{trace rule} )\\
& = &
\delta_s^{r}\left(\tau_{r,n+1}(\alpha)\right).
\end{eqnarray*}

Assume now that $\alpha= w_n t_{r,n+1}^k $. Then:
\begin{eqnarray*}
 \tau_{s,n+1}  \left(\phi_s^{r}(\alpha)\right)
& = &
  \tau_{s,n+1}  \left(\phi_s^{r}(w_n ) t_{r, n+1}^k\right)  \quad (k \text{ regarded modulo }  p^s) \\
& = &
x_k \, \tau_{s,n}  \left(\phi_s^{r}(w_n )\right) \\
& = &
x_k \, \delta_s^{r}\left(  \tau_{r,n} (w_n)\right) \quad (\text{induction hypothesis} )\\
& = &
\delta_s^{r}\left(  x_k\, \tau_{r,n} (w_n)\right) \\
& = &
\delta_s^{r}\left( \tau_{r,n+1}(w_n t_{r, n+1}^k)\right) \quad (\text{trace rule} )\\
& = &
\delta_s^{r}\left( \tau_{r,n+1}(\alpha)\right).
\end{eqnarray*}

Hence the proof is concluded.
\end{proof}

\begin{defn}\label{pvars}\rm
For  a $p$--adic integer $ \underleftarrow{a} = (a_1, a_2, \ldots) \not = 0$ we shall denote
$$
x_{\underleftarrow{a}}:= (x_{a_1}, x_{a_2}, \ldots ) \in  \varprojlim_{r\in{\Bbb N}} R[\mathfrak{X}_r]
$$
and we shall call $x_{\underleftarrow{a}}$ a $p$--{\it adic indeterminate}. Further, for an almost constant sequence
$a_i = (a_1,\ldots,a_{i-1}, a_i, a_i, \ldots)$ in ${\Bbb Z}$ we shall denote
$$
x_{a_i} := (x_{a_1},\ldots, x_{a_{i-1}},x_{a_i},x_{a_i}, \ldots ) \in  \varprojlim_{r\in{\Bbb N}} R[{\mathfrak X}_r]
$$
 and we shall say that $x_{a_i}$ is a {\it constant indeterminate}. Finally, we make  the convention $x_0:=1$.
\end{defn}

Let $r$,  $s$ and $v$ be positive integers such that $r\geq s\geq v$. By Lemmas~\ref{invrings} and~\ref{commute},
we have the following commutative diagram:
$$
\begin{diagram}
\node{\cdots}
\node{{\rm Y}_{p^{v},n}(u)}
\arrow{w,t}{} \arrow{s,l}{\tau_{v, n}} \node{{\rm Y}_{p^{s},n}(u)}\arrow{w,t}{\rho_v^{s}}
\arrow{s,l}{\tau_{s, n}}
\node{{\rm Y}_{p^r,n}(u)}
\arrow{w,t}{\rho_{s}^{r}} \arrow{s,l}{\tau_{r,n}} \node{\cdots}\arrow{w,t}{}
\\
\node{\cdots}
\node{R\left[\mathfrak{X}_v\right] }
\arrow{w,t}{}
\node{R\left[\mathfrak{X}_{s}\right] }\arrow{w,t}{\delta_v^s}
\node{R\left[\mathfrak{X}_{r}\right] }
\arrow{w,t}{\delta_s^r}
\node{\cdots}\arrow{w,t}{}
\end{diagram}
$$

The diagram above induces  a unique ring homomorphism  $\tau_{p^{\infty}, n} := \varprojlim_r \tau_{r, n}$,
$$
\tau_{p^{\infty}, n} :  {\rm Y}_{p^{\infty} ,n}(u)\longrightarrow\varprojlim_{r\in {\Bbb N}} R[\mathfrak{X}_r]
$$
Note that, by definition, $\tau_{p^{\infty}, n} = \left(\tau_{1, n}, \tau_{2, n}\tau_{3, n}, \ldots \,\right)$.

The natural inclusions (\ref{tower})  of Yokonuma--Hecke algebras induce, by construction, an inductive system  of $p$--adic Yokonuma--Hecke algebras:
$$
{\rm Y}_{p^{\infty},0}(u)  \, \subset {\rm Y}_{p^{\infty},1}(u) \subset {\rm Y}_{p^{\infty},2}(u) \subset \ldots
$$
(setting ${\rm Y}_{p^{\infty},0}(u):={\Bbb C}$). Let ${\rm Y}_{p^{\infty},\infty}(u)$ be  the associated inductive limit:
$$
{\rm Y}_{p^{\infty},\infty}(u) := \varinjlim_{n\in {\Bbb N}}{\rm Y}_{p^{\infty},n}(u)
$$

\begin{thm}\label{ptrace}
There exists a unique $p$--adic linear Markov trace $\tau_{p^{\infty}} = \left( \tau_{p^ {\infty}, n}\right)_{n\in {\Bbb N}}$\,,
$$
\tau_{p^{\infty}} :  {\rm Y}_{p^{\infty},\infty}(u) \longrightarrow \varprojlim_{r\in{\Bbb N}} R[{\mathfrak X}_r]
$$
such that:
$$
\begin{array}{rcll}
\tau_{p^ {\infty}, n}(\underleftarrow{a}\underleftarrow{b}) & = & \tau_{p^ {\infty}, n}(\underleftarrow{b}\underleftarrow{a}) \qquad &    \\
\tau_{p^ {\infty}, n} (1) & = & 1 & \\
\tau_{p^ {\infty}, n+1}({\underleftarrow{y}}g_n) & = & z\, \tau_{p^ {\infty}, n}({\underleftarrow{y}})  \qquad &  (\text{Markov  property} )  \\
\tau _{p^ {\infty}, n+1}({\underleftarrow{y}} {\bf t}_{n+1}^m) & = & x_m\, \tau_{p^ {\infty}, n}({\underleftarrow{y}}) \qquad  & (m \in {\Bbb Z}) \\
\tau_{p^ {\infty}, n+1} ({\underleftarrow{y}} {\bf t}_{n+1}^{\underleftarrow{m}}) & = & x_{\underleftarrow{m}} \, \tau_{p^ {\infty}, n}({\underleftarrow{y}})  \qquad  & ({\underleftarrow{m}} \in {\Bbb Z}_p)
\end{array}
$$
where $\underleftarrow{a},\underleftarrow{b},{\underleftarrow{y}} \in {\rm Y}_{p^{\infty},n}$.
\end{thm}

\begin{proof}
The trace $\tau_{p^ {\infty}}$ is unique by Theorem~\ref{trace} and, by  construction, it satisfies all  properties in the statement. Indeed, let us check the third and the fifth property. For ${\underleftarrow{y}}=(y_r)_r$, with $y_r \in {\rm Y}_{p^r,n}(u)$, and ${\underleftarrow{m}}=(m_r)_r \in {\Bbb Z}_p$ we have
$\tau_{p^{\infty}, n+1} ({\underleftarrow{y}} g_n) = \left(\tau_{r,n+1}(y_r g_n)\right)_r$. Hence, using Theorem \ref{trace}:
$$
\tau_{p^{\infty}, n+1} ({\underleftarrow{y}} g_n)= \left(z\, \tau_{r,n}(y_r)\right)_r
= z\, \left(\tau_{r,n}(y_r)\right)_r = z\, \tau_{p^{\infty}, n} ({\underleftarrow{y}}).
 $$
 Analogously, we have: $\tau_{p^{\infty}, n+1} ({\underleftarrow{y}} {\bf t}_{n+1}^{\underleftarrow{m}}) =
\left(\tau_{r,n+1}(y_r t_{r,n+1}^{m_r})\right)_r$. Then:
$$
\tau_{p^{\infty}, n+1} ({\underleftarrow{y}} {\bf t}_{n+1}^{\underleftarrow{m}}) =  \left(x_{m_r}\, \tau_{r,n}(y_r)\right)_r = (x_{m_r})_r \left(\tau_{r,n}(y_r)\right)_r
 =  x_{\underleftarrow{m}} \, \tau_{p^{\infty}, n}({\underleftarrow{y}}).
$$
\end{proof}
\begin{rem} \rm
We have the approximation:
$$
\tau_{p^ {\infty}, n}({\underleftarrow{y}}) = \lim_{r}\tau_{p^{\infty}, n}({\bf y}_r)
$$
for  ${\underleftarrow{y}}  = (y_r)_r = \lim_{r} {\bf y}_r \in{\rm Y}_{p^{\infty},n}$.
In this sense  $\tau_{p^{\infty}, n}({\underleftarrow{y}})$ can be viewed as an infinite series in $\varprojlim R[{\mathfrak X}_r]$:
$$
\tau_{p^{\infty}, n}({\underleftarrow{y}})  = \sum_{r=0}^{\infty}\tau_{p^{\infty}, n} ({\bf y}_r).
$$
Indeed,  the approximating sequence to $\tau_{p^{\infty}, n}({\underleftarrow{y}})$ by the values of $\tau_{p^{\infty}, n}$ on the elements ${\bf y}_r \in {\rm Y}_{n}(u)$ follows easily by Theorem~\ref{pprox} and our usual approximation arguments.
\end{rem}

\subsection{\it Computations}

We shall now give some computations of the traces ${\rm tr}_d$ and $\tau_{p^{\infty}}$.

\smallbreak

\noindent $\bullet$ For the element ${\bf t}^k\in {\rm Y}_{p^{\infty},1}(u)$ ($1$--strand braid with framing $k$) we have:
$$
\tau_{p^{\infty}} ({\bf t}^k) = (\tau_1(t_1^{k(\text{mod}\,\, p)}), \tau_2(t_2^{k(\text{mod}\,\, p^2)}), \ldots )=  (x_k, x_k,\ldots ) = x_k \in  R[{\mathfrak X}_r],
$$
by Definition~\ref{pvars}.
\bigbreak

\noindent $\bullet$ For the element ${\bf t}^{\underleftarrow{a}} = (t_1^{a_1}, t_2^{a_2}, \ldots) \in {\rm Y}_{p^{\infty} ,1}(u)$, where $\underleftarrow{a}=(a_r)_r$, we have:
$$
\tau_{p^{\infty}} ({\bf t}^{\underleftarrow{a}}) = (\tau_1(t_1^{a_1}), \tau_2(t_2^{a_2}), \ldots )
= (x_{a_1}, x_{a_2}, \ldots ) = x_{\underleftarrow{a}} \in \varprojlim R[{\mathfrak X}_r].
$$
Further, since $ {\bf t}^{\underleftarrow{a}} =  \lim_{r} {\bf t}^{a_r}$ we have the following trace approximation:
$$
\tau_{p^{\infty}}({\bf t}^{\underleftarrow{a}}) = \lim_{r} \tau_{p^{\infty}} ({\bf t}^{a_r})  = \lim_{r} x_{a_r} \in \varprojlim R[{\mathfrak X}_r]
$$

\noindent $\bullet$
For the $n$--strand identity braid with framings $k_1,\ldots,k_n\in {\Bbb Z}$ we have:
$$
{\rm tr}_d (t_1^{k_1}\ldots t_n^{k_n} )  =  x_{k_1}\ldots x_{k_n}
$$
$$
\tau_{p^{\infty}} ({\bf t}_1^{k_1}\ldots {\bf t}_n^{k_n}) =
\left(\tau_{r}(t_{r,1}^{k_1} \, \ldots \, t_{r,n}^{k_n})\right)_r =
x_{k_1}\ldots x_{k_n} =
\tau_{p^{\infty}} ({\bf t}_1^{k_1})\ldots \tau_{p^{\infty}}({\bf t}_n^{k_n} ) \in R[{\mathfrak X}_r]
$$

\bigbreak
\noindent$\bullet$
For the $n$--strand identity braid with framings $\underleftarrow{a_i} = (a_{ri})_r \in {\Bbb Z}_p$ we have:
\begin{eqnarray*}
\tau_{p^{\infty}} ({\bf t}_1^{\underleftarrow{a_1}}\ldots {\bf t}_n^{\underleftarrow{a_n}})
&  = &
 \left(\tau_{r}(t_{r,1}^{a_{r1}} \, \ldots \, t_{r,n}^{a_{rn}})\right)_r \\
 &  = &
(x_{a_{11}}\ldots x_{a_{1n}},x_{a_{21}}\ldots x_{a_{2n}}, \ldots ) \\
 & = &
(x_{a_{11}},x_{a_{21}},\ldots) \ldots (x_{a_{1n}},x_{a_{2n}},\ldots)  \\
 & = &
 x_{\underleftarrow{a_1}}\ldots x_{\underleftarrow{a_n} } \\
 & = &
 \tau_{p^{\infty}} ({\bf t}_1^{\underleftarrow{a_1}})\ldots \tau_{p^{\infty}}({\bf
t}_n^{\underleftarrow{a_n}}) \in \varprojlim R[{\mathfrak X}_r].
\end{eqnarray*}
Further, we have the approximation:
$$
\tau_{p^{\infty}}({\bf t}_1^{\underleftarrow{a_1}} \ldots {\bf t}_n^{\underleftarrow{a_n}}) =
\lim_{r} \tau_{p^{\infty}} ({\bf t}_1^{a_{r1}} \ldots  {\bf t}_n^{a_{rn}} ) =
\lim_{r} \left(x_{a_{r1}} \ldots  x_{a_{rn}} \right) \in \varprojlim R[{\mathfrak X}_r].
$$

\bigbreak
\noindent$\bullet$
For the elements $e_{d,i}\in {\rm Y}_{d,n}(u)$ and $e_{i} \in {\rm Y}_{p^{\infty},n}(u) \quad(i=1,\ldots , n-1)$ :
\begin{equation}\label{Ed}
E_d :={
\rm tr}_d (e_{d,i})  =
{\rm tr}_d \left(\frac{1}{d} \sum_{m=0}^{d-1}t_i^m t_{i+1}^{-m}\right) =
\frac{1}{d} \sum_{m=0}^{d-1}x_m x_{d-m}
\end{equation}

$$
\tau_{p^{\infty}}(e_{i}) =
\tau_{p^{\infty}} \left((e_{p,i}, e_{p^2,i}, \ldots )\right) = \left(\tau_{r}\left(e_{p^r,i}\right)\right)_r =
(E_p, E_{p^2},\ldots )=
(E_{p^r})_r \in \varprojlim R[{\mathfrak X}_r]
$$
Further, from (\ref{eiappx}) we have the approximation:
$$
\tau_{p^{\infty}}(e_{i}) = \lim_{r}\, \tau_{p^{\infty}}\left({\bf e}_{p^r,i}\right) = \lim_{r}\left(\frac{1}{p^r}\sum_{m=0}^{p^r-1}x_m x_{-m}\right)
=  \lim_{r} E_{p^r} \in \varprojlim R[{\mathfrak X}_r]
$$
\bigbreak
\noindent$\bullet$
For the elements $e_{d,i}g_i\in {\rm Y}_{d,n}(u)$ and $e_{i}g_i \in {\rm Y}_{p^{\infty},n}(u)\quad(i=1,\ldots , n-1)$ we have the following more general lemma.

\begin{lem}\label{engn}
Let $y\in {\rm Y}_{d,n}(u)$ and $\underleftarrow{y} = (y_r)_r =\lim_{r} {\bf y}_r \in {\rm Y}_{p^ {\infty},n}(u)$. Then:
\begin{enumerate}
\item[(i)] ${\rm tr}_d (y \, e_{d,n}g_n)  = {\rm tr}_d (y g_n) = z\, {\rm tr}_d (y)$
\item[(ii)]  $\tau_{p^{\infty}} (\underleftarrow{y} \, e_ng_n)  = \tau_{p^{\infty}} (\underleftarrow{y} g_n) = z\, \tau_{p^{\infty}} (\underleftarrow{y})$
\item[(iii)]
$\tau_{p^{\infty}}(\underleftarrow{y} \, e_ng_n) = z\, \lim_{r}\tau_{p^{\infty}} ({\bf y}_r).$
\end{enumerate}
\end{lem}
\begin{proof} (i) We have $y \, e_{d,n}g_n = \frac{1}{d} \sum_{m=0}^{d-1}y \,t_n^m t_{n+1}^{-m} g_n = \frac{1}{d} \sum_{m=0}^{d-1}y \,t_n^m g_n t_n^{-m}$
\vspace{.05in}
\noindent so, applying the trace yields the statement.

\vspace{.05in}
\noindent (ii) $\tau_{p^{\infty}}(\underleftarrow{y} \, e_ng_n) =
\left(\tau_{r}(y_r \, e_{p^r,n}g_n)\right)_r \stackrel{(i)}{=}
(z\,\tau_{r}(y_r))_r =
z\, \tau_{p^{\infty}}(\underleftarrow{y})=
\tau_{p^{\infty}}(\underleftarrow{y}g_n)$.

\vspace{.05in}
\noindent Finally, (iii) follows immediately from (ii) and Theorem~\ref{pprox}.
\end{proof}

\bigbreak
\noindent $\bullet$ For $g_i^2 \in {\rm Y}_{d,n}(u)$ and for $ g_i^2\in {\rm Y}_{\infty ,n}$, where $g_i^2= (g_i^2, g_i^2,\ldots)$ we have:

\vspace{.07in}
\noindent ${\rm tr}_d (g_i^2) = 1  + (u-1)z + (u-1)E_d$

\vspace{.07in}
\noindent $\tau_{p^{\infty}} (g_i^2 )= 1+(u-1)z + (u-1) (E_{p^r})_r  = 1+ (u-1)z + (u-1) \tau_{p^{\infty}} (e_{i}) \in \varprojlim R[{\mathfrak X}_r]$




\bigbreak
\noindent $\bullet$ For $g_i^3 \in {\rm Y}_{d,n}(u)$ and for $ g_i^3\in {\rm Y}_{\infty ,n}$ we have:

\vspace{.07in}
\noindent ${\rm tr}_d (g_i^3) = (u^2-u+1)z + (u^2-u)E_d$.

\vspace{.07in}
\noindent ${\rm tr}_d (g_i^{-3}) = (u^{-3}-u^{-2}+u^{-1})z + (u^{-3}-u^{-2}+u^{-1}-1)E_d$.

\vspace{.07in}
\noindent $\tau_{p^{\infty}} (g_i^3 )= (u^2-u+1)z + (u^2-u)(E_{p^r})_r \in \varprojlim R[{\mathfrak X}_r]$.

\vspace{.07in}
\noindent $\tau_{p^{\infty}} (g_i^{-3}) = (u^{-3}-u^{-2}+u^{-1})z + (u^{-3}-u^{-2}+u^{-1}-1)(E_{p^r})_r \in \varprojlim R[{\mathfrak X}_r]$.

\subsection{\it The Markov braid equivalence}

From the topological point of view, closing a framed braid gives rise to an oriented {\it framed link} and closing a $p$--adic framed braid gives rise to an oriented {\it $p$--adic framed link}, see Figure~\ref{plink}. By `closing' a braid $\beta$ we mean the standard closure, denoted $\widehat{\beta}$, where we join with simple disjoint arcs the corresponding top and bottom endpoints of the braid. Conversely, by the classical Alexander theorem (adapted to the various framed braid settings), an oriented framed link can be isotoped to the closure of a framed braid.

\smallbreak
Further, by the classical Markov theorem (also adapted to the various framed braid settings), isotopy classes of oriented framed links are in one--to--one correspondence with equivalence classes of framed braids. More precisely, we have the following  result, which is well--known for the case of classical framed links (see for example \cite{KS}), and which we also adapt here for the cases of modular framed links and $p$--adic framed links.

\begin{thm}[Markov equivalence for framed braids and $p$--adic framed braids]\label{markov}
 Isotopy classes of oriented framed links (resp. modular framed links) are in bijection with equivalence classes of framed braids in $\cup_{n\in{\Bbb N}}{\mathcal F}_{n}$ (resp. $\cup_{n\in{\Bbb N}}{\mathcal F}_{d,n}$). The equivalence relation is generated by the two moves:
\vspace{.07in}

\noindent (i) Conjugation: \ $\alpha \beta \sim \beta \alpha, \, \quad \alpha,\beta \in {\mathcal  F}_n$ (resp. ${\mathcal  F}_{d,n}$)
\vspace{.07in}

\noindent (ii) Markov move: $\alpha \sim \alpha{\sigma_n}^{\pm 1}, \quad \alpha \in {\mathcal  F}_n$ (resp. ${\mathcal  F}_{d,n}$)

\smallbreak

Further, isotopy classes of $p$--adic framed links  are in bijection with equivalence classes of $p$--adic framed braids in $\cup_{n\in{\Bbb N}}{\mathcal F}_{p^{\infty},n}$ under the following equivalence relation:

Two $p$--adic framed braids ${\underleftarrow{\alpha}} = (\alpha_r)_r, \ {\underleftarrow{\beta}} = (\beta_r)_r \in \cup_{n\in{\Bbb N}}{\mathcal F}_{p^{\infty},n}$ are equivalent if and only if for every $r$ the modular framed braids $\alpha_r$ and $\beta_r$ are Markov equivalent in $\cup_{n\in{\Bbb N}}{\mathcal F}_{d,n}$.

In view of the isomophisms (\ref{isoms}) the Markov equivalence of $p$--adic framed braids is generated by the moves:

\vspace{.07in}

\noindent (i) Conjugation: \ ${\underleftarrow{\alpha}} {\underleftarrow{\beta}} \sim
{\underleftarrow{\beta}} {\underleftarrow{\alpha}}, \, \quad {\underleftarrow{\alpha}}, {\underleftarrow{\beta}} \in {\mathcal F}_{p^{\infty},n}$
\vspace{.07in}

\noindent (ii) Markov move: ${\underleftarrow{\alpha}} \sim {\underleftarrow{\alpha}}{\sigma_n}^{\pm 1}, \quad {\underleftarrow{\alpha}} \in {\mathcal F}_{p^{\infty},n}$.
\end{thm}

According to Theorem~\ref{markov}, any invariant of oriented framed links has to agree on the closures of the braids $\alpha$, $\alpha{\sigma_n}$ and $\alpha{\sigma_n}^{- 1}$. Note the resemblance of the conjugation rule and the Markov property in Theorems~\ref{trace} and \ref{ptrace} with moves (i) and (ii) of Theorem~\ref{markov}.
 Having, now, present the recipe of Jones\cite{jo} we will try to define an invariant by re--scaling and normalization of the trace ${\rm tr}_d$ and the $p$--adic trace $\tau_{p^{\infty}}$.  In order to do that we need that the expression ${\rm tr}_d(\alpha g_n^{-1})$, for $\alpha \in {\rm Y}_{d,n}(u)$, factors through ${\rm tr}_d(\alpha)$, just like ${\rm tr}_d(\alpha g_n)$ does from the Markov property of the trace. Yet, we have:
\begin{equation}\label{traginverse}
{\rm tr}_d(\alpha g_n^{-1}) = {\rm tr}_d(\alpha g_n) + (u^{-1}-1){\rm tr}_d(\alpha e_{d,n})
 + (u^{-1}-1){\rm tr}_d(\alpha e_{d,n}g_n)
\end{equation}
Analogous requirements apply to $\tau_{p^{\infty}}(\alpha g_n^{-1})$, for $\alpha \in {\rm Y}_{p^{\infty},n}(u)$. Here we
have:
$$
\tau_{p^{\infty}} (\alpha g_n^{-1}) = \tau_{p^{\infty}}(\alpha g_n) + (u^{-1}-1) \tau_{p^{\infty}} (\alpha e_n)
 + (u^{-1}-1) \tau_{p^{\infty}}(\alpha e_ng_n)
$$
By Lemma~\ref{engn} (i) and (ii), we only need further that the traces ${\rm tr}_d$ and $\tau_{p^{\infty}}$
 satisfy the multiplicative properties:
\begin{equation}\label{mul}
{\rm tr}_d(\alpha e_{d,n}) = {\rm tr}_d(\alpha)\, {\rm tr}_d( e_{d,n}) \quad \alpha \in {\rm Y}_{d,n}(u)
\end{equation}
and
\begin{equation}\label{pmul}
\tau_{p^{\infty}}(\alpha e_n) = \tau_{p^{\infty}} (\alpha)\, \tau_{p^{\infty}}( e_n) \quad \alpha \in {\rm Y}_{p^{\infty},n}(u)
\end{equation}
With these properties we could then define framed link invariants using the same method as for defining the Jones polynomial\cite{jo}. Unfortunately, we do not have a  nice formula for ${\rm tr}_d(\alpha \, e_{d,n})$, and this causes a similar problem for $\tau_{p^{\infty}}(\alpha e_n)$. The reason is that the element $e_{d,n}$ involves the $n$th strand of the braid $\alpha$.

\section{The $E$--condition}

The goal of this section is to find conditions, so that equations (\ref{mul}) and (\ref{pmul}) hold. Since $d$ remains fixed throughout the section, we shall denote ${\rm tr}_d$ simply by ${\rm tr}$ and we will suppress the index $d$ from the framing generators of the algebras ${\rm Y}_{d,n}(u)$. We shall also suppress the values of the indices in the summation symbols.

\smallbreak
For $0\leq k\leq d-1$ we now define the elements:
\begin{equation}\label{editok}
e_{d,i}^{(k)}: = \frac{1}{d}\sum_{s=0}^{d-1}t_i^{k + s} t_{i+1}^{d-s}
\end{equation}
and also:
\begin{equation}\label{treditok}
E_{d}^{(k)} := {\rm tr}\left(e_{d,i}^{(k)}\right) = \frac{1}{d}\sum_{s=0}^{d-1}x_{k + s} x_{d-s}
\end{equation}
For example:
$E_3^{(2)} = \frac{1}{3}\left( 2x_2 + x_1^2\right)$. With the above notation $e_{d,i}^{(0)} = e_{d, i} $ and:
$$
E_{d}^{(0)} = {\rm tr}(e_{d,i}) := E_{d}.
$$
Note that in the definition of $E_{d}^{(k)}$ the sub--indices of the indeterminates  are regarded modulo $d$.

\begin{rem}\label{edikl}\rm
By a change of variable for $s$ it is easy to deduce the following useful formulas, stressing once more that the sub--indices of the indeterminates  are regarded modulo $d$.
$$
\frac{1}{d}\sum_{s=0}^{d-1}t_i^{k + s} t_{i+1}^{l-s} = e_{d,i}^{(k+l)}
\qquad {\rm \text and  }  \qquad
\frac{1}{d}\sum_{s=0}^{d-1}x_{k + s} x_{l-s} = E_{d}^{(k+l)} \qquad (k,l \in {\Bbb Z})
$$
\end{rem}

\subsection{\it Computing ${\rm tr}(\alpha e_{d,n})$}

By (\ref{basis})  every element  $\alpha \in {\rm Y}_{d,n}(u)$ is a unique linear combination of words in one of the following types:
$$
 w_{n-1}g_{n-1} \ldots g_i t_i^k \qquad \text{or} \qquad w_{n-1}t_n^k \qquad  \big(k\in {\Bbb Z}/d{\Bbb Z},    \ w_{n-1} \in {\rm Y}_{d,n-1}(u)\big)
$$
We shall now give some concrete computations.

\smallbreak
\noindent $\bullet$ For $n=1$ the only case is $\alpha = t_1^k$. So,
\begin{center}
${\rm tr}(\alpha) = x_k$
 \ and \ ${\rm tr}(\alpha e_{d,1}) = E_d^{(k)} = \frac{E_d^{(k)}}{x_k} {\rm tr}(\alpha)$.
\end{center}
\smallbreak
\noindent $\bullet$ For $n=2$ we have: $ {\rm (i)}\,\alpha = t_1^{l}t_2^{k}$ \ or \ ${\rm (ii)}\, \alpha = t_1^{l}g_1t_1^{k} $. Then:

\vspace{.1in}
\noindent \ (i)  \, ${\rm tr}(\alpha) = x_{k} x_l$ \quad \ \ and \quad ${\rm tr}(\alpha e_{d,2}) = x_{l} E_d^{(k)} = \frac{E_d^{(k)}}{x_k} {\rm tr}(\alpha)$.

\vspace{.1in}
\noindent (ii) \ \ ${\rm tr}(\alpha)= z\, x_{k + l}$ \quad and \quad ${\rm tr}(\alpha e_{d,2}) = z E_d^{(k + l)} = \frac{E_d^{(k + l)}}{x_{k + l}} {\rm tr}(\alpha)$.

\smallbreak
 In general we have the following results.

\begin{lem}\label{aone}
Let $\alpha = w_{n-1}t_n^{k}$ \, with \, $ w_{n-1}\in {\rm Y}_{d,n-1}(u)$. Then:
$$
{\rm tr}\left(\alpha \, e_{d,n} \right) = \frac{E_{d}^{(k)}}{x_k} \, {\rm tr}\left(\alpha\right).
$$
More generally: \qquad \qquad
${\rm tr}\left(\alpha \, e_{d,n}^{(m)} \right) = \frac{E_{d}^{(m+k)}}{x_k} \, {\rm tr}\left(\alpha\right)$.
\end{lem}

\begin{proof}
We prove the more general result. We have:
\begin{eqnarray*}
{\rm tr}\left(\alpha e_{d,n}^{(m)}\right)
 & = &
\frac{1}{d}\sum_{s}{\rm tr}\left( w_{n-1}t_n^kt_n^{m + s}t_{n+1}^{d-s}\right)\\
& = &
\frac{1}{d}\sum_{s}x_{d-s}{\rm tr}\left( w_{n-1}t_n^{m + k + s}\right)\\
& = &
\frac{1}{d}\sum_{s}x_{d-s}x_{m + k + s}{\rm tr}\left( w_{n-1}\right)\\
& = &
{\rm tr}( w_{n-1})\frac{1}{d}\sum_{s}x_{d-s}x_{m + k + s}
\, = \, {\rm tr}\left( w_{n-1}\right)E_d^{(m + k)}.
\end{eqnarray*}
On the other hand:  ${\rm tr}(\alpha) = x_{k}\, {\rm tr}\left(w_{n-1}\right)$.
\end{proof}

\begin{lem}\label{atwo}
Let $\alpha =w_{n-1} g_{n-1}\ldots g_it_i^{k} \in {\rm
Y}_{d,n}(u)$, where  $1\leq i\leq n-1$ and $w_{n-1}\in {\rm Y}_{d,n-1}(u)$. Then we have:
$$
{\rm tr}\left(\alpha \, e_{d,n} \right) =  z\, {\rm tr}(\alpha^{\prime} e_{d, n-1})
$$
where $\alpha^{\prime} := g_{n-2}\ldots g_it_i^{k}w_{n-1} \in {\rm Y}_{d,n-1}(u)$.
\end{lem}

\begin{proof}
We have:
\begin{eqnarray*}
{\rm tr}(\alpha e_{d,n}) & = & \frac{1}{d}\sum_s {\rm tr}(w_{n-1} g_{n-1}\ldots g_it_i^{k}t_n^s t_{n+1}^{d-s})\\
 & = & \frac{1}{d}\sum_sx_{d-s}{\rm tr}(w_{n-1} g_{n-1}\ldots g_it_i^{k}t_n^s)\\
& = & \frac{1}{d}\sum_sx_{d-s}{\rm tr}(w_{n-1}t_{n-1}^s g_{n-1}\ldots g_it_i^{k})\\
& = & \frac{z}{d}\sum_sx_{d-s}{\rm tr}(w_{n-1}t_{n-1}^s g_{n-2}\ldots g_it_i^{k})\\
& = & \frac{z}{d}\sum_sx_{d-s}{\rm tr}(g_{n-2}\ldots g_it_i^{k}w_{n-1}t_{n-1}^s ) \\
& = & \frac{z}{d}\sum_sx_{d-s}{\rm tr}(\alpha^{\prime}t_{n-1}^s )
 =  \frac{z}{d}\sum_s{\rm tr}(\alpha^{\prime}t_{n-1}^s t_n^{d-s})
  =  z\, {\rm tr}(\alpha^{\prime} e_{d, n-1}).
\end{eqnarray*}
\end{proof}

\noindent$\bullet$ For $n=3$ we have  for $\alpha$ the following possibilities:
$$
\begin{array}{clllll}
{\rm (i)} & t_1^{h}t_2^{l}t_3^{k} &
{\rm (ii)} & t_1^{h}t_2^{l}g_2t_2^k &
{\rm (iii)} & t_1^{h}t_2^{l}g_2g_1t_1^k \\
{\rm (iv)}& t_1^h g_1t_1^{l}t_3^{k} &
{\rm (v)} & t_1^h g_1t_1^{l}g_2t_2^k &
 {\rm (vi)} & t_1^h g_1t_1^{l}g_2g_1t_1^k
\end{array}
$$
Cases (i) and (iv) are applications of Lemma \ref{aone}.
 Cases (ii), (iii) and (v) show also factorizing through ${\rm tr}(\alpha)$. Indeed, by direct computations we obtain:
$$
{\rm (ii)}  \quad  {\rm tr} \left(\alpha e_{d,3}\right) =
\frac{E_d^{({k + l})}}{x_{k + l}} \, {\rm tr}(\alpha) \ \ \
$$
$$
{\rm (iii)} \quad  {\rm tr} \left(\alpha  e_{d,3}\right) =
\frac{E_d^{({h + k +  l})}} {x_{h + k + l}}\, {\rm tr}(\alpha )
$$
$$
{\rm (v)}  \quad  {\rm tr} \left(\alpha  e_{d,3}\right) =
\frac{E_d^{({h + k + l})}} {x_{h + k + l}}\, {\rm tr}(\alpha) \ \ \
$$
Notice that, even in the above simple cases where ${\rm tr}(\alpha e_{d,3})$ factors through ${\rm tr}(\alpha)$, the factors are not the same and they are different from $E_d$. It remains now to consider case (vi) for $\alpha$. Indeed we have:

\begin{eqnarray*}
{\rm (vi)} \ \  {\rm tr}(\alpha) & = &  z\, {\rm tr}(g_1^2t_2^{l}t_1^{h+k}) \\
&\stackrel{(\ref{quadr})}{=}  & z\, x_{l}x_{h+k} + (u-1)z E_d^{({h+k+l})} + (u-1)z^2
x_{h+k+l}, \ \ \ \ \ \ \ \ \ \ \ \ \ \ \ \ \ \ \ \
\end{eqnarray*}
while:
\begin{eqnarray*}
{\rm tr} (\alpha e_{d,3}) & = &  \frac{z}{d}\sum_{s}x_{d-s}{\rm tr}(g_1^2t_2^{l}t_1^{h+k+s}) \\
&\stackrel{(\ref{quadr})}{=}  &  z\, x_{l}E_d^{(h+k)} + \frac{(u-1)z}{d}\sum_{s}x_{d-s}E_d^{(h+k+l + s)} + (u-1)z^2 E_d^{(h+k+l)}.
\end{eqnarray*}

It is clear from the above that in order to have Eq.~\ref{mul} we must impose conditions on the set of indeterminates ${X}_d$ (recall (\ref{xd})). For example we have the following:

\begin{lem}\label{cyclic}
Let $c\in {\Bbb C}\setminus \{0\}$.
Setting $x_i = c^i$, we have:
$$
{\rm tr}(e_{d,n})=1, \quad {\rm tr}(e_{d,n}^{(k)})=c^k \quad {\rm \text and } \quad {\rm tr}(\alpha e_{d,n})  = {\rm tr}(\alpha) {\rm tr}(e_{d,n})
\qquad (\alpha\in {\rm Y}_{d,n}(u))
$$
\end{lem}
\begin{proof} The first two equalities follow from (\ref{Ed}) and (\ref{treditok}) by a direct computation. We shall prove the last one by induction. For $n=1$ we have $\alpha = t_1^k$. So, ${\rm tr}(\alpha) = c^k$
 \ and \ ${\rm tr}(\alpha e_{d,1}) = E_d^{(k)} = \frac{1}{d}\sum_{s=0}^{d-1}c^{k + s}
c^{d-s} =c^k$. Suppose the statement is true for any element in ${\rm Y}_{d,n-1}(u)$ and let $\alpha$ be an element of the inductive basis of ${\rm Y}_{d,n}(u)$.

If $\alpha = w_{n-1}t_n^k$ with $w_{n-1}\in {\rm Y}_{d,n-1}(u)$, then by Lemma~\ref{aone} we have: ${\rm tr}(\alpha e_{d,n}) = \frac{E_{d}^{(k)}}{x_k} \, {\rm tr}\left(\alpha\right) = \frac{c^k}{c^k} \, {\rm tr}\left(\alpha\right)$.

If $\alpha = g_{n-1}\ldots g_it_i^{k}w_{n-1}$ with $w_{n-1}\in {\rm Y}_{d,n-1}(u)$, then by Lemma~\ref{atwo} we have:
$ {\rm tr}(\alpha e_{d,n}) = z\, {\rm tr}(\alpha^{\prime} e_{d, n-1})$, where $\alpha^{\prime} = g_{n-2}\ldots g_it_i^{k}w_{n-1} \in {\rm Y}_{d,n-1}(u)$.
Using now the induction hypothesis on the word $\alpha^{\prime}$ we obtain:
${\rm tr}(\alpha e_{d,n}) = z\, {\rm tr}(\alpha^{\prime}) {\rm tr}(e_{d, n-1})=  {\rm tr}(g_{n-1}\alpha^{\prime}) {\rm tr}(e_{d, n})  = {\rm tr}(\alpha) {\rm tr}(e_{d, n})$.
\end{proof}

\begin{rem} \rm Unfortunately, the condition $x_i = c^i$ does not lead to an interesting framed link invariant from the topological viewpoint. For example:
$$
{\rm tr}_d (t_1^k t_2^l) = x_k x_l = c^{k+l} = x_{l+k} = {\rm tr}_d (t_1^{k+l} t_2^0)
$$
 but the closures of these two 2--stranded braids are not isotopic as framed (un)links of two components.
\end{rem}

\subsection{\it The $E$--system}

We shall now seek conditions on a set ${X}_d$ of $d-1$ non--zero complex number, other than those of Lemma~\ref{cyclic}, so that (\ref{mul}) is satisfied.

\begin{defn}\label{defcon}\rm
For $m= 0, \ldots ,d-1$, let ${\rm E}_{d}^{(m)}$ denote the polynomial
\begin{equation}\label{Epoly}
{\rm E}_{d}^{(m)} =\sum_{s=0}^{d-1}{\rm x}_{m+s}{\rm x}_{d-s}
\end{equation}
where, by definition ${\rm x}_0 = {\rm x}_d = 1$, and the sub--indices are regarded module $d$.
We say that the set of complex numbers ${X}_d =\{x_1, \ldots ,x_{d-1}\}$ satisfies the {\it $E$--condition} if
$x_1,\ldots , x_{d-1}$ satisfy the following {\it $E$--system} of non--linear equations in ${\Bbb C}$:
\begin{eqnarray*}\label{Esystem}
{\rm E}_d^{(1)} & = &  {\rm x}_1 {\rm E}_d^{(0)}\\
{\rm E}_d^{(2)} & = &  {\rm x}_2 {\rm E}_d^{(0)}\\
& \vdots & \\
{\rm E}_d^{(d-1)} & = &  {\rm x}_{d-1}{\rm E}_d^{(0)}
\end{eqnarray*}
Equivalently:
\begin{equation}\label{Esystem}
\sum_{s=0}^{d-1}{\rm x}_{m + s} {\rm x}_{d-s}  =  {\rm x}_m \sum_{s=0}^{d-1} {\rm x}_{s} {\rm x}_{d-s} \qquad (1\leq m \leq d-1)
\end{equation}
\end{defn}

Note that if ${X}_d$ satisfies the $E$--system then:
\begin{equation}\label{Econdition}
\frac{E_d^{(m)}}{x_m} = E_{d}= {\rm tr}(e_{d,i}) \qquad (1\leq m \leq d-1)
\end{equation}
Clearly, the $E$--condition guarantees a common factor, namely $E_d={\rm tr}(e_{d,n})$, at least for the cases of $\alpha$ where ${\rm tr}(\alpha\, e_{d,n})$ factors through ${\rm tr}(\alpha)$ (recall the case $n=3$). Surprisingly, we also have the following result.

\begin{thm}\label{thme}
If ${X}_d$ satisfies the $E$--condition then for all $\alpha \in
{\rm Y}_{d,n}(u)$ we have:
$$
{\rm tr}(\alpha e_{d,n}) = {\rm tr}(\alpha)\, {\rm
tr}( e_{d,n})
$$
\end{thm}

\begin{proof}
By the linearity of the trace it suffices to consider the case when  $\alpha$ is an element in the inductive basis (\ref{basis}) of ${\rm Y}_{d,n}(u)$. We proceed by induction on $n$. For $n=1$ we have: ${\rm tr}(\alpha e_{d,1}) = \frac{E_d^{(k)}}{x_k} {\rm tr}(\alpha) =  E_d\, {\rm tr}(\alpha) = {\rm tr}(\alpha)\, {\rm tr}( e_{d,1})$.
Suppose the statement is true for $n-1$, that is, for all elements in ${\rm Y}_{d,n-1}(u)$, and let   $\alpha$ be an element of the inductive basis of ${\rm Y}_{d,n}(u)$.

 If $\alpha = w_{n-1}t_n^k$ with $w_{n-1}\in {\rm Y}_{d,n-1}(u)$, then by Lemma~\ref{aone} we have:
$$
{\rm tr}(\alpha e_{d,n}) = \frac{E_{d}^{(k)}}{x_k} \, {\rm tr}\left(\alpha\right)
= \, E_d {\rm tr}\left(\alpha\right) = {\rm tr}(\alpha){\rm tr}(e_{d,n}).
$$

If  $\alpha = g_{n-1}\ldots g_it_i^{k}w_{n-1}\in {\rm
Y}_{d,n}(u)$ with $w_{n-1}\in {\rm Y}_{d,n-1}(u)$, then by Lemma~\ref{atwo} we have:
$ {\rm tr}(\alpha e_{d,n}) = z\, {\rm tr}(\alpha^{\prime} e_{d, n-1})$, where $\alpha^{\prime} = g_{n-2}\ldots g_it_i^{k}w_{n-1}$ in ${\rm Y}_{d,n-1}(u)$.
Using now the induction hypothesis on the word $\alpha^{\prime}$ we obtain:
 ${\rm tr}(\alpha e_{d,n}) = z\, {\rm tr}(\alpha^{\prime}) {\rm tr}(e_{d, n-1})=  {\rm tr}(g_{n-1}\alpha^{\prime}) {\rm tr}(e_{d, n})  = {\rm tr}(\alpha) {\rm tr}(e_{d, n})$.
\end{proof}

Next, we give a useful computational result using the $E$--condition.

\begin{lem}\label{cond}
For the set of indeterminates ${X}_d$ we have, assuming the $E$--condition:
$$
\frac{x_k}{d}\sum_{s}x_{d-s}E_d^{(k+s)} = \left[E_d^{(k)}\right]^2
\qquad ( k\in {\Bbb N}).
$$
Equivalently,
$$
x_k \, {\rm tr}\left(e_{d, n+1} e_{d,n}^{(k)}\right) = \left[{\rm
tr}\left(e_{d,n}^{(k)}\right)\right]^2.
$$
\end{lem}
\begin{proof}
Indeed:
$$
\frac{1}{d}\sum_{s}x_{d-s}E_d^{(k+s)} = \frac{1}{d}\sum_{s}x_{d-s}x_{k+s} E_d =
 E_d E_d^{(k)} = x_k^{-1} \left[ E_d^{(k)}\right]^2.
$$
\end{proof}

Let us see how exactly the $E$--condition works in the case (vi) of $n=3$, namely when   $\alpha= t_1^hg_1t_1^{l}g_2g_1t_1^k$. Recall:
\begin{center}
$
 {\rm tr}(\alpha) =
 z\, x_{l}x_{h+k} + (u-1)z \,E_d^{({h+l+k})} + (u-1)z^2 \,x_{h+l+k}.
$
\end{center}
Hence:
\begin{center}
$
{\rm tr}(\alpha){\rm tr}( e_{d,3}) =  z\, x_{l}x_{h+k} E_d + (u-1)z \,E_d^{({h+ l +k})} E_d + (u-1)z^2\, x_{h+ l +k}E_d.
$
\end{center}
On the other hand:
$$
{\rm tr} (\alpha e_{d,3})  =
z\, x_{l}E_d^{(h+k)} + \frac{(u-1)z}{d}\sum_{s}x_{d-s}E_d^{(h+ l +k+ s)} +
(u-1)z^2 \,E_d^{(h+ l +k)}.
$$
Then, from Lemma \ref{cond}:
$$
 {\rm tr} (\alpha e_{d,3})  =
 z\, x_{l}E_d^{(h+k)} +  (u-1)z\, \frac{ \left[E_d^{(h+ l+k)} \right]^2}{x_{h+ l+k}}\quad + (u-1)z^2\, E_d^{(h+ l+k)}.
$$
Applying now the $E$--condition to ${X}_d$ yields immediately: ${\rm tr} (\alpha e_{d,3}) = {\rm tr}(\alpha){\rm tr}( e_{d,3})$.

\subsection{\it Solutions of the $E$--system}

 The $E$--system  has always a not--all--zero solution. For example, we have the cyclic solution:
$$
{\rm x}_k= \zeta^k,
$$
where $\zeta$ is a primitive $d$th root of unity. Indeed:

\vspace{.07in}
${\rm E}_d^{(m)} = \sum_{s=0}^{d-1} {\rm x}_{m + s} {\rm x}_{d-s} = \zeta^m \quad \mbox{and} \quad  {\rm E}_d  = 1$, \quad so \quad ${\rm x}_m {\rm E}_d = {\rm E}_d^{(m)}$.
\vspace{.07in}

 The solution ${\rm x}_k= \zeta^k$ of the $E$--system is a special case of Lemma~\ref{cyclic}, so it is not interesting for our topological purposes.

\begin{rem} \rm
 It is worth observing at this point that the values ${\rm x}_i = c^i$ of Lemma~\ref{cyclic} do not comprise, in general, a solution of the $E$--system. For example, for $d=3$ we have the $E$--system:
$$
\begin{array}{lcr}
{\rm x}_1 + {\rm x}_2^2& = & 2 {\rm x}_1^2 {\rm x}_2 \\
{\rm x}_1^2 + {\rm x}_2& = & 2 {\rm x}_1 {\rm x}_2^2
\end{array}
$$
Substituting now ${\rm x}_i = c^i$ does not automatically satisfy the system equations.
\end{rem}

Beyond the above cyclic solution, for $d=3,4$ and $5$ we run the Mathematica program and we found  other  solutions of the $E$--system, for which:
$$
E_d = {\rm tr}(e_{d,i}) \neq 1, \,\,\text{for all} \,\,i.
$$
 For example,
in the case $d=3$ we have the non--trivial solutions:
$$
{\rm x}_1={\rm x}_2 = -\frac{1}{2}\quad \text{and}\quad
{\rm x}_1=\frac{1}{2}\left(-1 + i\sqrt{3} \right), \
{\rm x}_2=\frac{1}{4}\left( 1 + i\sqrt{3}\right),
$$
and also the solution where we take the conjugates in the previous one.

\smallbreak
 Consider now the set $\delta = \{\delta_j\}$ of the real numbers:
\begin{equation}\label{deltai}
\delta_j : = \frac{-(-1)^{j(d-1)}}{d-1}\qquad (j = 1,\ldots ,  d-1)
\end{equation}
and denote ${\rm E}_{d}^{(m)}(\delta)$ the evaluation of ${\rm E}_{d}^{(m)}$ at ${\rm x}_j = \delta_j$. According to (\ref{Epoly}) we have: ${\rm E}_{d}^{(0)} = 1 + \sum_{s=1}^{d-1}{\rm x}_{s}{\rm x}_{d-s}$. Then:
$$
{\rm E}_d^{(0)}(\delta)  = 1 + \sum_{s=1}^{d-1} \frac{(-1)^{s(d-1)}(-1)^{(d-s)(d-1)}}{(d-1)^2}
=
1 + \sum_{s=1}^{d-1} \frac{(-1)^{d(d-1)}}{(d-1)^2} = 1 + \frac{(d-1)}{(d-1)^2}.
$$
Hence:
\begin{center}
$
{\rm E}_{d}^{(0)}(\delta)  = \frac{d}{d-1}.
$
\end{center}

\begin{prop}
The set $\delta = \{\delta_j\}$ of Eq.~\ref{deltai} is a solution of the $E$--system (\ref{Esystem}). Moreover, for $d\not=2$ we have for this solution: $E_d = {\rm tr}(e_{d,i}) \neq 1$.
\end{prop}

\begin{proof}
 For $m= 1, \ldots , d-1$ we have:
$$
{\rm E}_{d}^{(m)} = 2{\rm x}_m + \sum_{s\not= 0, d-m}{\rm x}_{m+s}{\rm x}_{d-s}
$$
So:
\begin{eqnarray*}
{\rm E}_{d}^{(m)} (\delta)
&  = &
\frac{-2(-1)^{m(d-1)}}{d-1} + \frac{1}{(d-1)^2}\sum_{s\not= 0,d-m}(-1)^{(d+m)(d-1)}\\
& = &
\frac{-2(-1)^{m(d-1)}}{d-1} + \frac{(-1)^{(d-1)m}}{(d-1)^2}\sum_{s\not= 0,d-m}(-1)^{d(d-1)}\\
& = &
\frac{-2(-1)^{m(d-1)}}{d-1} + \frac{(-1)^{(d-1)m}}{(d-1)^2}(d-2)\\
& = &
\frac{-2(-1)^{m(d-1)}(d-1) + (-1)^{(d-1)m}(d-2)}{(d-1)^2}\\
& = &
\frac{-d(-1)^{m(d-1)}}{(d-1)^2} = \delta_m {\rm E}_{d}^{(0)}(\delta).
\end{eqnarray*}

Moreover, for $d\not=2$ we have:
\begin{center}
$E_d={\rm tr}(e_{d,i}) = \frac{1}{d}{\rm E}_{d}^{(0)}(\delta) = \frac{1}{d-1} \neq 1$.
\end{center}
\end{proof}

In the Appendix to this paper we give the general solution of the $E$--system, due to Paul G\'{e}rardin. Namely, for $a\in {\Bbb Z}/d{\Bbb Z}$ we denote  ${\rm exp}_a$ the exponential character of the group ${\Bbb Z}/d{\Bbb Z}$, that is:
$$
{\rm exp}_a (k): = \cos \frac{2\pi  a k}{d} + i \sin \frac{2\pi  a k}{d}\qquad (k\in {\Bbb Z}/d{\Bbb Z}).
$$
Then, the solutions of the $E$--system are parametrized by the non--empty subsets $S$ of ${\Bbb Z}/d{\Bbb Z}$. More precisely, a non--empty subset $S$ defines the solution $X_{d, S}=\{x_0, x_1, \ldots ,x_{d-1}\}$,  where:
$$
x_k = \frac{1}{\vert S\vert}\sum_{s\in S}{\rm exp}_s(k)\qquad (0\leq k\leq d-1).
$$
As we note in \cite{jula3} it is always $x_0=1$. Also,
\begin{equation}\label{valedi}
E_d = {\rm tr}\left( e_{d,i}\right) = \frac{ 1}{\vert S\vert}\qquad (1\leq i\leq n-1).
\end{equation}

\subsection{\it Lifting solutions to the $p$--adic level}

Let $r$ and $s$ be two positive integers, such that $r\geq s$. We shall prove now that a  solution of
 the $E$--system for $d= p^s$ lifts to a solution of the corresponding $E$--system for $d^{\prime}= p^r$. This is important for showing that there are also interesting solutions at the $p$--adic level.

Given $\zeta =(\zeta_1, \ldots ,\zeta_d)\in {\Bbb C}^d$ we define
$\zeta^{\prime} =(\zeta_1^{\prime}, \ldots
,\zeta_{d^{\prime}}^{\prime})\in {\Bbb C}^{d^{\prime}}$ as follows:
$$
\zeta_j^{\prime}= \left\{
\begin{array}{ll}
\zeta_i, & \text{for}\, i= 1, \ldots , d-1\\
\zeta_j, & \text
{for}\,i\equiv j \quad ({\rm mod}\, d).
\end{array}\right.
$$
Then we have the following commutative diagram:
\begin{equation}
\begin{diagram}\label{rsev}
\node{R\left[{\mathfrak X}_{r}\right]}
\arrow{e,t}{\delta_s^r} \arrow{se,b}{{\rm
ev}_{\zeta{\prime}}} \node{R\left[{\mathfrak X}_s\right]}
\arrow{s,r}{{\rm ev}_{\zeta}}\\
\node{} \node{R}
\end{diagram}
\end{equation}
where ${\rm ev}_c$ is the  evaluation  homomorphism  at $c\in {\Bbb C}^m$.

\begin{prop}\label{lifsol}
Let $d=p^s$ and  $d^{\prime}= p^r$ with $r\geq s$. If $\zeta$ is a solution of the system of equations
${\rm E}_d^{(k)}= {\rm x}_k {\rm E}_d \ (k=1,\ldots,d-1)$ then $\zeta^{\prime}$ is a
solution of the system of equations ${\rm E}_{d^{\prime}}^{(k)}=
{\rm x}_k {\rm E}_{d^{\prime}} \ (k=1,\ldots,d^{\prime}-1)$.
\end{prop}

\begin{proof}
The equation ${\rm E}_d^{(k)}= {\rm x}_k {\rm E}_d$ can be written  as:
${\rm tr}_d\left( e_{d,n}^{(k)}\right) = {\rm tr}_d\left(e_{d,n}t_{n+2}^k\right)$.
Now:  ${\rm ev}_{\zeta^{\prime}}\left({\rm
tr}_{d^{\prime}}\left(e_{d^{\prime},n}^{(k)}\right)\right)
 =
\left( {\rm ev}_{\zeta^{\prime}}\circ {\rm
 tr}_{d^{\prime}}\right)\left(e_{d^{\prime},n}^{(k)}\right)$,
and from  diagram (\ref{rsev}) we have:
\begin{eqnarray*}
{\rm ev}_{\zeta^{\prime}}\left({\rm
tr}_{d^{\prime}}\left(e_{d^{\prime},n}^{(k)}\right)\right)
 & = &
\left( {\rm ev}_{\zeta}\circ \delta_s^r\circ {\rm
 tr}_{d^{\prime}}\right)\left(e_{d^{\prime},n}^{(k)}\right)\\
 & = &
\left( {\rm ev}_{\zeta}\circ {\rm tr}_d \circ
 \phi_d^{d^{\prime}}\right)\left(e_{d^{\prime},n}^{(k)}\right) \quad (\text{Lemma \ref{commute}})\\
  & = &
{\rm ev}_{\zeta}\left( {\rm tr}_d \left(
 \phi_d^{d^{\prime}}\left(e_{d^{\prime},n}^{(k)}\right)\right)\right)\\
 & = &
{\rm ev}_{\zeta}\left( {\rm tr}_d \left(e_{d,n}^{(k)}\right)\right)\\
 & = &
{\rm ev}_{\zeta}\left( {\rm tr}_d \left(e_{d,n}t_{n+2}^k\right)\right)\quad (\text{Induction hypothesis})\\
 & = &
{\rm ev}_{\zeta}\left( {\rm
 tr}_d \left( \phi_d^{d^{\prime}}\left( e_{d^{\prime},n}t_{n+2}^k\right)\right)\right) \\
 & = &
{\rm ev}_{\zeta}\left( \delta^{r}_s  \left( {\rm
 tr}_{d^{\prime}}\left( e_{d^{\prime},n}t_{n+2}^k\right)\right)\right)\quad (\text{Lemma \ref{commute}})\\
 & = &
{\rm ev}_{\zeta^{\prime}} \left( {\rm
 tr}_{d^{\prime}}\left( e_{d^{\prime},n}t_{n+2}^k\right)\right)\quad (\text{Diagram (\ref{rsev})}).
\end{eqnarray*}
Hence $\zeta^{\prime}$ is a solution of the system:
$
{\rm E}_{d^{\prime}}^{(k)}= {\rm x}_k {\rm E}_{d^{\prime}}\quad (1\leq k\leq d^{\prime}-1).
$
\end{proof}

\section{Isotopy invariants of framed and $p$--adic framed links}

In this section we define an infinite family of isotopy
invariants of (modular) oriented framed links using the Markov traces of Theorem~\ref{trace} and the Markov equivalence of Theorem~\ref{markov}. Then, we also define an isotopy
invariant of classical oriented framed links and of $p$--adic oriented framed links using the $p$--adic Markov trace of Theorem~\ref{ptrace} and the $p$--adic Markov equivalence of Theorem~\ref{markov}.

\subsection{\it Framed link invariants from ${\rm tr}_d$}

Let $d\in {\Bbb N}$ and let ${X}_{d,S} =\{x_1, \ldots ,x_{d-1}\}$ be a solution of the  $E$--system parametrized by a non--empty  subset $S$ of ${\Bbb Z}/d{\Bbb Z}$. Then, using  Theorem \ref{thme} and Lemma~\ref{engn}(i) we can proceed with the factorization of ${\rm tr}_d(\alpha e_{d,n})$ in (\ref{traginverse}). Indeed, we have from (\ref{traginverse}):
\begin{eqnarray*}\label{Etraginverse}
{\rm tr}_d(\alpha g_n^{-1}) & = &
\left[z + (u^{-1}-1)E_d + (u^{-1}-1)z) \right] {\rm tr}_d(\alpha) \\
& = &  \frac{z - (u-1)E_d}{u}\, {\rm tr}_d(\alpha) = {\rm tr}_d(g_n^{-1}) \, {\rm tr}_d(\alpha) \qquad \qquad
\end{eqnarray*}

\noindent where $E_d$ was defined in (\ref{Ed}). For the value of $E_d$ under the $E$--condition, recall from (\ref{valedi}) that $E_d=\frac{1}{|S|}$. In analogy to the construction of the Jones polynomial\cite{jo}, we will first do  a re--scaling, by which $\alpha g_n$ and $\alpha g_n^{-1}$ will be assigned the same trace value for any $\alpha \in {\rm Y}_{d, n}(u)$.  More precisely, we define
\begin{equation}\label{omega}
\omega := \frac{z - (u-1)E_d}{uz}
\end{equation}
 so
\begin{equation}\label{Etraginverse2}
{\rm tr}_d(g_n^{-1}) = \frac{z - (u-1)E_d}{u} =\omega z.
\end{equation}
Then the re-scaling map:
$$
\sigma_i\mapsto \sqrt{\omega}g_i, \quad t_j^{s {(\rm mod}\, d)} \mapsto t_j^{s {(\rm mod}\, d)}
$$
defines the representation:
$$
\Omega_{d,\omega}: {\mathcal  F}_{d,n} \longrightarrow {\rm Y}_{d,n}(u)
$$
Moreover, composing with the natural projection from ${\mathcal  F}_{n}$ to ${\mathcal  F}_{d,n}$, the representation $\Omega_{d,\omega}$ lifts to a representation of ${\mathcal  F}_{n}$ to ${\rm Y}_{d,n}(u)$;  we retain for this the same notation, $\Omega_{d,\omega}$. Then we have the following:

\begin{defn}\label{gamma} \rm
 Given a solution of the $E$--system parametrized by a non--empty  subset $S$ of ${\Bbb Z}/d{\Bbb Z}$, for any framed braid $\alpha \in {\mathcal  F}_{n}$ we define for its closure $\widehat{\alpha}$:
$$
\Gamma_{d,S}(\widehat{\alpha}) := \left(\frac{1-\omega u}{\sqrt{\omega}(u-1)E_d}\right)^{n-1}
\left({\rm tr}_d \circ \Omega_{d,\omega}\right)(\alpha)
$$
Defining further the {\it exponential sum} $\epsilon(\alpha)$ of $\alpha$ as the algebraic sum of the exponents of the $\sigma_i$'s in $\alpha$ and denoting
$$
\Delta := \frac{1-\omega u}{\sqrt{\omega}(u-1)E_d} = \frac{1}{z\sqrt{\omega}}
$$
we can write:
$$
\Gamma_{d,S}(\widehat{\alpha}) = \Delta^{n-1} (\sqrt{\omega})^{\epsilon(\alpha)} {\rm tr}_d (\alpha)
$$
where the $\alpha$ in ${\rm tr}_d (\alpha)$ is the image of the braid $\alpha$ in ${\rm Y}_{d,n}(u)$ under the natural mapping: $\sigma_i \mapsto g_i$ and $t_j^s  \mapsto t_j^{s( {\rm mod}\, d)}$.
\end{defn}

Let now ${\mathcal L}$ denote the set of oriented  framed links and let ${\Bbb C}(z, x_1, \ldots , x_{d-1})$ be, as usual, the ring of rational functions on $z, {X}_{d,S}$ with complex coefficients. Then we have the following.

\begin{thm}\label{gammainv}
If the set ${X}_{d,S}$ satisfies the $E$--condition then the
map $\Gamma_{d,S}$ is an isotopy invariant of (modular) oriented  framed links:
$$
\begin{array}{ccccl}
\Gamma_{d,S} & : & {\mathcal L} & \longrightarrow &  {\Bbb C}(z, x_1, \ldots , x_{d-1})\\
& & \widehat{\alpha} & \mapsto & \Gamma_{d,S}(\widehat{\alpha})
\end{array}
$$
\end{thm}

\begin{proof}
By the classical Alexander theorem, any link can be isotoped to the closure of some braid. Then, by Theorem~\ref{markov}, in order to show that $\Gamma_{d,S}$ is constant on the isotopy class of the oriented framed link $\widehat{\alpha}$ for any $\alpha\in \cup_n {\mathcal F}_{n}$, it suffices to show that $\Gamma_{d,S}(\widehat{\alpha}) = \Gamma_{d,S}(\widehat{\alpha\sigma_n}) = \Gamma_{d,S}(\widehat{\alpha\sigma_n^{-1}})$ for  $\alpha \in {\mathcal F}_n$. The first equality is taken care by the coefficient $\Delta^{n-1}$ in Definition~\ref{gamma} and the second one by the re--scaling of the trace. More precisely, and since $\epsilon(\alpha \sigma_n) = \epsilon(\alpha )+1$ and $\epsilon(\alpha \sigma_n^{-1}) = \epsilon(\alpha )-1$, we obtain:
$$
\Gamma_{d,S}(\widehat{\alpha\sigma_n}) = \Delta^n \sqrt{\omega}^{\epsilon(\alpha g_n)}
{\rm tr}_d(\alpha g_n) = \Delta \sqrt{\omega} \, z \,\Gamma_{d,S} (\widehat{\alpha}),
$$
$$
\Gamma_{d,S}(\widehat{\alpha\sigma_n^{-1}}) = \Delta^n \sqrt{\omega}^{\epsilon(\alpha g_n^{-1})}
{\rm tr}_d(\alpha g_n^{-1}) = \Delta \sqrt{\omega}\, z \,\Gamma_{d,S} (\widehat{\alpha})
$$
Therefore, and since $\Delta \sqrt{\omega}\, z  = 1$, the proof of the Theorem is concluded.
\end{proof}

\subsection{\it Some computations}

In all computations that follow the integers $m$ in $x_m$ and $E_d^{(m)}$ are considered ${\rm mod}\, d$. Let ${X}_{d,S} =\{x_1, \ldots, x_{d-1}\}$ be a solution of the $E$--system parametrized by a non--empty subset $S$ of ${\Bbb Z}/d{\Bbb Z}$.

\smallbreak
\noindent $\bullet$ Clearly,  for the unknot ${\rm O}$ with framing zero we have $\Gamma_{d,S} ({\rm O})=1$. For the framed unknot ${\rm O}^k$ with framing $k \in {\Bbb Z}$ we have $\Gamma_{d,S} ({\rm O}^k)=x_{k}$.

\smallbreak
\noindent $\bullet$ Let ${\rm H}=\widehat{\sigma_1^2 t_1^k t_2^l}$ be the Hopf link with framings $k,l\in {\Bbb Z}$. We have $\epsilon(\sigma_1^2 t_1^k t_2^l) = 2$ and, using Remark~\ref{edikl}, we compute:
$$
\Gamma_{d,S}({\rm H}) = \Delta \omega\, {\rm tr}_d (g_1^2 t_1^k t_2^l) = \Delta \omega \left[x_lx_k +(u-1)E_d^{(k+l)} + (u-1)z\, x_{k+l} \right].
$$

\smallbreak
\noindent $\bullet$ Let ${\rm T}=\widehat{\sigma_1^3 t_1^k}$ be the right-handed trefoil with framing $k \in {\Bbb Z}$.  We have $\epsilon(\sigma_1^3 t_1^k)= 3$ and, using Lemma~\ref{powers} and Remark~\ref{edikl}, we compute:
$$
\Gamma_{d,S}({\rm T}) = \Delta \sqrt{\omega}^3 {\rm tr}_d (g_1^3 t_1^k) = \Delta \sqrt{\omega}^3 \left[(u^2 -u +1)z\, x_k + (u^2 -u) E_d^{(k)} \right].
$$

\smallbreak
\noindent $\bullet$ Let ${\rm T^{\prime}}=\widehat{\sigma_1^{-3} t_1^k}$ be the left--handed trefoil with framing $k \in {\Bbb Z}$.  We have $\epsilon(\sigma_1^{-3} t_1^k)= -3$ and, using Lemma~\ref{powers} and Remark~\ref{edikl}, we compute:
$$
\Gamma_{d,S}({\rm T^{\prime}}) = \Delta (\sqrt{\omega})^{-3} \left[(u^{-3} -u^{-2} +u^{-1})z\, x_k
+ (u^{-3} -u^{-2} +u^{-1}-1)E_d^{(k)} \right].
$$

\subsection{\it A skein relation for $\Gamma_{d,S}$}

Let $L_{+}$, $L_{-}$, $L_s$ and $L_{s\times}$, $s=0,\ldots,d-1$, be diagrams of oriented framed links, which are all identical, except near one crossing, where they differ by the ways indicated in Figure~\ref{skein}. Then we have:

\begin{prop}\label{skein}
The invariant $\Gamma_{d,S}$ satisfies the following skein relation:
$$
 \sqrt{\omega} \,\Gamma_{d,S} (L_{-})
= \frac{1}{\sqrt{\omega}} \, \Gamma_{d,S} (L_{+})  + \frac{u^ {-1}-1}{d} \sum_{s=0}^{d-1}\Gamma_{d,S} (L_s)
+ \frac{u^ {-1}-1}{d\sqrt{\omega}} \sum_{s=0}^{d-1}\Gamma_{d,S} (L_{s\times})
$$
\end{prop}

The above linear skein relation derives from (\ref{invrs}) and is diagrammatically related to Figure~\ref{g1invrs}, but with different coefficients.

\smallbreak

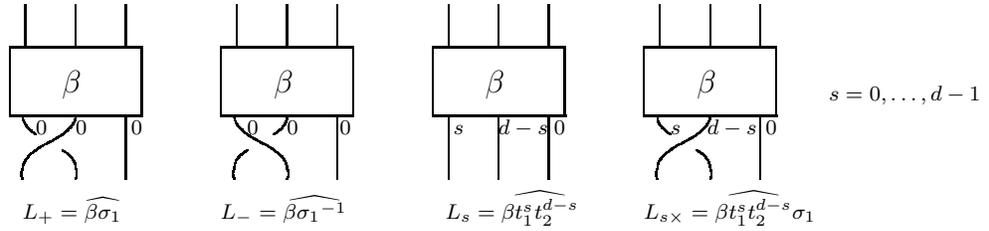
\begin{figure}[H]
 
\begin{picture}(330,80)

\put(20,53){$\beta$}

\qbezier(6,70)(6,78)(6,86) 
\qbezier(25,70)(25,78)(25,86)
\qbezier(44,70)(44,78)(44,86)

\qbezier(5,20)(4,22)(5,24) 
\qbezier(25,20)(26,22)(25,24) 

\qbezier(0,44)(0,59)(0,70)
\qbezier(50,44)(50,59)(50,70)

\qbezier(0,70)(25,70)(50,70)
\qbezier(0,44)(25,44)(50,44)
\qbezier(5,24)(5,29)(15,34)
\qbezier(15,34)(25,39)(25,44)
\qbezier(25,24)(25,28)(20,31)
\qbezier(10,37)(5,40)(5,44)
\qbezier(44,20)(44,32)(44,44)


\put(100,53){$\beta$}

\qbezier(86,70)(86,78)(86,86) 
\qbezier(105,70)(105,78)(105,86)
\qbezier(124,70)(124,78)(124,86)

\qbezier(85,20)(84,22)(85,24) 
\qbezier(105,20)(106,22)(105,24) 

\qbezier(80,44)(80,59)(80,70)
\qbezier(130,44)(130,59)(130,70)

\qbezier(80,70)(105,70)(130,70)
\qbezier(80,44)(105,44)(130,44)

\qbezier(85,24)(85,28)(90,31)
\qbezier(100,37)(105,40)(105,44)
\qbezier(105,24)(105,29)(95,34)
\qbezier(95,34)(85,39)(85,44)

\qbezier(124,20)(124,32)(124,44)


\put(180,53){$\beta$}

\qbezier(166,70)(166,78)(166,86) 
\qbezier(185,70)(185,78)(185,86)
\qbezier(204,70)(204,78)(204,86)
\qbezier(160,44)(160,59)(160,70)
\qbezier(210,44)(210,59)(210,70)

\qbezier(160,70)(185,70)(210,70)
\qbezier(160,44)(215,44)(210,44)

\qbezier(166,20)(166,32)(166,44) 
\qbezier(185,20)(185,32)(185,44)
\qbezier(204,20)(204,32)(204,44)

\put(260,53){$\beta$}

\qbezier(246,70)(246,78)(246,86) 
\qbezier(265,70)(265,78)(265,86)
\qbezier(284,70)(284,78)(284,86)

\qbezier(245,20)(244,22)(245,24) 
\qbezier(265,20)(266,22)(265,24) 

\qbezier(240,44)(240,59)(240,70)
\qbezier(290,44)(290,59)(290,70)

\qbezier(240,70)(265,70)(290,70)
\qbezier(240,44)(295,44)(290,44)
\qbezier(245,24)(245,29)(255,34)
\qbezier(255,34)(265,39)(265,44)
\qbezier(265,24)(265,28)(260,31)
\qbezier(250,37)(245,40)(245,44)
\qbezier(284,20)(284,32)(284,44)
\put(10,37){\tiny{$0$}}
\put(25,37){\tiny{$0$}}
\put(46,37){\tiny{$0$}}

\put(90,37){\tiny{$0$}}
\put(105,37){\tiny{$0$}}
\put(125,37){\tiny{$0$}}

\put(168,37){\tiny{$s$}}
\put(185,37){\tiny{$d-s$}}
\put(206,37){\tiny{$0$}}

\put(250,37){\tiny{$s$}}
\put(264,37){\tiny{$d-s$}}
\put(286,37){\tiny{$0$}}

\put(310,50){\tiny{$s=0, \ldots ,d-1$}}
\put(5,5){\tiny{$L_{+}=\widehat{\beta \sigma_1}$}}
\put(80,5){\tiny{$L_{-}=\widehat{\beta {\sigma_1}^{-1}}$}}
\put(165,5){\tiny{$L_s=\widehat{\beta t_1^st_2^{d-s}}$}}
\put(240,5){\tiny{$L_{s\times}=\widehat{\beta t_1^st_2^{d-s}}\sigma_1$}}

\end{picture}
\caption{The links in the skein relation}\label{skein}
\end{figure}

\begin{proof}
The proof is standard. By the Alexander theorem for framed links we may assume that $L_+$ is in braided form and that $L_+ = \widehat{\beta \sigma_i}$ for some $\beta\in {\mathcal  F}_n$ and for some $i$.  Then:
\begin{center}
$
L_- = \widehat{\beta \sigma_i^{-1}}, \qquad L_s = \widehat{\beta t_i^s t_{i+1}^{d-s}}, \qquad L_{s\times} = \widehat{\beta t_i^s t_{i+1}^{d-s} \sigma_i}  \quad (s=0,\ldots, d-1).
$
\end{center}
We now apply relation (\ref{invrs}) for the $g_i^{-1}$ in the expression for $\Gamma_{d,S} (L_{-})$, we recall (\ref{edi}) and we note that $\epsilon(\beta \sigma_i^{-1}) = \epsilon(\beta)-1$, $\epsilon(\beta \sigma_i) = \epsilon(\beta)+1$, $\epsilon(\beta t_{i}^s t_{i+1}^{d-s}) = \epsilon(\beta)$ and $\epsilon(\beta t_{i}^s t_{i+1}^{d-s} \sigma_i) = \epsilon(\beta)+1$, to obtain:
\begin{eqnarray*}
\Gamma_{d,S} (L_{-})
&  = &
\Delta^{n-1} (\sqrt{\omega})^{\epsilon(\beta \sigma_i^{-1})} {\rm tr}_d (\beta g_i^{-1}) \\
& \stackrel{(\ref{invrs})}{=} &
\Delta^{n-1} (\sqrt{\omega})^{\epsilon(\beta)-1}
\left[ {\rm tr}_d (\beta g_i) + (u^{-1} -1) {\rm tr}_d (\beta e_{d,i}) + (u^{-1} -1) {\rm tr}_d (\beta e_{d,i}g_i) \right] \\
& = &
\Delta^{n-1} (\sqrt{\omega})^{\epsilon(\beta)}
\left[ \frac{1}{\sqrt{\omega}} {\rm tr}_d (\beta g_i) +
\frac{u^{-1} -1}{\sqrt{\omega}} {\rm tr}_d (\beta e_{d,i}) + \frac{u^{-1} -1}{\sqrt{\omega}} {\rm tr}_d (\beta e_{d,i} g_i) \right] \\
& \stackrel{(\ref{edi})}{=} &
\frac{1}{\omega} \Delta^{n-1} (\sqrt{\omega})^{\epsilon(\beta)+1} {\rm tr}_d (\beta g_i)
 +
\frac{u^{-1} -1}{d\sqrt{\omega}} \sum_{s=0}^{d-1} \Delta^{n-1} (\sqrt{\omega})^{\epsilon(\beta)} {\rm tr}_d (\beta t_i^s t_{i+1}^{d-s}) \\
& + &
\frac{u^{-1} -1}{d\omega} \sum_{s=0}^{d-1} \Delta^{n-1} (\sqrt{\omega})^{\epsilon(\beta)+1}  {\rm tr}_d (\beta t_i^s t_{i+1}^{d-s} g_i) \\
& = &
\frac{1}{\omega} \Gamma_{d,S} (L_{+})
 +
\frac{u^{-1} -1}{d\sqrt{\omega}} \sum_{s=0}^{d-1} \Gamma_{d,S} (L_s)
+
\frac{u^{-1} -1}{d\omega} \sum_{s=0}^{d-1} \Gamma_{d,S} (L_{s\times}).
\end{eqnarray*}
Equivalently, multiplying by $\sqrt{\omega}$ we finally obtain:
$$
 \sqrt{\omega} \,\Gamma_{d,S} (L_{-})
= \frac{1}{\sqrt{\omega}} \, \Gamma_{d,S} (L_{+})  + \frac{u^ {-1}-1}{d} \sum_{s=0}^{d-1}\Gamma_{d,S} (L_s)
+ \frac{u^ {-1}-1}{d\sqrt{\omega}} \sum_{s=0}^{d-1}\Gamma_{d,S} (L_{s\times}).
$$
\end{proof}

\begin{rem} \rm
Note that, by (\ref{edikrels}) we have $\beta e_{d,i} g_i = \beta g_i e_{d,i} $. Thus:
\begin{center}
$
{\rm tr}_d (\beta e_{d,i} g_i) = {\rm tr}_d (\beta g_i e_{d,i} ) =  {\rm tr}_d (e_{d,i} \beta g_i ) =
\frac{1}{d} \sum_{s=0}^{d-1} {\rm tr}_d (t_i^s t_{i+1}^{d-s} \beta g_i ).
$
\end{center}

Hence $L_{s\times}$ could be also considered to be the link $\widehat{t_i^s t_{i+1}^{d-s} \beta  g_i}$. Also, $L_s = \widehat{\beta t_i^s t_{i+1}^{d-s}} = \widehat{t_i^s t_{i+1}^{d-s} \beta }$. These observations help see the extra framings in $L_{s\times} $ and $L_{s}$, that come from $t_i^s t_{i+1}^{d-s}$, at the top of the corresponding braid.
\end{rem}

\begin{rem} \rm
If we restrict our interest to classical braids and we consider them as framed braids with all framings zero we could also `kill' the framing on the level of the algebra $Y_{d,n}(u)$ by taking $d=1$. Then, as we mentioned earlier in the paper (Remarks~\ref{YtoH} and~\ref{ocneanu}), the algebra $Y_{1,n}(u)$ coincides with the classical Iwahori--Hecke algebra and the trace ${\rm tr}_1$ coincides with the Ocneanu trace. Hence, the invariant $\Gamma_{1}$ coincides with the HOMFLYPT polynomial and the skein relation boils down to the well--known skein relation of the HOMFLYPT polynomial.
\end{rem}

\begin{rem} \rm
As mentioned in the Introduction, in \cite{jula3} we  represented the classical braid group $B_n$ in the Yokonuma--Hecke algebra $Y_{d,n}(u)$ by treating the framing generators just as formal elements. So, $\Gamma_{d,S}$ can be seen as an invariant of classical knots. But, then, a skein relation has no topological interpretation.  As we showed in \cite{jula3}, in that case $\Gamma_{d,S}$ satisfies a `closed' cubic relation, which factors to the quadratic relation of the Iwahori--Hecke algebra ${\rm H}_n(u)$. Further, in \cite{ChLa} we show that for generic values of the parameters $u,z$ the invariants $\Gamma_{d,S}$ do not coincide with the HOMFLYPT polynomial. Yet, our computational data \cite{CJKL} seem to indicate that these invariants do not distinguish more or less knot pairs than the HOMFLYPT polynomial. Also, in \cite{CJKL} we show that the invariants $\Gamma_{d,S}$ (but not the traces ${\rm tr}_d$) have the multiplicative property on connected sums.

Even if, eventually, we just obtain invariants topologically equivalent to the HOMFLYPT polynomial, the Yokonuma--Hecke algebras comprise the only known examples of algebras  with topological applications to different knot categories, which support Markov traces that do not re--scale directly according to topological braid equivalence.
\end{rem}

\subsection{\it A link invariant for classical and $p$--adic framed links}

In this subsection we construct an invariant for classical oriented framed links (with no modular restriction on the framing) and for $p$--adic oriented framed links, by lifting the values of the invariants $\Gamma_{d,S}$ to the $p$--adic context.

Let ${\mathcal L}_{p^{\infty}}$ denote the set of $p$--adic oriented framed links. For positive integers $r, s$ such that $r\geq s$, the connecting ring epimorphism $\delta_s^r$ (recall (\ref{delta})) yields a connecting epimorphism $\Xi_s^r$ from the ring of rational functions ${\Bbb C}(z,{\mathfrak X}_{r})$ to the ring of rational functions ${\Bbb C}(z,{\mathfrak X}_s)$. It is a routine to prove  the following lemma.

\begin{lem}\label{gammalift}
For all $r\geq s\geq v$, the following diagram is commutative.
$$
\begin{diagram}
\node{\cdots}
\node{{\mathcal L}_{p^{\infty}}}
\arrow{w,t}{} \arrow{s,l}{\Gamma_{p^v}} \node{{\mathcal L}_{p^{\infty}}}\arrow{w,t}{{\rm Id}}
\arrow{s,l}{\Gamma_{p^s}}
\node{{\mathcal L}_{p^{\infty}}}
\arrow{w,t}{{\rm Id}} \arrow{s,l}{\Gamma_{p^r}} \node{\cdots}\arrow{w,t}{}
\\
\node{\cdots}
\node{{\Bbb C}\left(z, {\mathfrak  X}_v\right)}
\arrow{w,t}{}
\node{{\Bbb C}\left(z, {\mathfrak  X}_s\right)}\arrow{w,t}{\Xi_v^s}
\node{{\Bbb C}\left(z, {\mathfrak  X}_r\right)}
\arrow{w,t}{\delta_s^r}
\node{\cdots}\arrow{w,t}{}
\end{diagram}
$$
\end{lem}
The ring $\varprojlim R[{\mathfrak X}_r]$  turns out to be an integral domain. We shall also denote $R_{p^{\infty}}$ the field of fractions of $\varprojlim R[{\mathfrak X}_r]$. Taking now inverse limits in the diagram of Lemma~\ref{gammalift} we obtain the map:
\begin{center}
$
\Gamma_{p^{\infty}} := \varprojlim_{r\in{\Bbb N}}\Gamma_{p^r}
$
\end{center}

\begin{thm}\label{padicinv}
If for all $r$ the set ${\mathfrak X}_{r}$ satisfies the $E$--condition, then the map
$$
\begin{array}{ccccl}
\Gamma_{p^{\infty}} & : &  {\mathcal L}_{p^{\infty}} & \longrightarrow &R_{p^{\infty}} \\
& & \widehat{\underleftarrow{\alpha}} & \mapsto & (\Gamma_{p}(\widehat{\alpha_1}), \Gamma_{p^2}(\widehat{\alpha_2}),\ldots)
\end{array}
$$
for any ${\underleftarrow{\alpha}}= (\alpha_r)_r \in \cup_n{\mathcal F}_{p^{\infty}, n}$ is constant on the equivalence classes defined by the $p$--adic version of Theorem~\ref{markov}.
In particular, restricting to the set ${\mathcal L}$ of classical oriented framed links, we have that the map
$$
\begin{array}{ccccl}
\Gamma_{p^{\infty}} & : &  {\mathcal L} & \longrightarrow &R_{p^{\infty}} \\
& & \widehat{\alpha} & \mapsto & (\Gamma_{p}(\widehat{\alpha}), \Gamma_{p^2}(\widehat{\alpha}),\ldots)
\end{array}
$$
for any $\alpha \in \cup_n{\mathcal F}_n$ is an invariant of oriented framed links.
\end{thm}

\begin{proof}
By  Proposition~\ref{lifsol} we have non-trivial solutions of the $E$--system in the $p$--adic context.
Let now  ${\underleftarrow{\beta}}= (\beta_r)_r, \, {\underleftarrow{\alpha}}= (\alpha_r)_r  \in \cup_n{\mathcal F}_{p^{\infty}, n}$ be  Markov equivalent $p$--adic framed braids. Then, according to Theorem~\ref{markov}, we have that in each position the modular framed braids $\beta_r$ and $\alpha_r$ are Markov equivalent. So, $\Gamma_{p^r}(\widehat{\beta_r}) = \Gamma_{p^r}(\widehat{\alpha_r})$, hence $\Gamma_{p^{\infty}}(\widehat{\underleftarrow{\alpha}}) = \Gamma_{p^{\infty}}(\widehat{\underleftarrow{\beta}})$.

Moreover, restricting $\Gamma_{p^{\infty}}$ to the set ${\mathcal L}$ of classical oriented framed links we have that $\Gamma_{p^{\infty}}$ is also an isotopy invariant of classical oriented framed links. Recall that the classical framed braid group is represented in the $p$--adic Yokonuma--Hecke algebra $Y_{p^{\infty},n}$ without modular restriction on the framing and this fact is preserved in the invariant $\Gamma_{p^{\infty}}$.
\end{proof}

\section{Appendix: The $E$--system \\ by Paul G\'{e}rardin}

We interpret each polynomial in  (4.1)
$$
\sum_{0\le s<d}  x_{m+s}x_{d-s}, \ 0\le m <d
$$
in the $d$ complex numbers $x_0,x_1,...,x_{d-1}$ as the value at $m$ of the convolution product by itself of the element $x : s\mapsto x_s$ in the complex algebra $\Bbb C[D]$ of the cyclic group $D=\Bbb Z/d\Bbb Z$ : the convolution product $f\ast g$ of two elements $f,g\in  \Bbb C[D] $ is
$$f\ast g : w\mapsto \sum_{u+v=w}f(u)g(v),$$
the sum being on the set of $(u,v)\in D\times D$ with sum $w$.\par

The algebra $\Bbb C[D]$ is commutative algebra with unit $\delta_0$, the characteristic function of the unit element $0\in D$. It is the direct sum of its simple ideals $\Bbb Ce_a, a\in D$, the $e_a$'s being the characters of the group $D$ : $$e_a : u\mapsto e^{2\pi iau/d}.$$

They satisfy the following relations : $e_a\ast e_b$ is $de_a$ for $a=b$ and $0$ otherwise, so that   the $e_a/d,a\in D$ are its elementary idempotents. \smallbreak

The algebra $\Bbb C[D]$ has another product, with unit $e_0$, given by the product of values:
$$f g : w\mapsto f(w)g(w),$$
and is the direct sum of its simple ideals $\Bbb C\delta_a, a\in D$, where  $\delta_a$ is the characteristic function of the   element $a\in D$ ; they are also the elementary idempotents for this structure.
\smallbreak

The Fourier transform on $\Bbb C[D]$:
$$
f\mapsto \widehat f : v\mapsto (f\ast e_v)(0) = \sum_{u\in D}f(u)e_v(-u)
$$
is a linear automorphism.  In particular, $\widehat{\delta_a}=e_{-a}, \widehat{e_a}=d\delta _a, a\in D$. Its inverse is $f\mapsto\big(u\mapsto{1\over d}\widehat {f}(-u)\big),$
which means \ $\displaystyle\widehat{\!\!\widehat {f}}(u)=d\,f(-u)$\par
The Fourier transform sends the convolution product to the product of values:
$$
\widehat{f\ast g}=\widehat{f}\;\widehat{g},
$$
hence $\widehat{f g}=d^{-1}\widehat{f}\ast\widehat{g}$.
\medbreak

The $E$--condition (4.2) can now be written as:
$$x\ast x=(x\ast x)(0)\,x.$$

To solve the $E$--system $x(0)=1, x\ast x=(x\ast x)(0)\,x$, we use Fourier transform to obtain:
$$
x(0)=1, \ \widehat{x}\, ^2=(x\ast x)(0)\, \widehat{x}
$$
If $(x\ast x)(0)=0$, then  $\widehat{x}\, ^2=0$ so $\widehat{x}$ is $0$ and also is $x$, which is excluded by the condition $x(0)=1$. Now, the equation says that the function $ \widehat{x}$  is constant on its support $S$ where it is  $(x\ast x)(0)$.  As the characteristic function  of $S$ is, up to the factor $d$,  the Fourier transform of the sum of $e_a, a\in S$, we have shown that
$$
\displaystyle x=(x\ast x)(0){1\over d}\sum_{s\in S}e_s
$$
As $x(0)=1$, we have $\displaystyle (x\ast x)(0){1\over d}|S|=1$, with $|S|$ the cardinality of $S$.

Finally, we have proved that the solutions of the  $E$--system are the functions $x_S$ parametrized by the non--empty subsets $S$ of the cyclic group $D$ of order $d$ as follows:
$$
x_S={1\over |S|}\sum_{s\in S}e_s
$$
For $S=D$, it is the trivial solution $\delta_0$. The complement of the support of any non trivial solution is another solution. In particular, each element $a\in D$ defines two solutions of the $E$--system : one is the character $e_a$, the other is given by $\displaystyle {e_a\over 1-d}$ outside $0$. When the order $d$ is even, we can take $a=d/2$, this gives the solution $\displaystyle u\mapsto {(-1)^u\over 1-d}, u\neq 0$.



\end{document}